\documentclass[11pt,thmsa]{article}
\usepackage{amsfonts}

\usepackage{graphicx}
\usepackage{amsmath}

\newtheorem{theorem}{Theorem}[section]

\newtheorem{corollary}[theorem]{Corollary}

\newtheorem{definition}[theorem]{Definition}

\newtheorem{lemma}[theorem]{Lemma}

\newtheorem{proposition}[theorem]{Proposition}
\newtheorem{remark}[theorem]{Remark}

\newenvironment{proof}[1][Proof]{\textbf{#1.} }{\ \rule{0.5em}{0.5em}}
\input{tcilatex}

\begin{document}

\title{Wigner measures in the discrete setting: high-frequency analysis of sampling
\& reconstruction operators}
\author{Fabricio Maci\`{a} \\
DMA - \'{E}cole Normale Sup\'{e}rieure,\\
45, rue d'Ulm,\\
75230 Paris cedex 05,\\
France.\\
email: \texttt{fabricio.macia@ens.fr}}
\date{}
\maketitle

\begin{abstract}
The goal of this article is that of understanding how the oscillation and
concentration effects developed by a sequence of functions in $\mathbb{R}^{d}
$ are modified by the action of Sampling and Reconstruction operators on
regular grids. Our analysis is performed in terms of Wigner and defect
measures, which provide a quantitative description of the high frequency
behavior of bounded sequences in $L^{2}\left( \mathbb{R}^{d}\right) $. We
actually present explicit formulas that make possible to compute such
measures for sampled/reconstructed sequences. As a consequence, we are able
to characterize sampling and reconstruction operators that preserve or
filter the high-frequency behavior of specific classes of sequences. The
proofs of our results rely on the construction and manipulation of Wigner
measures associated to sequences of discrete functions.


\end{abstract}

\tableofcontents

\section{Introduction}

\subsection[Statement of the problem]{Statement of the problem: oscillation and concentration under
the effect of sampling and reconstruction}

A central problem in Numerical Analysis and Signal Theory is that of
reconstructing a function $u\left( x\right) $ defined in $\mathbb{R}^{d}$
from a discrete set of measurements taken on an uniform grid of step size $h$%
. This discrete values are typically obtained by applying to the function $u$
a \textbf{sampling operator }$S_{\varphi }^{h}$ of the following type:
\begin{equation*}
S_{\varphi }^{h}u\left( n\right) :=\frac{1}{h^{d}}\int_{\mathbb{R}%
^{d}}u\left( x\right) \overline{\varphi \left( \frac{x}{h}-n\right) }dx,
\end{equation*}
for some \textbf{sampling function }$\varphi $. One then tries to recover $u$
by means of a \textbf{reconstruction} (or \textbf{interpolation})\textbf{\
operator} $T_{\psi }^{h}$ through the formula:
\begin{equation}
T_{\psi }^{h}S_{\varphi }^{h}u\left( x\right) :=\sum_{n\in \mathbb{Z}%
^{d}}S_{\varphi }^{h}u\left( n\right) \psi \left( \frac{x}{h}-n\right) ,
\label{recon}
\end{equation}
where $\psi $ is some fixed \textbf{reconstruction function}. This process
usually only provides an approximation of the original function $u$, with an
error that vanishes as $h$ tends to zero. Such reconstruction schemes have
been the object of intensive study both from the point of view of
Approximation Theory and Numerical Analysis.

Here we shall be concerned with the \textbf{high-frequency} approximation
properties of those operators, that is, we shall study how a reconstruction
scheme such as (\ref{recon}) is able to capture (or filter) \textbf{%
oscillation} and \textbf{concentration}-like phenomena on the functions it
is intended to approximate. More generally, me shall be interested in
clarifying how the high frequency behavior of a sequence of reconstructed
functions depends on the profiles $\varphi$, $\psi$ and the sampling rate $h$
chosen.

Before giving a more precise statement of our objectives, let us first
illustrate the above discussion with two specific examples: consider $%
f_{k}\left( x\right) :=k^{d/2}\rho \left( k\left( x-x_{0}\right) \right) $
and $g_{k}\left( x\right) :=\rho \left( x\right) e^{ikx\cdot \xi ^{0}}$ with
$\rho \in L^{2}\left( \mathbb{R}^{d}\right) $; the sequence $\left(
f_{k}\right) $ concentrates around the point $x_{0}$ as $k\rightarrow \infty
$, whereas $\left( g_{k}\right) $ oscillates in the direction $\xi ^{0}$.
The results we shall present in this paper are aimed to understand to what
extent the sequences $\left( T_{\psi }^{h_{k}}S_{\varphi
}^{h_{k}}f_{k}\right) $ and $\left( T_{\psi }^{h_{k}}S_{\varphi
}^{h_{k}}g_{k}\right) $ reproduce the same behavior as $\left( f_{k}\right) $
and $\left( g_{k}\right) $ (i.e., if concentration and oscillation persist),
for a given sequence $\left( h_{k}\right) $ of positive reals that tends to
zero (the sampling steps) and some choice of $\varphi $ and $\psi $.

Perhaps, the simplest convenient setting to formulate our results is
provided by the notion of \textbf{defect measure}, an object that gives a
quantitative description of what we shall understand by concentration and
oscillation effects and whose definition we next recall. Let $\left(
u_{k}\right) $ be a weakly converging sequence in the space $L^{2}\left(
\mathbb{R}^{d}\right) $; denote by $u$ its weak limit and remark that the
densities $\left| u_{k}-u\right| ^{2}$ are uniformly bounded in $L^{1}\left(
\mathbb{R}^{d}\right) $. Helly's compactness Theorem then ensures that some
subsequence $\left( \left| u_{k_{n}}-u\right| ^{2}\right) $ weakly converges
in the set of positive Radon measures;\footnote{%
From now on, we shall use the term \emph{measure }as an abbreviation of the
longer \emph{Radon measure}. Recall that the space of Radon measures $%
\mathcal{M}\left( \mathbb{R}^{d}\right) $ is identified, by Riesz's Theorem,
with the space of continuous linear functionals on $C_{c}\left( \mathbb{R}%
^{d}\right) $.} or, in other words, that there exists a positive measure $%
\nu $ on $\mathbb{R}^{d}$ such that
\begin{equation*}
\int_{\mathbb{R}^{d}}\phi\left( x\right) \left| u_{k_{n}}\left( x\right)
-u\left( x\right) \right| ^{2}dx\rightarrow\int_{\mathbb{R}^{d}}\phi\left(
x\right) d\nu\left( x\right) \qquad\text{as }n\rightarrow\infty\text{,}
\end{equation*}
for every $\phi\in C_{c}\left( \mathbb{R}^{d}\right) $. When the above
convergence takes place without extracting a subsequence we say that $\nu$
is the \textbf{defect measure }of the sequence $\left( u_{k}\right) $.

Immediately from this definition one deduces the following general
principle: if $\nu$ is the defect measure of a sequence $\left( u_{k}\right)
$ and $\omega\subset\mathbb{R}^{d}$ is a bounded Borel set, then there is an
equivalence between $\nu\left( \omega\right) =0$ and the fact that $%
u_{k}|_{\omega}$ converges strongly to $u|_{\omega}$ in $L^{2}\left(
\omega\right) $. Thus, the support of $\nu$ is precisely the set where
strong convergence fails, that is, the set where oscillations and
concentrations take place.

But defect measures are also able to detect concentration and oscillatory
phenomena and give quantitative information about them. Consider the
sequences $\left( f_{k}\right) $, $\left( g_{k}\right) $ previously defined;
they both weakly converge to zero in $L^{2}\left( \mathbb{R}^{d}\right) $
and it is easy to check that their respective defect measures are $\left\|
\rho\right\| _{L^{2}\left( \mathbb{R}^{d}\right) }^{2}\delta_{x_{0}}$ and $%
\left| \rho\left( x\right) \right| ^{2}dx$. Notice that, in the first case,
the defect measure actually captures the concentration of the sequence
around the point $x=x_{0}$. In the complementary of that point, where the
sequence converges strongly to zero, the measure vanishes. In the second
example, the defect measure is uniformly distributed on $\mathbb{R}^{d}$,
this being consistent with the fact that strong convergence does not take
place in any subset of $\mathbb{R}^{d}$.

Let us point out that the analysis of concentration and oscillation effects
developed by a sequence of functions is a central issue in many problems of
the Calculus of Variations and Partial Differential Equations. A number of
applications of defect measures may be found in the analysis of variational
problems with loss of compactness performed by P.-L. Lions in \cite{Li85a,
Li85b}.\footnote{%
We also refer to L.C. Evans' notes \cite{Ev} for an exposition of some
additional applications as well as a discussion of other measure-theoretical
objects (such as, for example, Young measures) designed to study the failure
of strong convergence.}

Consider a sequence $\left( u_{k}\right) $, weakly converging to zero in $%
L^{2}\left( \mathbb{R}^{d}\right) $; sample it using a profile $\varphi$ and
form the reconstructed sequence
\begin{equation*}
v_{k}:=T_{\psi}^{h_{k}}S_{\varphi}^{h_{k}}u_{k},
\end{equation*}
for some given $\psi$ and some sequence $\left( h_{k}\right) $ of positive
reals tending to zero. The functions $v_{k}$ are bounded in $L^{2}\left(
\mathbb{R}^{d}\right) $ and tend weakly to zero provided $\varphi$ and $\psi$
satisfy suitable hypotheses (see Lemma \ref{Lemma norm TphiUh} in Section
\ref{Sec samprec} below). Suppose furthermore that the densities $\left|
v_{k}\right| ^{2}$ weakly converge to the defect measure $\nu_{\varphi,\psi}$%
.

One of the main issues addressed in this article is that of understanding
the relations existing between the defect measure $\nu _{\varphi ,\psi }$,
the profiles $\varphi $, $\psi $ and the sequences $\left( u_{k}\right) $, $%
\left( h_{k}\right) $. Among these, we point out:\medskip

\textbf{A. }Is there a formula, valid for any sequence $\left( u_{k}\right) $%
, relating $\nu _{\varphi ,\psi }$ to the defect measure $\nu $ only in
terms of the profiles $\varphi $ and $\psi $?\medskip

\textbf{B. }Given $\left( u_{k}\right) $, characterize the profiles $\varphi
$ and $\psi $ such that $\nu _{\varphi ,\psi }=0$. This is the problem of
\textbf{filtering} since, as we have discussed before, $\nu _{\varphi ,\psi
}=0$ is equivalent to the strong convergence to zero of the sequence $\left(
T_{\psi }^{h_{k}}S_{\varphi }^{h_{k}}u_{k}\right) $.\medskip

\textbf{C. }Similarly, characterize the profiles $\varphi $ and $\psi $ such
that $\nu _{\varphi ,\psi }=\nu $ for a given $\left( u_{k}\right) $.\medskip

\textbf{D. }Finally, characterize the profiles that give $\nu _{\varphi
,\psi }=\nu $ \emph{for every} $\left( u_{k}\right) $. \medskip

We shall prove that the answer to question \textbf{A }is negative. This is
due to the fact that the measure $\nu _{\varphi ,\psi }$ is sensitive to the
characteristic directions of oscillation of the sequence $\left(
u_{k}\right) $, whereas $\nu $ is unable to distinguish them. As we have
seen above, the defect measure of the oscillating sequence $\left(
g_{k}\right) $ equals $\left| \rho \left( x\right) \right| ^{2}dx$
independently of the vector $\xi ^{0}$; that is not the case for $\nu
_{\varphi ,\psi }$. Indeed, under additional assumptions on $\varphi $ and $%
\psi $ we prove (see Theorem \ref{Thm comp Wm} and Corollary \ref{Cor comp
Dm}):
\begin{equation*}
\nu _{\varphi ,\psi }\left( x\right) =\sum_{k\in \mathbb{Z}^{d}}\left|
\widehat{\psi }\left( \xi ^{0}+2\pi k\right) \right| ^{2}\left| \widehat{%
\varphi }\left( \xi ^{0}\right) \right| ^{2}\left| \rho \left( x\right)
\right| ^{2}dx.
\end{equation*}
Thus the measure $\nu _{\varphi ,\psi }$ is $\xi ^{0}$-dependent and cannot
be expressed solely in terms of $\nu $, $\varphi $ and $\psi $. Note that $%
\nu _{\varphi ,\psi }$ is identically zero as soon as any of $\sum_{k\in
\mathbb{Z}^{d}}\left| \widehat{\psi }\left( \xi ^{0}+2\pi k\right) \right|
^{2}$ or $\widehat{\varphi }\left( \xi ^{0}\right) $ is null. Analogously,
the profiles that give $\nu _{\varphi ,\psi }=\nu $ are precisely those
which satisfy
\begin{equation*}
\sum_{k\in \mathbb{Z}^{d}}\left| \widehat{\psi }\left( \xi ^{0}+2\pi
k\right) \right| ^{2}\left| \widehat{\varphi }\left( \xi ^{0}\right) \right|
^{2}=1.
\end{equation*}

Therefore, in order to understand how $\nu_{\varphi,\psi}$ is built, we must
have at our disposal an object that is able to distinguish between
oscillatory phenomena at different directions.

\subsection{Wigner measures}

This refinement is provided by the theory of \textbf{Wigner measures}.%
\footnote{%
This object is present in the work of E.P. Wigner on semiclassical quantum
mechanics \cite{Wig}. Recently, Wigner measures have gained interest since
the works of P. G\'{e}rard \cite{Ge91c}, P.-L. Lions \& Th. Paul \cite
{Li-Pau}, P. Markowich, N. Mauser \& F. Poupaud \cite{M-M-P} among others.
Related objects are the \textbf{Microlocal defect measures }or $\mathbf{H}$-%
\textbf{measures},\textbf{\ }introduced independently by P. G\'{e}rard \cite
{Ge91b} and L. Tartar \cite{Ta90}.} Given a bounded sequence in $L^{2}\left(
\mathbb{R}^{d}\right) $ one associates to it a measure $\mu\left(
x,\xi\right) $ on $\mathbb{R}^{d}\times\mathbb{R}^{d}$ which describes the
concentration and oscillation effects (these are the respective roles of the
variables $x$ and $\xi$) occurring at some characteristic length-scale. This
measure takes into account the characteristic speeds as well as the
directions of propagation of oscillations. One way of defining them consists
in replacing the density $\left| u\left( x\right) \right| ^{2}$ involved in
the definition of the defect measure by the phase space (microlocal)
density:
\begin{equation}
m^{\varepsilon}\left[ u\right] \left( x,\xi\right) :=\frac{1}{\left(
2\pi\varepsilon\right) ^{d}}\overline{u\left( x\right) }\widehat{u}\left(
\xi/\varepsilon\right) e^{ix\cdot\xi/\varepsilon},  \label{meu}
\end{equation}
where $\widehat{u}$ is the Fourier transform of $u$ and $\varepsilon$ is a
positive constant. The $\left( 2\pi\right) ^{-d}$ factor in the definition
of $m^{\varepsilon}\left[ u\right] $ is placed to have:
\begin{equation}
\int_{\mathbb{R}^{d}}m^{\varepsilon}\left[ u\right] \left( x,\xi\right)
d\xi=\left| u\left( x\right) \right| ^{2},\qquad\int_{\mathbb{R}%
^{d}}m^{\varepsilon}\left[ u\right] \left( x,\xi\right) dx=\frac{\left|
\widehat{u}\left( \xi/\varepsilon\right) \right| ^{2}}{\left(
2\pi\varepsilon\right) ^{d}}.  \label{marginals m}
\end{equation}
Thus, the function $m^{\varepsilon}\left[ u\right] $ may be looked at as
joint physical space-Fourier space ``density'', in spite of the fact that $%
m^{\varepsilon}\left[ u\right] $ is not positive in general. However, limits
of these quantities are positive measures:

\begin{theorem}
\label{Thm definition wm}Let $\left( u_{k}\right) $ be a bounded sequence in
$L^{2}\left( \mathbb{R}^{d}\right) $ and let $\left( \varepsilon _{k}\right)
$ be a sequence of positive numbers tending to zero. Then it is possible to
extract a subsequence $\left( u_{k_{n}}\right) $ such that, for every test
function $a\in \mathcal{S}\left( \mathbb{R}^{d}\times \mathbb{R}^{d}\right) $%
,
\begin{equation}
\lim_{n\rightarrow \infty }\int_{\mathbb{R}^{d}\times \mathbb{R}^{d}}a\left(
x,\xi \right) m^{\varepsilon _{k_{n}}}\left[ u_{k_{n}}\right] \left( x,\xi
\right) dxd\xi =\int_{\mathbb{R}^{d}\times \mathbb{R}^{d}}a\left( x,\xi
\right) d\mu \left( x,\xi \right) ,  \label{limit definition wm}
\end{equation}
where $\mu $ is a finite positive measure on $\mathbb{R}^{d}\times \mathbb{R}%
^{d}$.
\end{theorem}

A measure $\mu \in \mathcal{M}_{+}\left( \mathbb{R}^{d}\times \mathbb{R}%
^{d}\right) $ is called the \textbf{Wigner measure} of the sequence $\left(
u_{k}\right) $ at scale $\left( \varepsilon _{k}\right) $ whenever the limit
(\ref{limit definition wm}) holds without extracting a subsequence.
Different proofs of Theorem \ref{Thm definition wm} may be found in \cite
{Ge-Lei, Li-Pau, Ge96}. Let us point out that other quadratic densities may
used to define Wigner measures. For instance, in \cite{Li-Pau} $\mu $ is
obtained by replacing $m^{\varepsilon }\left[ u\right] $ in the limit (\ref
{limit definition wm}), by the more familiar \textbf{Wigner transform}:
\begin{equation}
w^{\varepsilon }\left[ u\right] \left( x,\xi \right) :=\int_{\mathbb{R}%
^{d}}u\left( x-\varepsilon \frac{p}{2}\right) \overline{u\left(
x+\varepsilon \frac{p}{2}\right) }e^{ip\cdot \xi }\frac{dp}{\left( 2\pi
\right) ^{d}}.  \label{definition WT}
\end{equation}
It is also possible to consider \textbf{Wave-packet (Husimi) transforms}. Of
course, all this methods are equivalent (the same limit is obtained), cf.
the discussion in \cite{Ge-Lei}.

The Wigner measure encodes all the information contained in the defect
measure provided the sequence $\left( u_{k}\right) $ oscillates at
frequencies of the order of $\varepsilon_{k}^{-1}$. More precisely (see \cite
{Ge-Lei, Li-Pau}):

\begin{proposition}
\label{Prop e-oscill}If $\mu $ is the Wigner measure at scale $\left(
\varepsilon _{k}\right) $ of a sequence $\left( u_{k}\right) $ and $\nu $ is
the measure obtained as the weak limit in $\mathcal{M}_{+}\left( \mathbb{R}%
^{d}\right) $ of the densities $\left| u_{k}\right| ^{2}dx$, then the
identity
\begin{equation*}
\nu \left( x\right) =\int_{\mathbb{R}^{d}}\mu \left( x,d\xi \right)
\end{equation*}
holds provided $\left( u_{k}\right) $ is $\mathbf{\varepsilon }_{k}$\textbf{%
-oscillatory}:
\begin{equation}
\limsup_{k\rightarrow \infty }\int_{\left| \xi \right| >R/\varepsilon
_{k}}\left| \widehat{u_{k}}\left( \xi \right) \right| ^{2}d\xi \rightarrow
0\qquad \text{as }R\rightarrow \infty .  \label{def e-oscillating}
\end{equation}
\end{proposition}

Notice that condition (\ref{def e-oscillating}) actually expresses that the
energy of the Fourier transform of $u_{k}$ is concentrated in a ball of
radius $R/\varepsilon_{k}$, which should be understood as the requirement
that the sequence $\left( u_{k}\right) $ does not oscillate at length scales
finer than $\varepsilon_{k}$.

To illustrate this discussion it may be helpful to look at explicit
computations. The Wigner measure at scale $\left( \varepsilon _{k}\right) $
of the concentrating sequence $\left( f_{k}\right) $ defined at the
beginning of this section is given by:
\begin{equation}
\mu \left( x,\xi \right) =\left\{
\begin{array}{ll}
\left\| \rho \right\| _{L^{2}\left( \mathbb{R}^{d}\right) }^{2}\delta
_{x_{0}}\left( x\right) \otimes \delta _{0}\left( \xi \right) & \text{if }%
\varepsilon _{k}k\rightarrow 0,\medskip \\
\delta _{x_{0}}\left( x\right) \otimes \left| \widehat{\rho }\left( \xi
\right) \right| ^{2}\dfrac{d\xi }{\left( 2\pi \right) ^{d}} & \text{if }%
\varepsilon _{k}=k^{-1},\medskip \\
0 & \text{if }\varepsilon _{k}k\rightarrow \infty ,
\end{array}
\right.  \label{examp con}
\end{equation}
while for the oscillating sequence $\left( g_{k}\right) $ it can be checked
to be:
\begin{equation}
\mu \left( x,\xi \right) =\left\{
\begin{array}{ll}
\left| \rho \left( x\right) \right| ^{2}dx\otimes \delta _{0}\left( \xi
\right) & \text{if }\varepsilon _{k}k\rightarrow 0,\medskip \\
\left| \rho \left( x\right) \right| ^{2}dx\otimes \delta _{\xi ^{0}}\left(
\xi \right) & \text{if }\varepsilon _{k}=k^{-1},\medskip \\
0 & \text{if }\varepsilon _{k}k\rightarrow \infty .
\end{array}
\right.  \label{examp osc}
\end{equation}
These examples show the importance of the choice of the scale $\left(
\varepsilon _{k}\right) $. When this scale is taken to be coarser than the
characteristic length-scale $k^{-1}$ of oscillation/concentration, it is no
longer true that the projection on the first component of their Wigner
measures coincides with the defect measure. On the other hand, in the case $%
\varepsilon _{k}k\rightarrow 0$ (the scale chosen is much smaller than the
actual oscillation scale) the Wigner measure is not able to capture the
direction of oscillation. Hence, to obtain a complete description, the scale
$\left( \varepsilon _{k}\right) $ must be taken of the same order than that
of the oscillations.

Wigner measures turn out to be the correct tool for comparing the high
frequency behavior of the sequences $\left( u_{k}\right) $ and $\left(
T_{\psi}^{h_{k}}S_{\varphi}^{h_{k}}u_{k}\right) $.

\subsection{Computation of Wigner and defect measures}

Given a sequence of sampling steps $\left( h_{k}\right) $, it is clear that
the functions $T_{\psi}^{h_{k}}S_{\varphi}^{h_{k}}u_{k}$ will not develop
oscillation and concentration effects of characteristic sizes asymptotically
smaller that $h_{k}$. Most commonly, these functions will form an $h_{k}$%
-oscillatory sequence;\footnote{%
However, this may fail for some pathological examples (see paragraph \ref{ce
osc}).} consequently, only Wigner measure at scales coarser or of the same
order than $\left( h_{k}\right) $ will be considered.

In order to establish explicit formulas, we shall require additional
hypotheses on $\varphi $, $\psi $ and on the Wigner measures involved.
Nevertheless, in order to simplify the statement of our results, in this
introduction we shall impose the following (more restrictive) condition on
the admissible profiles:
\begin{equation}
\left| \gamma \left( x\right) \right| \leq C\left( 1+\left| x\right| \right)
^{-d-\varepsilon },\text{\qquad for every }x\in \mathbb{R}^{d}\text{ and
some }C,\varepsilon >0.  \label{strong hyp}
\end{equation}
More general results may be found in Section \ref{Sec Sr}.

We prove the following:

\begin{theorem}
\label{Thm comp Wm}Let $\varphi $, $\psi $ satisfy (\ref{strong hyp}).
Suppose $\left( u_{k}\right) $ is a bounded sequence in $L^{2}\left( \mathbb{%
R}^{d}\right) $ and that $\mu $ is its Wigner measure at scale $\left(
h_{k}\right) $. Suppose moreover that the measures
\begin{equation}
\left| \widehat{\varphi }\left( \xi +2\pi n\right) \right| ^{2}\mu \left(
x,\xi +2\pi n\right)  \label{shan cond}
\end{equation}
are mutually singular for $n\in \mathbb{Z}^{d}$.

Then the Wigner measure at scale $\left( h_{k}\right) $ of the sequence $%
\left( T_{\psi }^{h_{k}}S_{\varphi }^{h_{k}}u_{k}\right) $ is given by:
\begin{equation*}
\mu _{\varphi ,\psi }\left( x,\xi \right) =\left| \widehat{\psi }\left( \xi
\right) \right| ^{2}\sum_{k\in \mathbb{Z}^{d}}\left| \widehat{\varphi }%
\left( \xi +2\pi n\right) \right| ^{2}\mu \left( x,\xi +2\pi n\right) .
\end{equation*}
\end{theorem}

From this, one deduces:

\begin{corollary}
\label{Cor comp Dm}If, moreover, $\left| T_{\psi }^{h_{k}}S_{\varphi
}^{h_{k}}u_{k}\right| ^{2}dx$ weakly converges to a measure $\nu _{\varphi
,\psi }$ then:
\begin{equation*}
\nu _{\varphi ,\psi }\left( x\right) =\int_{\mathbb{R}^{d}}\sum_{k\in
\mathbb{Z}^{d}}\left| \widehat{\psi }\left( \xi +2\pi k\right) \right|
^{2}\left| \widehat{\varphi }\left( \xi \right) \right| ^{2}\mu \left(
x,d\xi \right) .
\end{equation*}
\end{corollary}

This shows, in particular, that a formula relating $\nu _{\varphi ,\psi }$
and the weak limit $\nu $ of $\left| u_{k}\right| ^{2}dx$ does not exist
unless $\left( u_{k}\right) $ is $h_{k}$-oscillatory and $\mu $ is of the
form $\nu \left( x\right) \otimes \sigma \left( \xi \right) $. It also shows
that $\nu =\nu _{\varphi ,\psi }$ if and only if $\sum_{k\in \mathbb{Z}%
^{d}}\left| \widehat{\psi }\left( \xi +2\pi k\right) \right| ^{2}\left|
\widehat{\varphi }\left( \xi \right) \right| ^{2}=1$ for $\mu $-almost every
$\xi \in \mathbb{R}^{d}$. Consequently, there do not exist profiles $\varphi
$, $\psi $ satisfying (\ref{strong hyp}) such that $\nu $ equals $\nu
_{\varphi ,\psi }$ for every $h_{k}$-oscillatory sequence $\left(
u_{k}\right) $.

On the other hand, Theorem \ref{Thm comp Wm} implies that question \textbf{A
}above does have a positive answer in terms of Wigner measures, at least
when restricted to the class of sequences which satisfy (\ref{shan cond}).
That condition, roughly speaking, imposes a restriction on the size of the
region in frequency space where an admissible sequence fails to converge
strongly to zero. Below, we shall compare it with that appearing in
Shannon's sampling Theorem.

The above results will be obtained as corollaries of the more general
Theorems \ref{Thm main} and \ref{Thm defect}. Profiles that belong to
negative-order Sobolev spaces or that fail to satisfy the localization
hypothesis (\ref{strong hyp}) are allowed. However, this will require to
impose compatibility conditions on the Wigner measure $\mu $.

As an illustration of the range of results that will be obtained in this
more general setting, we present an \textbf{asymptotic version }of \textbf{%
Shannon's sampling} \textbf{Theorem}.\footnote{%
See paragraph \ref{sbs examples} for a statement of Shannon's original
sampling Theorem.} It corresponds to taking as sampling profile $\varphi
=\delta _{0}$, the Dirac delta at the origin, and as reconstruction function
$\widehat{\psi }:=\mathbf{1}_{Q}$, where $Q:=\left[ -\pi ,\pi \right) ^{d}$.
Notice that $S_{\delta _{0}}^{h}u\left( n\right) =u\left( hn\right) $ is the
discretization operator, whereas the $T_{\psi }^{h}$ corresponds to
band-limited reconstruction.

\begin{theorem}
\label{Thm Shan}Let $\left( u_{k}\right) $ be a bounded sequence in $%
L^{2}\left( \mathbb{R}^{d}\right) $ and denote by $\mu $ its Wigner measure
at scale $\left( h_{k}\right) $. Suppose, in addition, that $u_{k}\in
H^{s}\left( \mathbb{R}^{d}\right) $ for some $s>d/2$ and
\begin{equation}
\begin{array}{ll}
\text{i)} & \left( 1-h_{k}^{2}\Delta _{x}\right) ^{s/2}u_{k}\text{\quad are
uniformly bounded in }L^{2}\left( \mathbb{R}^{d}\right) .\smallskip \\
\text{ii)} & \mu \left( \mathbb{R}^{d}\times \left( \partial Q+2\pi n\right)
\right) =0\text{\quad for }n\in \mathbb{Z}^{d}.\smallskip \\
\text{iii)} & \mu \left( x,\xi +2\pi n\right) ,\text{\quad }n\in \mathbb{Z}%
^{d},\text{\quad are mutually singular measures.}
\end{array}
\label{hyp Shan}
\end{equation}
Then, the Wigner measure at scale $\left( h_{k}\right) $ of $\left( T_{\psi
}^{h_{k}}S_{\delta _{0}}^{h_{k}}u_{k}\right) $ is
\begin{equation}
\mu _{\delta _{0},\psi }\left( x,\xi \right) =\mathbf{1}_{Q}\left( \xi
\right) \sum_{n\in \mathbb{Z}^{d}}\mu \left( x,\xi +2\pi n\right) .
\label{rep}
\end{equation}
Moreover, if $\left| T_{\psi }^{h_{k}}S_{\delta _{0}}^{h_{k}}u_{k}\right|
^{2}dx$ and $\left| u_{k}\right| ^{2}dx$ weakly converge to $\nu _{S}$ and $%
\nu $, respectively, then
\begin{equation*}
\nu _{S}\left( x\right) =\int_{\mathbb{R}^{d}}\mu \left( x,d\xi \right) =\nu
\left( x\right) .
\end{equation*}
\end{theorem}

Thus, unlike the operators considered in Theorem \ref{Thm comp Wm}, the
composition of discretization and band-limited reconstruction preserves the
defect measure for a large class of sequences.

Notice that, by the Sobolev imbedding Theorem, $S_{\delta _{0}}^{h_{k}}u_{k}$
is well-defined. Actually, (\ref{hyp Shan}.i) ensures that the sequence of
discretizations is square-summable and, consequently, that $\left( T_{\psi
}^{h_{k}}S_{\delta _{0}}^{h_{k}}u_{k}\right) $ is bounded in $L^{2}\left(
\mathbb{R}^{d}\right) $ and $h_{k}$-oscillatory (for a more complete result,
we refer to Lemma \ref{Lemma norm TphiUh}). Condition (\ref{hyp Shan}.ii)
appears because $\widehat{\psi }$ is not continuous; we shall discuss its
necessity in paragraph \ref{sbs nec}. Finally, (\ref{hyp Shan}.iii) should
be understood as the analog of Shannon's original band-limited condition in
this context.

To conclude this short description, let us present how the above results may
be refined when the sequence $\left( u_{k}\right) $ is known to be $%
\varepsilon_{k}$-oscillatory and the sampling rate $\left( h_{k}\right) $ is
taken to satisfy $h_{k}/\varepsilon_{k}\rightarrow0$. As it can be expected,
much more precision is gained:

\begin{theorem}
\label{Thm Wm hlesse}Suppose $\varphi $, $\psi $ satisfy (\ref{strong hyp})
and $\left( u_{k}\right) $ is an $\varepsilon _{k}$-oscillatory, bounded
sequence in $L^{2}\left( \mathbb{R}^{d}\right) $. If $\mu $ is its Wigner
measure at scale $\left( \varepsilon _{k}\right) $ then the corresponding
measure of the sequence $\left( T_{\psi }^{h_{k}}S_{\varphi
}^{h_{k}}u_{k}\right) $ is
\begin{equation*}
\mu _{\varphi ,\psi }=\left| \widehat{\psi }\left( 0\right) \right|
^{2}\left| \widehat{\varphi }\left( 0\right) \right| ^{2}\mu .
\end{equation*}
Moreover, if the densities a $\left| T_{\psi }^{h_{k}}S_{\varphi
}^{h_{k}}u_{k}\right| ^{2}dx$ and $\left| u_{k}\right| ^{2}dx$ weakly
converge to $\nu _{\varphi ,\psi }$ and $\nu $ respectively then
\begin{equation*}
\nu _{\varphi ,\psi }\left( x\right) =\sum_{n\in \mathbb{Z}^{d}}\left|
\widehat{\psi }\left( 2\pi n\right) \right| ^{2}\left| \widehat{\varphi }%
\left( 0\right) \right| ^{2}\nu \left( x\right) .
\end{equation*}
\end{theorem}

This Theorem holds under much more general conditions on $\varphi $ and $%
\psi $ (see Theorem \ref{Thm main2}) and gives a positive answer to question
\textbf{A }provided we consider only $\varepsilon _{k}$-oscillatory
sequences.

An immediate consequence of the above result is that zero-mean sampling
profiles $\varphi $ (i.e. with $\widehat{\varphi }\left( 0\right) =0$, as a
wavelet, for instance) completely filter any oscillations that occur at
scales much coarser than the sampling rate $h_{k}$. For such a profile, $\nu
_{\varphi ,\psi }=0$ for every $\varepsilon _{k}$-oscillatory sequence. An
analogous phenomenon occurs for reconstruction profiles satisfying $\widehat{%
\psi }\left( 2\pi n\right) =0$ for every $n\in \mathbb{Z}^{d}$.

On the other hand, a sufficient condition to have equality between $%
\nu_{\varphi,\psi}$ and $\nu$ is that $\left| \widehat{\varphi}\left(
0\right) \right| =\left| \widehat{\psi}\left( 0\right) \right| =1$ and $%
\left| \widehat{\psi}\left( 2\pi n\right) \right| =0$ for $n\neq0$.

\subsection{Strategy of proof: Wigner measures in the discrete setting}

The proof of the results we have presented above will be achieved by
analyzing separately the sampling and reconstruction operators $S_{\varphi
}^{h}$ and $T_{\psi }^{h}$. In order to develop, it is necessary to deal
with a concept of \textbf{Wigner measure associated to a sequence of
discrete functions}. We shall introduce it by means of a discrete analogous
of the transform $m^{\varepsilon }\left[ \cdot \right] $. We detail this in
the following paragraph.

To a discrete square-summable function $U\in L^{2}\left( h\mathbb{Z}%
^{d}\right) $, where $L^{2}\left( h\mathbb{Z}^{d}\right) $ stands for the
space of the functions $U$ defined on $\mathbb{Z}^{d}$ with values in $%
\mathbb{C}$ such that the norm
\begin{equation*}
\left\| U\right\| _{h}:=\left( h^{d}\sum_{n\in \mathbb{Z}^{d}}\left|
U_{n}\right| ^{2}\right) ^{1/2}
\end{equation*}
is finite, we associate:
\begin{equation}
M^{\varepsilon }\left[ U\right] \left( x,\xi \right) :=\dfrac{h^{2d}}{\left(
2\pi \varepsilon \right) ^{d}}\dsum\limits_{m\in \mathbb{Z}^{d}}\overline{%
U_{m}}\widehat{U}\left( \frac{h}{\varepsilon }\xi \right) e^{im\cdot \left(
h/\varepsilon \right) \xi }\delta _{hm}\left( x\right) .
\label{definition MeUh}
\end{equation}
Here, $\delta _{hm}$ is the Dirac mass centered at the point $hm$ and $%
\widehat{U}$ denotes the discrete Fourier transform:
\begin{equation*}
\widehat{U}\left( \xi \right) :=\sum_{n\in \mathbb{Z}^{d}}U_{n}e^{-in\cdot
\xi },
\end{equation*}
which, as is well known, is a $2\pi \mathbb{Z}^{d}$-periodic function in $L_{%
\text{loc}}^{2}\left( \mathbb{R}^{d}\right) $. The discrete transform $%
M^{\varepsilon }\left[ U\right] $ may be related to the continuous $%
m^{\varepsilon }\left[ u\right] $ by noticing that
\begin{equation}
M^{\varepsilon }\left[ U\right] =m^{\varepsilon }\left[ T_{\delta _{0}}^{h}U%
\right] ,\quad \text{where\quad }T_{\delta _{0}}^{h}U\left( x\right)
=h^{d}\sum_{k\in \mathbb{Z}^{d}}U_{n}^{h}\delta _{hn}\left( x\right) .
\label{relation MUh and mTUh}
\end{equation}
This is meaningful, since $m^{\varepsilon }\left[ u\right] $ is well-defined
for any tempered distribution $u\in \mathcal{S}^{\prime }\left( \mathbb{R}%
^{d}\right) $.

In order to simplify our language we make the following definition:

\begin{definition}
Let $h=\left( h_{k}\right) $ be a scale. We shall call a sequence $\left(
U^{h_{k}}\right) $ $\mathbf{h}_{k}$-\textbf{bounded }if and only if $%
U^{h_{k}}\in L^{2}\left( h_{k}\mathbb{Z}^{d}\right) $ and $\left\|
U^{h_{k}}\right\| _{h_{k}}\leq C$ for every $k\in \mathbb{N}$.
\end{definition}

One has the following convergence result (which is not a direct consequence
of Theorem \ref{Thm definition wm}):

\begin{proposition}
\label{Thm definition wmD}Let $\left( h_{k}\right) $, $\left( \varepsilon
_{k}\right) $ be scales such that $\left( h_{k}/\varepsilon _{k}\right) $ is
bounded\textbf{\ }and let $\left( U^{h_{k}}\right) $ be an $h_{k}$-bounded
sequence of discrete functions. Then $\left( M^{\varepsilon _{k}}\left[
U^{h_{k}}\right] \right) $ is bounded in $\mathcal{S}^{\prime }\left(
\mathbb{R}^{d}\times \mathbb{R}^{d}\right) $ and given any of its convergent
subsequences $\left( U^{h_{k_{n}}}\right) $ there exists a positive measure $%
\mu $ such that,
\begin{equation}
\lim_{n\rightarrow \infty }\left\langle M^{\varepsilon _{k_{n}}}\left[
U^{h_{k_{n}}}\right] ,a\right\rangle _{\mathcal{S}^{\prime }\times \mathcal{S%
}}=\int_{\mathbb{R}^{d}\times \mathbb{R}^{d}}a\left( x,\xi \right) d\mu
\left( x,\xi \right) ,  \label{limit definition wmD}
\end{equation}
for every $a\in \mathcal{S}\left( \mathbb{R}^{d}\times \mathbb{R}^{d}\right)
$.
\end{proposition}

This will be proved as a Corollary of the more general Proposition \ref{Prop
Wm ThU}, which in turn follows from the analysis of Wigner measures in
negative-order Sobolev spaces that is performed in Section \ref{Sec Ap}. As
in the continuous setting, we say that a measure $\mu $ is the \textbf{%
Wigner measure at scale }$\left( \varepsilon _{k}\right) $ of a sequence of
discrete functions $\left( U^{h_{k}}\right) $ if the limit (\ref{limit
definition wmD}) holds for the whole sequence.

\begin{remark}
i) When $\left( h_{k}/\varepsilon _{k}\right) $ is unbounded, it may happen
that $M^{\varepsilon _{k_{n}}}\left[ U^{h_{k_{n}}}\right] $ is not bounded
in $\mathcal{S}^{\prime }\left( \mathbb{R}^{d}\times \mathbb{R}^{d}\right) $%
.\smallskip

ii) If $h_{k}/\varepsilon _{k}\rightarrow c>0$ then $\mu $ is not finite.
Indeed, it is periodic (with respect to the lattice $\left( 2\pi /c\right)
\mathbb{Z}^{d}$) in the $\xi $ variable.\smallskip

iii) However, when $h_{k}/\varepsilon _{k}\rightarrow 0$, the Wigner measure
$\mu $ is finite, as in the continuous case.
\end{remark}

With this tool at our disposal, we are able to compare the Wigner measure of
a sequence of discrete functions $\left( U^{h_{k}}\right) $ with that of a
reconstructed sequence $\left( T_{\psi }^{h_{k}}U^{h_{k}}\right) $.
Analogously, we may compute the Wigner measures of sequences of sampled
discrete functions $\left( S_{\varphi }^{h_{k}}u_{k}\right) $ in terms of
those corresponding to the original sequence $\left( u_{k}\right) $. These
are respectively the contents of Theorems \ref{Thm WM TUh Uh} and \ref{Thm
WM Suh uh}.

\subsection{Plan of the article}

Results and assumptions concerning the operators $S_{\varphi }^{h}$ and $%
T_{\psi }^{h}$ are collected in Section \ref{Sec samprec}.

In Section \ref{Sec h=e}, the problem of computing Wigner measures for
sequences of sampled or reconstructed functions is addressed. Formulas for
Wigner measures at scales of the same order than the sampling/reconstruction
step $\left( h_{k}\right) $ are presented in Theorems \ref{Thm WM TUh Uh}
and \ref{Thm WM Suh uh}. Theorems \ref{Thm comp Wm} and \ref{Thm Shan} then
easily follow from those two results. We also point out the relationships
existing between these Wigner measures and the concept of \textbf{Wigner
series }introduced in \cite{M-M-P, G-M-M-P}.

The problem of the computation of defect measures of sequences of the from $%
\left( T_{\psi }^{h_{k}}U^{h_{k}}\right) $ is considered in Section \ref{Sec
defect}; the main results are presented in Proposition \ref{Proposition
defect} and Corollary \ref{Cor defect}.

In Section \ref{Sec hlesse} we investigate Wigner measures at scales $\left(
\varepsilon _{k}\right) $ satisfying $h_{k}/\varepsilon _{k}\rightarrow 0$.
Explicit formulas are presented in Theorems \ref{Thm WM TUh Uh eh} and \ref
{Thm WM Suh uh eh}, from which Theorem \ref{Thm Wm hlesse} immediately
follows.

The composition of sampling and reconstruction is studied in Section \ref
{Sec Sr}, the main results of this article are proven there.

Finally, Section \ref{Sec Ap} contains the elements from the Theory of
Wigner measures on which the proofs of most of the results of this article
are based on. Propositions \ref{Prop definition wm l2loc} and \ref{Prop
localization}, which extend the Theory of Wigner measures to sequences in
Sobolev spaces of negative order, are systematically used throughout this
paper.

\section{Notations and conventions}

We briefly present some notation that will be used throughout this article.

$B\left( x;R\right) $ will denote the open ball with radius $R$ of $\mathbb{R%
}^{d}$ centered at the point $x$. We shall set
\begin{equation*}
Q:=\left[ -\pi,\pi\right) ^{d},
\end{equation*}
and $\mathbf{1}_{A}$ will denote the characteristic function of a set $%
A\subseteq\mathbb{R}^{d}$.

We write $\Gamma $ to denote de lattice $2\pi \mathbb{Z}^{d}$. A function $f$
defined on $\mathbb{R}^{d}$ is $\Gamma $-periodic if $f\left( x+\gamma
\right) =f\left( x\right) $ for every $\gamma \in \Gamma $ and every $x\in
\mathbb{R}^{d}$.

We adopt the following convention for the Fourier transform:
\begin{equation*}
\widehat{u}\left( \xi\right) :=\int_{\mathbb{R}^{d}}u\left( x\right)
e^{-ix\cdot\xi}dx.
\end{equation*}

Given a measurable function $\varphi \left( \xi \right) $, the\textbf{\
Fourier multiplier }of symbol $\varphi $ is the operator $\varphi \left(
D_{x}\right) $ formally defined by
\begin{equation*}
\varphi \left( D_{x}\right) u\left( x\right) :=\int_{\mathbb{R}^{d}}\varphi
\left( \xi \right) \widehat{u}\left( \xi \right) e^{ix\cdot \xi }\frac{d\xi
}{\left( 2\pi \right) ^{d}}=\check{\varphi}\ast u\left( x\right) ,
\end{equation*}
$\check{\varphi}$ being the inverse Fourier transform of $\varphi $.

A particularly important Fourier multiplier is the \textbf{Bessel potential }%
$\left\langle D_{x}\right\rangle $, of symbol
\begin{equation*}
\left\langle \xi \right\rangle :=\left( 1+\left| \xi \right| ^{2}\right)
^{1/2}.
\end{equation*}
Next, we recall the definition of some function spaces.

As usual, $\mathcal{S}\left( \mathbb{R}^{d}\right) $ denotes the space of
\textbf{rapidly decreasing functions} and $\mathcal{S}^{\prime}\left(
\mathbb{R}^{d}\right) $ stands for its dual, the space of \textbf{tempered
distributions}.

Given $r\in\mathbb{R}$, $H^{r}\left( \mathbb{R}^{d}\right) $, the \textbf{%
Sobolev space }of order $r$, consists of the distributions $u\in\mathcal{S}%
^{\prime}\left( \mathbb{R}^{d}\right) $ such that $\left\langle
D_{x}\right\rangle ^{r}u\in L^{2}\left( \mathbb{R}^{d}\right) $.

The weighted space $L^{2}\left( \mathbb{R}^{d};\left\langle x\right\rangle
^{r}\right) $ is that of the functions $u\in L_{\text{loc}}^{1}\left(
\mathbb{R}^{d}\right) $ such that
\begin{equation*}
\left\| u\right\| _{L^{2}\left( \mathbb{R}^{d};\left\langle x\right\rangle
^{r}\right) }:=\left( \int_{\mathbb{R}^{d}}\left| u\left( x\right) \right|
^{2}\left\langle x\right\rangle ^{r}dx\right) ^{1/2}<\infty .
\end{equation*}
The analogous definition is understood for $L^{\infty }\left( \mathbb{R}%
^{d};\left\langle x\right\rangle ^{r}\right) $.

By $C^{\infty}\left( \mathbb{R}^{d};\left\langle x\right\rangle ^{r}\right) $
we intend the space of functions $u\in C^{\infty}\left( \mathbb{R}%
^{d}\right) $ such that
\begin{equation*}
\left\| \partial_{x}^{\alpha}u\right\| _{L^{\infty}\left( \mathbb{R}%
^{d};\left\langle x\right\rangle ^{r}\right) }<\infty\qquad\text{for every
multiindex }\alpha\in\mathbb{N}^{d}\text{.}
\end{equation*}

$C_{0}\left( \mathbb{R}^{d}\right) $ denotes the spaces of continuous
functions on $\mathbb{R}^{d}$ vanishing at infinity.

Given an open set $\Omega\subseteq\mathbb{R}^{d}$, $\mathcal{M}_{+}\left(
\Omega\right) $ is the set of \textbf{positive Radon measures }on $\Omega$,
which can be identified through Riesz's Theorem to the set of positive
functionals on $C_{c}\left( \Omega\right) $, the space of continuous
functions on $\Omega$ with compact support.

Finally, in order to lighten our writing,

\begin{center}
\emph{we shall write }$\mathcal{S}$\emph{\ and }$\mathcal{S}^{\prime}$\emph{%
\ instead of }$\mathcal{S}\left( \mathbb{R}_{x}^{d}\times \mathbb{R}%
_{\xi}^{d}\right) $\emph{\ and }$\mathcal{S}^{\prime}\left( \mathbb{R}%
_{x}^{d}\times\mathbb{R}_{\xi}^{d}\right) $\emph{\ respectively.}
\end{center}

For a measurable function $f:\mathbb{R}^{d}\rightarrow\mathbb{C}$, we use
the notation
\begin{equation*}
D_{f}:=\left\{ x\in\mathbb{R}^{d}\text{ }:\text{ }f\text{ is not continuous
at }x\right\} .
\end{equation*}

An important, perhaps non-standard, definition is that of a scale:

\begin{definition}
\label{definition scale}A \textbf{scale }$\left( \varepsilon _{k}\right) $
is a sequence of positive numbers that tends to zero as $k\rightarrow \infty
$.
\end{definition}

Given two scales $\left( h_{k}\right) $ and $\left( \varepsilon _{k}\right) $%
, the notations $h_{k}\ll \varepsilon _{k}$ and $h_{k}\sim \varepsilon _{k}$
will be used to indicate that $\lim_{k\rightarrow \infty }h_{k}/\varepsilon
_{k}=0$ and $\lim_{k\rightarrow \infty }h_{k}/\varepsilon _{k}=c>0$
respectively.

\section{\label{Sec samprec}Sampling and reconstruction}

\subsection{\label{sbs examples}Definitions and examples}

The sampling and reconstruction operators we are going to consider are next
described. Given a distribution $\varphi \in \mathcal{S}^{\prime }\left(
\mathbb{R}^{d}\right) $ we set for every $n\in \mathbb{Z}^{d}$ and $h>0$,
\begin{equation*}
\varphi _{n}^{h}\left( x\right) :=\varphi \left( \frac{x}{h}-n\right) .
\end{equation*}
The \textbf{reconstruction }(or \textbf{synthesis})\textbf{\ operator} $%
T_{\varphi }^{h}$, acting on discrete functions $U$ of $\mathbb{Z}^{d}$ is
defined to be
\begin{equation}
T_{\varphi }^{h}U\left( x\right) :=\sum_{n\in \mathbb{Z}^{d}}U_{n}\varphi
_{n}^{h}\left( x\right) .  \label{definition TphiUh}
\end{equation}
This expression is well-defined for finitely supported discrete functions.
When $\varphi $ is a continuous function such that $\varphi \left( 0\right)
=1$ and $\varphi \left( k\right) =0$ for $k\in \mathbb{Z}^{d}\setminus
\left\{ 0\right\} $ then $T_{\varphi }^{h}U$ is actually a function that
\textbf{interpolates} the discrete values $U_{n}$ on the grid $h\mathbb{Z}%
^{d}$, i.e. $T_{\varphi }^{h}U\left( hn\right) =U_{n}$ for all $n\in \mathbb{%
Z}^{d}$.

Analogously, the \textbf{sampling }(or \textbf{analysis})\textbf{\ operator }%
$S_{\varphi }^{h}$, a priori only acting on functions $u\in \mathcal{S}%
\left( \mathbb{R}^{d}\right) $, is defined as follows: $S_{\varphi }^{h}u$
is the discrete function given by
\begin{equation*}
S_{\varphi }^{h}u\left( n\right) :=h^{-d}\left\langle \overline{\varphi
_{n}^{h}},u\right\rangle _{\mathcal{S}^{\prime }\left( \mathbb{R}^{d}\right)
\times \mathcal{S}\left( \mathbb{R}^{d}\right) }.
\end{equation*}
When $\varphi =\delta _{0}$, we obtain the usual \textbf{discretization}
operator: $S_{\delta _{0}}^{h}u\left( n\right) =h^{-d}u\left( hn\right) $
for every $n\in \mathbb{Z}^{d}$.

Indeed, these sampling/reconstruction schemes include several well-known
such procedures on regular grids. Among many others we may cite:\bigskip

$\bullet $ \textbf{Cardinal }$\mathbf{B}$-\textbf{Splines. }The $B$-spline
of order zero is the function $\varphi \left( x\right) :=\mathbf{1}_{\left[
-1/2,1/2\right] ^{d}}\left( x\right) $; the function $T_{\varphi }^{h}U$ is
just the piecewise constant interpolation of the discrete function $U$ on
the grid $h\mathbb{Z}^{d}$. The $B$-spline of order $1$,
\begin{equation*}
\varphi \left( x\right) =\mathbf{1}_{\left[ -1/2,1/2\right] ^{d}}\ast
\mathbf{1}_{\left[ -1/2,1/2\right] ^{d}}=\prod_{j=1}^{d}\left( 1-\left|
x_{j}\right| \right) _{+},
\end{equation*}
gives rise to the piecewise linear interpolation operator. Analogously, $B$%
-splines of order $r\in \mathbb{N}$ are defined iterating this convolution $%
r $ times. These are $C^{r-1}\left( \mathbb{R}^{d}\right) $ functions
supported in $\left[ -r/2,r/2\right] ^{d}$, taking the value $1$ at the
origin. More details may be found, for instance, in \cite{dB}.\medskip

$\bullet $ \textbf{Band-limited sampling/reconstruction. }This corresponds
to the profile
\begin{equation*}
\varphi \left( \xi \right) :=\prod_{j=1}^{d}\func{sinc}\left( \xi
_{j}\right) ,
\end{equation*}
where the \textbf{cardinal sine function} is defined by
\begin{equation*}
\func{sinc}\left( t\right) :=\frac{\sin \pi t}{\pi t}.
\end{equation*}
It is easy to check that $\widehat{\varphi }\left( \xi \right) =\mathbf{1}%
_{Q}\left( \xi \right) $. This profile is relevant because of \textbf{%
Shannon's sampling Theorem}: \emph{a function }$u$\emph{\ belongs to the
space}
\begin{equation*}
V^{h}:=\left\{ u\in L^{2}\left( \mathbb{R}^{d}\right) :\limfunc{supp}%
\widehat{u}\subset \left[ -\pi /h,\pi /h\right) ^{d}\right\} =\text{range}%
\left( T_{\varphi }^{h}\right)
\end{equation*}
\emph{if and only if}
\begin{equation*}
u=\sum_{n\in \mathbb{Z}^{d}}u\left( hn\right) \varphi _{n}^{h}.
\end{equation*}
In particular, such functions are determined by their values on the grid $h%
\mathbb{Z}^{d}$.\medskip

$\bullet $ \textbf{Wavelets. }Take again $h_{k}:=2^{-k}$ for every $k\in
\mathbb{Z}$. A function $\psi \in L^{2}\left( \mathbb{R}^{d}\right) $ is a
wavelet provided $\left\{ \psi _{n}^{h_{k}}:n\in \mathbb{Z}^{d},k\in \mathbb{%
Z}\right\} $ is an orthonormal basis of $L^{2}\left( \mathbb{R}^{d}\right) $%
. For more details on wavelets and the closely related \textbf{%
MultiResolution Analyses}, the reader may see \cite{H-W, Me1}.\bigskip

Additional examples and references (from the viewpoint of Signal Theory),
may be found in the survey \cite{U}.

\subsection{Boundedness properties}

In order to ensure that the sampling and reconstruction operators are
bounded, we shall make the assumption (\textbf{\ref{BP}}) below:
\begin{equation}
\begin{array}{c}
\text{There exist }s\in \mathbb{R}\text{ and }B>0\text{ such that }\varphi
\in H^{s}\left( \mathbb{R}\right) \text{ and}\medskip \\
\tau _{\left\langle D_{x}\right\rangle ^{s}\varphi }\left( \xi \right)
:=\dsum\limits_{k\in \mathbb{Z}^{d}}\left| \left\langle \xi +2\pi
k\right\rangle ^{s}\widehat{\varphi }\left( \xi +2\pi k\right) \right|
^{2}\leq B\text{\quad for a.e. }\xi \in \mathbb{R}^{d}.
\end{array}
\tag{\QTR{bf}{BP}}  \label{BP}
\end{equation}

\begin{lemma}
\label{Lemma norm TphiUh}Suppose $\varphi \in \mathcal{S}^{\prime }\left(
\mathbb{R}^{d}\right) $. Then the following are equivalent:\bigskip

i) $\varphi $ satisfies (\ref{BP}).\bigskip

ii) There exists $B>0$ such that
\begin{equation}
\left\| \left\langle hD_{x}\right\rangle ^{s}T_{\varphi }^{h}U\right\|
_{L^{2}\left( \mathbb{R}^{d}\right) }\leq \sqrt{B}\left\| U\right\|
_{L^{2}\left( h\mathbb{Z}^{d}\right) }  \label{estimate norm TphiUh}
\end{equation}
holds uniformly for $h>0$ and $U\in L^{2}\left( h\mathbb{Z}^{d}\right) $%
.\bigskip

iii) There exists $B>0$ such that
\begin{equation}
\left\| S_{\varphi }^{h}u\right\| _{L^{2}\left( h\mathbb{Z}^{d}\right) }\leq
\sqrt{B}\left\| \left\langle hD_{x}\right\rangle ^{-s}u\right\|
_{L^{2}\left( \mathbb{R}^{d}\right) }  \label{estimate norm Sphiuh}
\end{equation}
holds uniformly for $h>0$ and $u\in H^{-s}\left( \mathbb{R}^{d}\right) $%
.\bigskip

Moreover, whenever i), ii) or iii) is fulfilled, the smallest constant $B$
for which any of the above assertion holds is precisely $\left\| \tau
_{\left\langle D_{x}\right\rangle ^{s}\varphi }\right\| _{L^{\infty }\left(
Q\right) }$.
\end{lemma}

\begin{proof}
To see why i) and ii) are equivalent, first observe that, given any $\varphi
\in H^{s}\left( \mathbb{R}^{d}\right) $ the following identity holds
\begin{equation}
T_{\varphi }^{h}=\left\langle hD_{x}\right\rangle ^{-s}T_{\left\langle
D_{x}\right\rangle ^{s}\varphi }^{h}.  \label{Thphi factorized}
\end{equation}
To check this, simply notice that
\begin{equation*}
\widehat{T_{\varphi }^{h}U}\left( \xi \right) =h^{d}\widehat{\varphi }\left(
h\xi \right) \sum_{n\in \mathbb{Z}^{d}}U_{n}e^{-ihn\cdot \xi }=\widehat{%
\varphi }\left( h\xi \right) h^{d}\widehat{U}\left( h\xi \right) ,
\end{equation*}
hence
\begin{equation*}
\widehat{T_{\varphi }^{h}U}\left( \xi \right) =\left\langle h\xi
\right\rangle ^{-s}\left\langle h\xi \right\rangle ^{s}\widehat{\varphi }%
\left( h\xi \right) h^{d}\widehat{U}\left( h\xi \right) =\widehat{%
\left\langle hD_{x}\right\rangle ^{-s}T_{\left\langle D_{x}\right\rangle
^{s}\varphi }^{h}U}\left( \xi \right) .
\end{equation*}
Since $\left\langle D_{x}\right\rangle ^{s}\varphi \in L^{2}\left( \mathbb{R}%
^{d}\right) $ and
\begin{equation*}
\left\| \left\langle hD_{x}\right\rangle ^{s}T_{\varphi }^{h}U\right\|
_{L^{2}\left( \mathbb{R}^{d}\right) }=\left\| T_{\left\langle
D_{x}\right\rangle ^{-s}\varphi }^{h}U\right\| _{L^{2}\left( \mathbb{R}%
^{d}\right) }
\end{equation*}
it suffices to deal with the case $s=0$. But it is a well-known result (see
for instance \cite{dBo-dVo-R, Ron}) that for $\varphi \in L^{2}\left(
\mathbb{R}^{d}\right) $, i) and ii) are equivalent and that $\left\|
T_{\varphi }^{h}\right\| =\left\| \tau _{\varphi }\right\| _{L^{\infty
}\left( \mathbb{R}^{d}\right) }$ whenever $T_{\varphi }^{h}$ is bounded.

Statements ii) and iii) are equivalent because of the following duality
relation:
\begin{equation*}
\left( \left\langle hD_{x}\right\rangle ^{s}T_{\varphi }^{h}U,\left\langle
hD_{x}\right\rangle ^{-s}u\right) _{L^{2}\left( \mathbb{R}^{d}\right)
}=\left( U,S_{\varphi }^{h}u\right) _{L^{2}\left( h\mathbb{Z}^{d}\right) },
\end{equation*}
which holds for every $u\in H^{-s}\left( \mathbb{R}^{d}\right) $ and $U\in
L^{2}\left( h\mathbb{Z}^{d}\right) $. This is simple to check:

\begin{align*}
\left( \left\langle hD_{x}\right\rangle ^{s}T_{\varphi }^{h}U,\left\langle
hD_{x}\right\rangle ^{-s}u\right) _{L^{2}\left( \mathbb{R}^{d}\right)
}=\sum_{n\in \mathbb{Z}^{d}}U_{n}\int_{\mathbb{R}^{d}}\left\langle
hD_{x}\right\rangle ^{s}\varphi _{n}^{h}\left( x\right) \overline{%
\left\langle hD_{x}\right\rangle ^{-s}u\left( x\right) }dx \\
=\sum_{n\in \mathbb{Z}^{d}}U_{n}\left\langle \varphi _{n}^{h},\overline{u}%
\right\rangle _{H^{s}\left( \mathbb{R}^{d}\right) \times H^{-s}\left(
\mathbb{R}^{d}\right) } \\
=h^{d}\sum_{n\in \mathbb{Z}^{d}}U_{n}^{h}\overline{S_{\varphi }^{h}u\left(
n\right) }.
\end{align*}
\medskip
\end{proof}

\begin{remark}
\label{Rmk boundedness ThU}For $s\leq 0$, estimate (\ref{estimate norm
TphiUh}) implies that
\begin{equation}
\left\| \left\langle \varepsilon D_{x}\right\rangle ^{s}T_{\varphi
}^{h}U^{h}\right\| _{L^{2}\left( \mathbb{R}^{d}\right) }\leq \sqrt{B}\left\|
U^{h}\right\| _{L^{2}\left( h\mathbb{Z}^{d}\right) },
\label{estimate norm TphiUh epsilon}
\end{equation}
a soon as $h/\varepsilon \leq 1$, as it can be easily checked taking Fourier
transforms.
\end{remark}

A sufficient condition for (\ref{BP}) in terms of decay on $\varphi $ is
next given:

\begin{lemma}
\label{Lemma suffBP}Suppose $\varphi \in H^{s}\left( \mathbb{R}^{d}\right) $
satisfies, for some $\varepsilon >0$,
\begin{equation}
\int_{\mathbb{R}^{d}}\left| \left\langle D_{x}\right\rangle ^{s}\varphi
\left( x\right) \right| ^{2}\left( 1+\left| x\right| \right) ^{d+\varepsilon
}dx<\infty  \label{strong shyp}
\end{equation}
Then $\widehat{\varphi }$ and $\tau _{\left\langle D_{x}\right\rangle
^{s}\varphi }$ are continuous functions. In particular, (\ref{BP}) always
holds for such a $\varphi $.
\end{lemma}

\begin{proof}
It follows the lines of \cite{Me1}, Lemma II.7. Under condition (\ref{strong
shyp}), $\left\langle \xi \right\rangle ^{s}\widehat{\varphi }\in
H^{d/2+\varepsilon /2}\left( \mathbb{R}^{d}\right) $; Sobolev's imbedding
Theorem then ensures that $\left\langle \xi \right\rangle ^{s}\widehat{%
\varphi }$ is a continuous function and hence so is $\widehat{\varphi }$.
The continuity of $\tau _{\left\langle D_{x}\right\rangle ^{s}\varphi }$ is
a consequence of the fact that, whenever $\chi \in C_{c}^{\infty }\left(
\mathbb{R}^{d}\right) $ satisfies $\sum_{n\in \mathbb{Z}^{d}}\left| \chi
\left( \xi +2\pi n\right) \right| \geq 1$, the expression
\begin{equation*}
\left[ \sum_{n\in \mathbb{Z}^{d}}\left\| u\chi \left( \cdot +2\pi n\right)
\right\| _{H^{s}\left( \mathbb{R}^{d}\right) }^{2}\right] ^{1/2}
\end{equation*}
defines an equivalent norm in $H^{s}\left( \mathbb{R}^{d}\right) $, $s\geq 0$%
. This actually proves that
\begin{equation*}
\sum_{n\in \mathbb{Z}^{d}}\sup_{\xi \in \mathbb{R}^{d}}\left| \left\langle
\xi \right\rangle ^{s}\widehat{\varphi }\left( \xi \right) \chi \left( \xi
+2\pi n\right) \right| ^{2}<\infty .
\end{equation*}
In particular, the series defining $\tau _{\left\langle D_{x}\right\rangle
^{s}\varphi }$ is uniformly convergent and the claim then follows.\medskip
\end{proof}

Condition (\ref{strong shyp}) automatically holds for profiles $\varphi $
such that
\begin{equation}
\left| \left\langle D_{x}\right\rangle ^{s}\varphi \left( x\right) \right|
\leq C\left( 1+\left| x\right| \right) ^{-d-\varepsilon },\text{\qquad for
every }x\in \mathbb{R}^{d}\text{ and some }C,\varepsilon >0;
\label{strong shyp2}
\end{equation}
in particular, the hypothesis (\ref{strong hyp}) we assumed in the
introduction implies (\ref{BP}) for $s=0$.

Now we can prove a general result from which Proposition \ref{Thm definition
wmD} immediately follows:

\begin{proposition}
\label{Prop Wm ThU}Suppose $\varphi $ satisfies (\ref{BP}) and we are given
scales $\left( h_{k}\right) $, $\left( \varepsilon _{k}\right) $ such that $%
\left( h_{k}/\varepsilon _{k}\right) $ is bounded.\textbf{\ }If $\left(
U^{h_{k}}\right) $ is an $h_{k}$-bounded sequence of discrete functions then
the distributions $m^{\varepsilon _{k}}\left[ T_{\varphi }^{h_{k}}U^{h_{k}}%
\right] $ are uniformly bounded in $\mathcal{S}^{\prime }$. Moreover, the
limit of any weakly convergent subsequence is a positive measure.
\end{proposition}

The proof this is a direct consequence of Remark \ref{Rmk boundedness ThU}
and the general result established in Proposition \ref{Prop definition wm
l2loc}.

\subsection{Bases and projections}

Below, we recall some results from Approximation Theory that will be needed
in the sequel. These results deal with the range in $H^{s}\left( \mathbb{R}%
^{d}\right) $ of the reconstruction operator $T_{\varphi }^{h}$, which we
denote $V_{\varphi }^{h}$.

The space $V_{\varphi }^{h}$ is a \textbf{Principal Shift Invariant (PSI) }%
space. When any of the conditions of Lemma \ref{Lemma norm TphiUh} are
satisfied, the family $\left\{ h^{-d/2}\varphi _{n}^{h}:n\in \mathbb{Z}%
^{d}\right\} $ is said to form a \textbf{Bessel system }for $V_{\varphi
}^{h} $.

The next Lemma clarifies how the function $\tau _{\left\langle
D_{x}\right\rangle ^{s}\varphi }$ characterizes further basis properties of
the functions $\varphi _{n}^{h}$.

\begin{lemma}
\label{Lemma Riesz ort}Let $\varphi \in \mathcal{S}^{\prime }\left( \mathbb{R%
}^{d}\right) $ satisfy (\ref{BP}). Then\bigskip

i) $\left\{ h^{-d/2}\varphi _{n}^{h}:n\in \mathbb{Z}^{d}\right\} $ is an
orthonormal basis of $V_{\varphi }^{h}$ if and only if
\begin{equation*}
\tau _{\left\langle D_{x}\right\rangle ^{s}\varphi }\left( \xi \right) =1%
\text{\qquad for a.e. }\xi \in \mathbb{R}^{d}\text{.}
\end{equation*}

ii) $\left\{ h^{-d/2}\varphi _{n}^{h}:n\in \mathbb{Z}^{d}\right\} $ is a
\textbf{Riesz basis} \footnote{%
This means that there exist constants $A,B>0$ such that
\begin{equation*}
A\left\| U\right\| _{L^{2}\left( h\mathbb{Z}^{d}\right) }^{2}\leq \left\|
T_{\varphi }^{h}U\right\| _{H^{s}\left( \mathbb{R}^{d}\right) }^{2}\leq
B\left\| U\right\| _{L^{2}\left( h\mathbb{Z}^{d}\right) }^{2}
\end{equation*}
for all $U\in L^{2}\left( h\mathbb{Z}^{d}\right) $. This is equivalent to
the existence of a linear isomorphism $R:V^{h}\rightarrow V^{h}$ such that $%
\left\{ h^{-d/2}R\varphi _{n}^{h}:n\in \mathbb{Z}^{d}\right\} $ forms an
orthonormal basis of $H^{s}\left( \mathbb{R}^{d}\right) $. This property is
sometimes also referred as that $\left( \varphi _{n}^{h}\right) _{n\in
\mathbb{Z}^{d}}$ form a \textbf{stable frame }in $H^{s}\left( \mathbb{R}%
^{d}\right) $.} of $V_{\varphi }^{h}$ if and only if there exist constants $%
A,B>0$ such that
\begin{equation*}
A\leq \tau _{\left\langle D_{x}\right\rangle ^{s}\varphi }\left( \xi \right)
\leq B\text{\qquad for a.e. }\xi \in \mathbb{R}^{d}\text{.}
\end{equation*}
\end{lemma}

\begin{proof}
The operator $\left\langle D_{x}\right\rangle ^{s}:H^{s}\left( \mathbb{R}%
^{d}\right) \rightarrow L^{2}\left( \mathbb{R}^{d}\right) $ is unitary.
Hence $\left\{ h^{-d/2}\varphi _{n}^{h}:n\in \mathbb{Z}^{d}\right\} $ is an
orthonormal (resp. Riesz) basis of $V_{\varphi }^{h}$ if and only if $%
\left\{ h^{-d/2}\left( \left\langle D_{x}\right\rangle ^{s}\varphi \right)
_{n}^{h}:n\in \mathbb{Z}^{d}\right\} $ is an orthonormal (resp. Riesz) basis
of the range of $T_{\left\langle D_{x}\right\rangle ^{s}\varphi }^{h}$.
Thus, the Lemma needs only to be proved for profiles $\varphi \in
L^{2}\left( \mathbb{R}^{d}\right) $ and this is a well-known result (see,
for instance, \cite{Ron}).\bigskip
\end{proof}

We shall also need the following expression for the orthogonal projection
onto $V_{\varphi}^{h}$:

\begin{lemma}
\label{Lemma ort}Let $\varphi \in \mathcal{S}^{\prime }\left( \mathbb{R}%
^{d}\right) $ satisfy (\ref{BP}). The orthogonal projection $P_{\varphi
}^{h}:H^{s}\left( \mathbb{R}^{d}\right) \rightarrow V_{\varphi }^{h}$ equals
$P_{\varphi }^{h}=T_{\varphi }^{h}S_{\widetilde{\left\langle
D_{x}\right\rangle ^{s}\varphi }}^{h}\left\langle hD_{x}\right\rangle ^{s}$,
where, for $f\in L^{2}\left( \mathbb{R}^{d}\right) $, $\tilde{f}\in
L^{2}\left( \mathbb{R}^{d}\right) $ is defined by:
\begin{equation*}
\widehat{\tilde{f}}\left( \xi \right) :=\left\{
\begin{array}{ll}
\dfrac{\widehat{f}\left( \xi \right) }{\tau _{f}\left( \xi \right) }, &
\quad \text{if }\tau _{f}\left( \xi \right) \neq 0,\medskip \\
0 & \quad \text{otherwise.}
\end{array}
\right.
\end{equation*}
\end{lemma}

\begin{proof}
The proof of the result for $s=0$ may be found in \cite{dBo-dVo-R}, Theorem
2.9. We can reduce ourselves to this case by noticing that
\begin{equation*}
P_{\varphi }^{h}=\left\langle hD_{x}\right\rangle ^{-s}P_{\left\langle
hD_{x}\right\rangle ^{s}\varphi }^{h}\left\langle hD_{x}\right\rangle ^{s},
\end{equation*}
since, as we have seen in (\ref{Thphi factorized}), the range of $T_{\varphi
}^{h}$ equals that of $\left\langle hD_{x}\right\rangle ^{-s}T_{\left\langle
hD_{x}\right\rangle ^{s}\varphi }^{h}$ and $\left\langle hD_{x}\right\rangle
^{s}$ is an orthogonal mapping. Using the $L^{2}$-result we obtain:
\begin{equation*}
P_{\varphi }^{h}=\left\langle hD_{x}\right\rangle ^{-s}T_{\left\langle
hD_{x}\right\rangle ^{s}\varphi }^{h}S_{\widetilde{\left\langle
D_{x}\right\rangle ^{s}\varphi }}^{h}\left\langle hD_{x}\right\rangle
^{s}=T_{\varphi }^{h}S_{\widetilde{\left\langle D_{x}\right\rangle
^{s}\varphi }}^{h}\left\langle hD_{x}\right\rangle ^{s},
\end{equation*}
as claimed.\medskip
\end{proof}

\section{\label{Sec h=e}High frequency analysis: $h\sim \protect\varepsilon $%
}

\subsection{Reduction to the case $h=\protect\varepsilon$}

In this section we analyze the effect of sampling and reconstruction on
Wigner measures at scales $\left( \varepsilon_{k}\right) $, of the same
order of the sampling/reconstruction rate $\left( h_{k}\right) $ (i.e., such
that $\left( h_{k}/\varepsilon_{k}\right) $ is bounded).

First notice that it suffices to treat the case $\varepsilon_{k}=h_{k}$; the
more general one can be obtained by a proper rescaling. This is due to the
following identity:
\begin{equation*}
m^{\varepsilon}\left[ u\right] \left( x,\xi\right) =\left( h/\varepsilon
\right) ^{d}m^{h}\left[ u\right] \left( x,\left( h/\varepsilon\right)
\xi\right) ,
\end{equation*}
which clearly implies:

\begin{lemma}
Suppose $h_{k}/\varepsilon _{k}\rightarrow c>0$. Then $m^{\varepsilon _{k}}%
\left[ u_{k}\right] $ converges in $\mathcal{S}^{\prime }$ if and only if $%
m^{h_{k}}\left[ u_{k}\right] $ does. Their respective limits $\mu _{c}$ and $%
\mu $ are related through:
\begin{equation}
\mu _{c}\left( x,\xi \right) =c^{d}\mu \left( x,c\xi \right) .  \label{muc}
\end{equation}
\end{lemma}

When $h_{k}=\varepsilon _{k}$, the transforms $M^{h_{k}}\left[ U^{h_{k}}%
\right] $ are $\Gamma $-periodic in the variable $\xi $; hence, so are their
limiting Wigner measures

\subsection{Sampling}

We start by exploring the effect of sampling on the structure of Wigner
measures. The computation of the Wigner measure at scale $\left(
h_{k}\right) $ of a sequence of samples $\left( S_{\varphi
}^{h_{k}}u_{k}\right) $ is done in the following theorem; it is applicable
whenever the hypothesis (\textbf{\ref{h-osc}}) below is fulfilled:
\begin{equation}
\limfunc{essup}\limits_{\xi \in Q}\sum_{\left| n\right| \geq R}\left|
\left\langle \xi +2\pi n\right\rangle ^{s}\widehat{\varphi }\left( \xi +2\pi
n\right) \right| ^{2}\rightarrow 0\qquad \text{as }R\rightarrow \infty .
\tag{\QTR{bf}{D}}  \label{h-osc}
\end{equation}
Notice that profiles with the property (\ref{strong shyp}) immediately
verify (\ref{h-osc}).

Before stating our result, it is important to notice that the Fourier
transform of a profile $\varphi $ satisfying condition (\ref{BP}) is an
element of $L_{\text{loc}}^{2}\left( \mathbb{R}^{d}\right) $. In particular,
it is only defined modulo a set of zero Lebesgue measure. Thus, when dealing
with pointwise properties of $\widehat{\varphi }$, we shall systematically
assume that a precise representative of the class of $\widehat{\varphi }$
has, once for all, been chosen.

For instance, the Wigner measures $\mu $ of the sequences $\left(
u_{k}\right) $ in Theorem \ref{Thm WM Suh uh} below, will be assumed to
satisfy conditions (\textbf{\ref{singular measures}}) and (\textbf{\ref{ND}}%
):
\begin{equation}
\left| \widehat{\varphi }\left( \xi +2\pi n\right) \right| ^{2}\mu \left(
x,\xi +2\pi n\right) ,\quad n\in \mathbb{Z}^{d},\quad \text{are mutually
singular measures}.  \tag{\QTR{bf}{MS}}  \label{singular measures}
\end{equation}
\begin{equation}
\mu \left( \mathbb{R}^{d}\times \overline{D_{\widehat{\varphi }}}\right) =0,
\tag{\QTR{bf}{ND}}  \label{ND}
\end{equation}
where, recall, $D_{\widehat{\varphi }}$ stands for the set of discontinuity
points of $\widehat{\varphi }$. These conditions must be understood to hold
for the same representative of $\widehat{\varphi }$.

\begin{theorem}
\label{Thm WM Suh uh}Let $\left( h_{k}\right) $ be a scale and take $\varphi
$ satisfying (\ref{BP}) and (\ref{h-osc}). Let $\left( u_{k}\right) $ be a
sequence in $H^{-s}\left( \mathbb{R}^{d}\right) $ such that $\left(
\left\langle h_{k}D_{x}\right\rangle ^{-s}u_{k}\right) $ is bounded and
suppose that $m^{h_{k}}\left[ u_{k}\right] $ converges to a Wigner measure $%
\mu $ that fulfills (\ref{ND}), (\ref{singular measures}).$\ $

Then $M^{h_{k}}\left[ S_{\varphi }^{h_{k}}u_{k}\right] $ converges to the
Wigner measure $\mu ^{\varphi }$ given by:
\begin{equation}
\mu ^{\varphi }\left( x,\xi \right) =\sum_{n\in \mathbb{Z}^{d}}\left|
\widehat{\varphi }\left( \xi +2\pi n\right) \right| ^{2}\mu \left( x,\xi
+2\pi n\right) .  \label{relation meuh and MeSUh h=e}
\end{equation}
\end{theorem}

\begin{remark}
i) As pointed out above, formula (\ref{muphi and mu h equal e}) holds for
the same precise representative of the Fourier transform $\widehat{\varphi }$
which was chosen in (\ref{ND}) and (\ref{singular measures}). \smallskip

ii) The necessity of hypotheses (\ref{ND}) and (\ref{singular measures})
will be discussed in paragraph \ref{sbs nec}.\smallskip

iii) Condition (\ref{h-osc}) may be replaced by the assumption that $\left(
\left\langle h_{k}D_{x}\right\rangle ^{-s}u_{k}\right) $ is $h_{k}$%
-oscillatory. This will be made clear in the proof of the Theorem.\smallskip

iv) The boundedness of $\left\langle h_{k}D_{x}\right\rangle ^{-s}u_{k}$
implies that $\left( u_{k}\right) $ is $h_{k}$-oscillatory.
\end{remark}

The proof of this Theorem is postponed to the end of this section.

The expression (\ref{relation meuh and MeSUh h=e}) may be related to the
concept of \textbf{Wigner series} introduced in \cite{M-M-P, G-M-M-P}.
Recall that given $u\in \mathcal{S}^{\prime }\left( \mathbb{R}^{d}\right) $,
the \textbf{Wigner series} of $u$ at scale $\varepsilon $ is defined by:

\begin{equation*}
w_{S}^{\varepsilon }\left[ u\right] \left( x,\xi \right) :=\frac{1}{\left(
2\pi \right) ^{d}}\sum_{n\in \mathbb{Z}^{d}}u\left( x-\varepsilon \pi
n\right) \overline{u}\left( x+\varepsilon \pi n\right) e^{in\cdot \xi }.
\end{equation*}
It is easy to check that $w_{S}^{\varepsilon }\left[ u\right] \left( x,\xi
\right) =\sum_{n\in \mathbb{Z}^{d}}w^{\varepsilon }\left[ u\right] \left(
x,\xi +2\pi n\right) $.\footnote{%
See (\ref{definition WT}) for the definition of the Wigner transform $%
w^{\varepsilon }\left[ u\right] $.}

When $\left( u_{k}\right) $ is bounded in $L^{2}\left( \mathbb{R}^{d}\right)
$, $\varepsilon _{k}$-oscillatory and possesses a Wigner measure at scale $%
\left( \varepsilon _{k}\right) $ then the following relation holds:
\begin{equation}
\lim_{k\rightarrow \infty }\int_{\mathbb{R}^{d}\times \mathbb{R}^{d}}a\left(
x,\xi \right) w_{S}^{\varepsilon _{k}}\left[ u_{k}\right] \left( x,\xi
\right) dxd\xi =\int_{\mathbb{R}^{d}\times \mathbb{R}^{d}}\sum_{n\in \mathbb{%
Z}^{d}}a\left( x,\xi +2\pi n\right) d\mu \left( x,\xi \right) ,
\label{ws period}
\end{equation}
for $a\in \mathcal{S}$, see \cite{Bra}.

Theorem \ref{Thm WM Suh uh} has a simple interpretation in terms of Wigner
series: the measure $\mu ^{\varphi }$ may be obtained as the limit of the
Wigner series
\begin{equation*}
w_{S}^{h_{k}}\left[ \widehat{\varphi }\left( h_{k}D_{x}\right) u_{k}\right] .
\end{equation*}
This is due to the fact that, under any of the hypotheses (\ref{h-osc}), $%
\left( \widehat{\varphi }\left( h_{k}D_{x}\right) u_{k}\right) $ is $h_{k}$%
-oscillatory. Besides, as a consequence of Proposition \ref{Prop
localization}, the Wigner measure at scale $\left( h_{k}\right) $ of $\left(
\widehat{\varphi }\left( h_{k}D_{x}\right) u_{k}\right) $ is given by $%
\left| \widehat{\varphi }\left( \xi \right) \right| ^{2}\mu \left( x,\xi
\right) $. The assertion then follows from (\ref{ws period}).

As was already mentioned in the Introduction, condition (\ref{singular
measures}) is a restriction on the support of the measure $\left| \widehat{%
\varphi }\left( \xi \right) \right| ^{2}\mu \left( x,\xi \right) $. Two
extremal cases in which it is trivially satisfied are the following:\medskip

i) $\widehat{\varphi }|_{\mathbb{R}^{d}\setminus Q}\equiv 0$, in this case (%
\ref{singular measures}) holds independently of what $\mu $ is.\medskip

ii) The sequence $\left( u_{k}\right) $ is \textbf{asymptotically
band-limited } i.e. its Wigner measures at scale $\left( h_{k}\right) $ is
concentrated on the cube $\overline{Q}$. For those sequences, condition (\ref
{singular measures}) only involves the behavior of $\mu $ on the boundary $%
\partial Q$: it essentially expresses that the restrictions of $\mu $ to
parallel sides of $\partial Q$ do not overlap (i.e. are mutually singular).
A sufficient condition for this is, for instance,
\begin{equation}
\limsup_{k\rightarrow \infty }\int_{\mathbb{R}^{d}\setminus Q_{R}}\left|
\widehat{u_{k}}\left( \frac{\xi }{h_{k}}\right) \right| ^{2}\frac{d\xi }{%
\left( 2\pi h_{k}\right) ^{d}}\rightarrow 0\qquad \text{as }R\rightarrow
\infty ,  \label{definition S h-oscillatory}
\end{equation}
where $Q_{R}:=\left[ -\pi ,\pi -1/R\right) ^{d}$.

\begin{remark}
In any of the above cases, we have:
\begin{equation*}
\mathbf{1}_{Q}\left( \xi \right) \mu ^{\varphi }\left( x,\xi \right) =\left|
\widehat{\varphi }\left( \xi \right) \right| ^{2}\mu \left( x,\xi \right) .
\end{equation*}
Hence, the restriction of $\mu ^{\varphi }$ to $\mathbb{R}^{d}\times Q$
coincides with $\mu $ if and only if $\left| \widehat{\varphi }\left( \xi
\right) \right| ^{2}=1$ for $\mu $-almost every $\xi \in \overline{Q}$.
\end{remark}

The specific choice $\varphi=\delta_{0}$ corresponds to the analysis of
\textbf{discretization}, for then $S_{\delta_{0}}^{h}u\left( n\right)
=u\left( hn\right) $. Theorem \ref{Thm WM Suh uh} takes the following simple
form:

\begin{corollary}
\label{Corollary discret}Let $\left( h_{k}\right) $ be a scale and let $%
\left( u_{k}\right) $ be a sequence in $H^{s}\left( \mathbb{R}^{d}\right) $,
for some $s>d/2$, such that $\left( \left\langle h_{k}D_{x}\right\rangle
^{s}u_{k}\right) $ is bounded. If $\mu $ is its Wigner measure at scale $%
\left( h_{k}\right) $ and the measures $\mu \left( x,\xi +2\pi n\right) $
are mutually singular then Wigner measure $\mu ^{\delta _{0}}$ corresponding
to the sequence of discretizations is the periodization:
\begin{equation*}
\mu ^{\delta _{0}}\left( x,\xi \right) =\sum_{n\in \mathbb{Z}^{d}}\mu \left(
x,\xi +2\pi n\right) .
\end{equation*}
In other words, $\mu ^{\delta _{0}}$ is the limit of the Wigner series $%
w_{S}^{h_{k}}\left[ u_{k}\right] $.
\end{corollary}

This Corollary is particularly useful in the explicit computation of Wigner
measures for discrete functions. As an example, consider the concentrating
and oscillating sequences we defined in the Introduction, $f_{k}\left(
x\right) =k^{d/2}\rho \left( k\left( x-x_{0}\right) \right) $ and $%
g_{k}\left( x\right) :=\rho \left( x\right) e^{ikx\cdot \xi ^{0}}$ with $%
\rho \in L^{2}\left( \mathbb{R}^{d}\right) $. Using identities (\ref{examp
con}) and (\ref{examp osc}) we obtain, for $\left( f_{k}\right) $ and $%
\left( g_{k}\right) $ respectively:
\begin{equation*}
\mu ^{\delta _{0}}\left( x,\xi \right) =\delta _{x_{0}}\left( x\right)
\otimes \sum_{n\in \mathbb{Z}}\left| \widehat{\rho }\left( \xi +2\pi
n\right) \right| ^{2}\dfrac{d\xi }{\left( 2\pi \right) ^{d}},
\end{equation*}
if, for instance, $\limfunc{supp}\widehat{\rho }\subset Q$, and
\begin{equation}
\left| \rho \left( x\right) \right| ^{2}dx\otimes \sum_{n\in \mathbb{Z}%
^{d}}\delta _{\xi ^{0}+2\pi n}\left( \xi \right) ,  \label{examp osc D}
\end{equation}
with no assumption on $\rho $.

\subsection{Reconstruction}

Now we deal with the reconstruction operator $T_{\varphi }^{h}$; it modifies
the high-frequency behavior of a sequence of discrete functions in the
following way:

\begin{theorem}
\label{Thm WM TUh Uh}Let $\left( h_{k}\right) $ be a scale and $\left(
U^{h_{k}}\right) $ be an $h_{k}$-bounded sequence; take $\varphi $
satisfying (\ref{BP}). If $M^{h_{k}}\left[ U^{h_{k}}\right] $ converges to
the Wigner measure $\mu $ which verifies (\ref{ND}) then $m^{h_{k}}\left[
T_{\varphi }^{h_{k}}U^{h_{k}}\right] $ converges to a Wigner measure $\mu
_{\varphi }$ given by:
\begin{equation}
\mu _{\varphi }\left( x,\xi \right) =\left| \widehat{\varphi }\left( \xi
\right) \right| ^{2}\mu \left( x,\xi \right) .
\label{muphi and mu h equal e}
\end{equation}
\end{theorem}

The proof of Theorem \ref{Thm WM TUh Uh} is based on explicit formulas for
the Fourier transforms of $T_{\varphi }^{h}U$. As we have already seen,
\begin{equation}
\widehat{T_{\varphi }^{h}U}\left( \xi \right) =\widehat{\varphi }\left( h\xi
\right) h^{d}\widehat{U}\left( h\xi \right) ,  \label{FT of ThphiUh}
\end{equation}
for any $U\in L^{2}\left( h\mathbb{Z}^{d}\right) $. The following Remark
ensures that Proposition \ref{Prop localization} can be applied in the proof
below.

\begin{remark}
\label{Rmk phihat in linf}If $\varphi \in H^{s}\left( \mathbb{R}^{d}\right) $
satisfies (\ref{BP})\textbf{\ }then $\widehat{\varphi }\in L^{\infty }\left(
\mathbb{R}^{d};\left\langle \xi \right\rangle ^{s}\right) $.
\end{remark}

\begin{proof}[Proof of Theorem \ref{Thm WM TUh Uh}]
Just notice that (\ref{FT of ThphiUh}) can be rewritten as:
\begin{equation*}
T_{\varphi }^{h}U^{h}=\widehat{\varphi }\left( hD_{x}\right) T_{\delta
_{0}}^{h}U^{h}.
\end{equation*}
The hypotheses made on $\varphi $ and $\mu $ allow us to apply Proposition
\ref{Prop localization} (see Remark \ref{Rmk phihat in linf}) and
conclude\medskip
\end{proof}

Identity (\ref{muphi and mu h equal e}) expresses how the measure $\mu $ is
modulated by the profile $\varphi $; the necessity of the hypothesis (\ref
{ND}) for this result is discussed in paragraph \ref{sbs nec} as well.

Since $\mu$ is $\Gamma$-periodic in $\xi$, formula (\ref{muphi and mu h
equal e}) suggests that $\mu$ may be compared to the periodization of $%
\mu_{\varphi }$ with respect to the variable $\xi$.

\begin{corollary}
\label{Corollary period}Let $\varphi $, $\left( U^{h_{k}}\right) $, $\mu $
and $\mu _{\varphi }$ be as in Theorem \ref{Thm WM TUh Uh}. Then the
periodization
\begin{equation}
\mu _{\varphi ,s}\left( x,\xi \right) :=\sum_{n\in \mathbb{Z}%
^{d}}\left\langle \xi +2\pi n\right\rangle ^{2s}\mu _{\varphi }\left( x,\xi
+2\pi n\right)  \label{periods}
\end{equation}
is a well-defined\footnote{%
The limit defining the sum (\ref{periods}) is understood to exist for the
weak convergence of measures in $\mathcal{M}_{+}\left( \mathbb{R}^{d}\times
\mathbb{R}^{d}\right) $.} measure, $\Gamma $-periodic in $\xi $, that
satisfies:
\begin{equation}
\mu _{\varphi ,s}\left( x,\xi \right) =\tau _{\left\langle
D_{x}\right\rangle ^{s}\varphi }\left( \xi \right) \mu \left( x,\xi \right) .
\label{iden ws dwm}
\end{equation}
In particular:\smallskip

i) If $\tau _{\left\langle D_{x}\right\rangle ^{s}\varphi }\left( \xi
\right) =1$ except for $\xi $ in a set of zero $\mu $-measure then $\mu
_{\varphi ,s}=\mu $.\smallskip

ii) $\tau _{\left\langle D_{x}\right\rangle ^{s}\varphi }\equiv 1$ if and
only if the identity $\mu _{\varphi ,s}=\mu $ holds for every sequence $%
\left( U^{h_{k}}\right) $.
\end{corollary}

\begin{proof}
Since $\left| \left\langle \xi \right\rangle ^{s}\widehat{\varphi }\left(
\xi \right) \right| ^{2}$ is a nonnegative continuous function, the series
defining $\tau _{\left\langle D_{x}\right\rangle ^{s}\varphi }\left( \xi
\right) $ converges absolutely for every $\xi $ on the support of $\mu $
(which consists of continuity points for $\widehat{\varphi }\left( \xi
\right) $). Thus, by the dominated convergence Theorem:
\begin{equation*}
\int_{\mathbb{R}^{d}\times \mathbb{R}^{d}}a\left( x,\xi \right) \tau
_{\varphi ,s}\left( \xi \right) d\mu \left( x,\xi \right) =\sum_{n\in
\mathbb{Z}^{d}}\int_{\mathbb{R}^{d}\times \mathbb{R}^{d}}a\left( x,\xi
\right) \left| \left\langle \xi +2\pi n\right\rangle ^{s}\widehat{\varphi }%
\left( \xi +2\pi n\right) \right| ^{2}d\mu \left( x,\xi \right)
\end{equation*}
for every $a\in C_{c}\left( \mathbb{R}^{d}\times \mathbb{R}^{d}\right) $.
Taking now into account (\ref{muphi and mu h equal e}) and the fact that $%
\mu $ is $\Gamma $-periodic in $\xi $, we find that
\begin{equation*}
\int_{\mathbb{R}^{d}\times \mathbb{R}^{d}}a\left( x,\xi \right) \tau
_{\varphi ,s}\left( \xi \right) d\mu \left( x,\xi \right) =\sum_{n\in
\mathbb{Z}^{d}}\int_{\mathbb{R}^{d}\times \mathbb{R}^{d}}a\left( x,\xi
\right) \left\langle \xi +2\pi n\right\rangle ^{2s}d\mu _{\varphi }\left(
x,\xi +2\pi n\right)
\end{equation*}
and the first part of the result follows.

Statement i) as well as the ``only if'' part of ii) are trivial. To obtain
the necessity in ii), just consider sequences of discrete functions whose
Wigner measures are of the form $\mu \left( x,\xi \right) =\nu \left(
x\right) \otimes \sum_{n\in \mathbb{Z}^{d}}\delta _{\xi ^{0}+2\pi n}$ (as (%
\ref{examp osc D}), for instance). Clearly, for $\mu _{\varphi ,s}=\mu $ to
hold for such a measure, we must have $\tau _{\left\langle
D_{x}\right\rangle ^{s}\varphi }\left( \xi ^{0}\right) =1$.\medskip
\end{proof}

\begin{remark}
i) Because of Lemma \ref{Lemma Riesz ort}, if relation $\mu _{\varphi
,s}=\mu $ holds for every $h_{k}$-bounded sequence of discrete functions
then the profile $\varphi $ has the property: $\left\{ h^{-d/2}\varphi
_{n}^{h}:n\in \mathbb{Z}^{d}\right\} $ is an orthonormal family in $%
H^{s}\left( \mathbb{R}^{d}\right) $ for every $h>0$.\smallskip

ii) However, the converse is not true, if $\varphi $ gives rise to an
orthonormal family then $\tau _{\left\langle D_{x}\right\rangle ^{s}\varphi
}\left( \xi \right) =1$ holds outside a set of null Lebesgue measure. If $%
\mu $ is supported on that set, identity $\mu _{\varphi ,s}=\mu $ may not
hold.
\end{remark}

As in the preceding section, our result has an interpretation in terms of
Wigner series. Under the conditions of Theorem \ref{Thm WM TUh Uh}, the
measure $\mu _{\varphi ,s}$ may be obtained as the limit as $k\rightarrow
\infty $ of the functions
\begin{equation*}
w_{S}^{h_{k}}\left[ \left\langle h_{k}D_{x}\right\rangle ^{s}T_{\varphi
}^{h_{k}}U^{h_{k}}\right] ,
\end{equation*}
provided $\left( \left\langle h_{k}D_{x}\right\rangle ^{s}T_{\varphi
}^{h_{k}}U^{h_{k}}\right) $ is $h_{k}$-oscillatory. Note however, that this
may not be the case for certain profiles $\varphi $ (see paragraph \ref{ce
osc}).

In particular, Corollary \ref{Corollary period} shows that the limits of $%
w_{S}^{h_{k}}\left[ T_{\varphi }^{h_{k}}U^{h_{k}}\right] $ and $M^{h_{k}}%
\left[ U^{h_{k}}\right] $ coincide if we chose, for instance, $\varphi :=%
\mathbf{1}_{\left[ -1/2,1/2\right) ^{d}}$.

\subsection{\label{sbs nec}The necessity of the hypotheses of Theorems \ref
{Thm WM Suh uh} and \ref{Thm WM TUh Uh}}

Formulas (\ref{muphi and mu h equal e}) and (\ref{relation meuh and MeSUh
h=e}) may not hold when $\widehat{\varphi }$ is not continuous and the
Wigner measure $\mu $ does not vanish on the closure of the set of
discontinuity points $D_{\widehat{\varphi }}$. We illustrate this with two
one-dimensional examples where
\begin{equation*}
\varphi \left( x\right) =\frac{\sin \pi x}{\pi x}.
\end{equation*}
We will chose $\mathbf{1}_{Q}$ as the representative of $\widehat{\varphi }$
for which the counterexamples will be built.\bigskip

\textbf{1. Necessity of condition (\ref{ND})} \textbf{in Theorem \ref{Thm WM
TUh Uh}}. Take $U^{h}$ to be the sequence discrete function of $L^{2}\left( h%
\mathbb{Z}^{d}\right) $ given by their Fourier transforms:
\begin{equation*}
\widehat{U^{h}}\left( \xi \right) :=\frac{1}{h}\sum_{n\in \mathbb{Z}}\mathbf{%
1}_{\left( -1,1\right) }\left( \frac{\xi -\left( 2n+1\right) \pi }{h}\right)
.
\end{equation*}
Then, denoting by $\mu $ the Wigner measure at scale $h$ of $\left(
U^{h}\right) $,
\begin{equation*}
\left| \widehat{\varphi }\left( \xi \right) \right| ^{2}\mu \left( x,\xi
\right) =\frac{\sin ^{2}\left( x\right) }{\pi ^{2}x^{2}}dx\otimes \delta
_{-\pi }\left( \xi \right) .
\end{equation*}
This measure differs from $\mu _{\varphi }$, which is given by:
\begin{equation*}
\mu _{\varphi }\left( x,\xi \right) =\frac{\sin ^{2}\left( x/2\right) }{\pi
^{2}x^{2}}dx\otimes \left[ \delta _{\pi }\left( \xi \right) +\delta _{-\pi
}\left( \xi \right) \right] .
\end{equation*}

\begin{remark}
i) The particular choice of the representative of $\widehat{\varphi }$ does
not play a role. Theorem \ref{Thm WM TUh Uh} still fails if we take as
representative of $\widehat{\varphi }$ the characteristic functions of $%
\left( -\pi ,\pi \right) ^{d}$ or $\left[ -\pi ,\pi \right] ^{d}$.\smallskip

ii) In particular, this example shows that even the two projections on $x$
and $\xi $ of the measures $\mu $ and $\mu _{\varphi }$ may differ.\medskip

iii) This also shows that the periodization in $\xi $ of $\mu _{\varphi }$
does not necessarily coincide with $\mu $, even when $\tau _{\varphi }=1$ as
is the case here. Thus the conclusion of Corollary \ref{Corollary period}
may fail when $\widehat{\varphi }$ is not continuous.
\end{remark}

Our counterexample to Theorem \ref{Thm WM Suh uh} is essentially the same as
the previous one:\bigskip

\textbf{2.}\emph{\ }\textbf{Necessity of condition (\ref{ND})} \textbf{in
Theorem \ref{Thm WM Suh uh}}. Define
\begin{equation*}
\widehat{v^{h}}\left( \xi \right) :=\mathbf{1}_{\left( -1,1\right) }\left(
\xi +\pi /h\right) .
\end{equation*}
Then, denoting by $\mu $ the Wigner measure at scale $h$ of $\left(
v^{h}\right) $,
\begin{equation*}
\sum_{n\in \mathbb{Z}}\left| \widehat{\varphi }\left( \xi +2\pi n\right)
\right| ^{2}\mu \left( x,\xi +2\pi n\right) =\frac{\sin ^{2}\left( x\right)
}{\pi ^{2}x^{2}}dx\otimes \sum_{n\in \mathbb{Z}}\delta _{\left( 2n+1\right)
\pi }\left( \xi \right)
\end{equation*}
and this is different from $\mu ^{\varphi }$, which is precisely:
\begin{equation*}
\mu ^{\varphi }\left( x,\xi \right) =\frac{\sin ^{2}\left( x/2\right) }{\pi
^{2}x^{2}}dx\otimes \sum_{n\in \mathbb{Z}}\delta _{\left( 2n+1\right) \pi
}\left( \xi \right) .
\end{equation*}

Finally, we investigate hypothesis (\ref{singular measures}). Now we set $%
\varphi :=\delta _{0}$.\bigskip

\textbf{3.}\emph{\ }\textbf{Necessity of condition (\ref{singular measures})}
\textbf{in Theorem \ref{Thm WM Suh uh}. }Define
\begin{equation*}
\widehat{v^{h}}\left( \xi\right) :=\mathbf{1}_{Q}\left( h\xi\right)
\sum_{n\in\mathbb{Z}}\mathbf{1}_{\left( -1,1\right) }\left( \xi-\left(
2n+1\right) \pi\right) .
\end{equation*}
Clearly, as in our first example, the periodization of the Wigner measure of
$\left( v^{h}\right) $ is:
\begin{equation*}
\sum_{k\in\mathbb{Z}^{d}}\mu\left( x,\xi+2\pi n\right) =\frac{\sin
^{2}\left( x/2\right) }{\pi^{2}x^{2}}dx\otimes\sum_{k\in\mathbb{Z}%
}\delta_{\left( 2n+1\right) \pi}\left( \xi\right) .
\end{equation*}
However, the sequence of discretizations $\left(
S_{\delta_{0}}^{h}v^{h}\right) $ has the following one:
\begin{equation*}
\mu^{\delta_{0}}\left( x,\xi\right) =\frac{\sin^{2}\left( x\right) }{%
\pi^{2}x^{2}}dx\otimes\sum_{n\in\mathbb{Z}}\delta_{\left( 2n+1\right) \pi
}\left( \xi\right) .
\end{equation*}

The proof of these counterexamples easily follows from (\ref{examp con}),
identity (\ref{relation meuh and meFuh}) and Lemma \ref{Lemma orth}.

\subsection{A Poisson summation formula and proof of Theorem \ref{Thm WM Suh
uh}}

The computation of the Fourier transform of $S_{\varphi }^{h}u$ is given by
the following identity:

\begin{lemma}
\label{Lemma FT Shu}Let $\varphi $ satisfy (\ref{BP}) and $u\in H^{-s}\left(
\mathbb{R}^{d}\right) $. Then the Fourier transform of $S_{\varphi }^{h}u$
is:
\begin{equation}
h^{d}\sum_{n\in \mathbb{Z}^{d}}S_{\varphi }^{h}u\left( n\right) e^{-ihn\cdot
\xi }=\sum_{n\in \mathbb{Z}^{d}}\overline{\widehat{\varphi }\left( h\xi
+2\pi n\right) }\widehat{u}\left( \xi +\frac{2\pi }{h}n\right) ,
\label{generalized poisson fmla}
\end{equation}
the convergence of the first series being in $L_{\text{loc}}^{2}\left(
\mathbb{R}^{d}\right) $ while the second takes place in $L_{\text{loc}%
}^{1}\left( \mathbb{R}^{d}\right) $.
\end{lemma}

\begin{proof}
Begin by noticing that $\overline{\widehat{\varphi }}\widehat{u}\in
L^{1}\left( \mathbb{R}^{d}\right) $ and thus
\begin{equation*}
\Pi ^{h}\left( \xi \right) :=\sum_{n\in \mathbb{Z}^{d}}\overline{\widehat{%
\varphi }\left( h\xi +2\pi n\right) }\widehat{u}\left( \xi +2\pi /hn\right)
\end{equation*}
is a well-defined $\left( 2\pi /h\right) \mathbb{Z}^{d}$-periodic $L_{\text{%
loc}}^{1}\left( \mathbb{R}^{d}\right) $ function, the series defining it
being absolutely convergent in $L_{\text{loc}}^{1}\left( \mathbb{R}%
^{d}\right) $. We can compute its Fourier coefficients:
\begin{align*}
\int_{\left[ -\pi /h,\pi /h\right) ^{d}}\Pi ^{h}\left( \xi \right)
e^{ihn\cdot \xi }\frac{h^{d}d\xi }{\left( 2\pi \right) ^{d}}=\sum_{k\in
\mathbb{Z}^{d}}\int_{Q}\overline{\widehat{\varphi }\left( \xi +2\pi k\right)
}\widehat{u}\left( \frac{\xi +2\pi k}{h}\right) e^{in\cdot \xi }\frac{d\xi }{%
\left( 2\pi \right) ^{d}} \\
=\int_{\mathbb{R}^{d}}\overline{\widehat{\varphi }\left( \xi \right) }%
\widehat{u}\left( \frac{\xi }{h}\right) e^{in\cdot \xi }\frac{d\xi }{\left(
2\pi \right) ^{d}} \\
=\int_{\mathbb{R}^{d}}\overline{h^{d}\widehat{\varphi }\left( h\xi \right)
e^{-ihn\cdot \xi }}\widehat{u}\left( \xi \right) \frac{d\xi }{\left( 2\pi
\right) ^{d}} \\
=\left\langle \overline{\varphi _{n}^{h}},u\right\rangle _{\mathcal{S}%
^{\prime }\times \mathcal{S}}=h^{d}S_{\varphi }^{h}u\left( n\right) .
\end{align*}
Lemma \ref{Lemma norm TphiUh} proves that $S_{\varphi }^{h}u$ is
square-summable and, consequently,
\begin{equation*}
\Pi ^{h}\left( \xi \right) =\sum_{n\in \mathbb{Z}^{d}}h^{d}S_{\varphi
}^{h}u\left( n\right) e^{-ihn\cdot \xi },
\end{equation*}
the sum being understood in the $L^{2}$-sense. This is precisely formula (%
\ref{generalized poisson fmla}).\medskip
\end{proof}

\begin{remark}
Identity (\ref{generalized poisson fmla}) may be viewed as a generalization
of \textbf{Poisson summation formula}. Taking as $\varphi $ the Dirac delta $%
\delta _{0}$, we obtain:
\begin{equation*}
h^{d}\sum_{n\in \mathbb{Z}^{d}}u\left( hn\right) e^{-ihn\cdot \xi
}=\sum_{n\in \mathbb{Z}^{d}}\widehat{u}\left( \xi +\frac{2\pi }{h}n\right) ,
\end{equation*}
for every $u\in H^{s}\left( \mathbb{R}^{d}\right) $ with $s>d/2$.
\end{remark}

\begin{proof}[Proof of Theorem \ref{Thm WM Suh uh}]
The proof will be done in two steps:\medskip

\textbf{Step 1: }We first establish the result for sequences such that $%
\widehat{\varphi }\left( \xi \right) \widehat{u_{k}}\left( \xi /h_{k}\right)
$ has support in a ball $B\left( 0;R\right) $ for every $k\in \mathbb{N}$.
We claim that the following formula holds:
\begin{equation*}
T_{\delta _{0}}^{h_{k}}S_{\varphi }^{h_{k}}u_{k}\left( x\right)
=\sum_{\left| n\right| \leq R+\pi \sqrt{d}}e^{-2\pi in\cdot x/h_{k}}%
\overline{\widehat{\varphi }}\left( h_{k}D_{x}\right) u_{k}\left( x\right) .
\end{equation*}
This is obtained by applying the inverse Fourier transform at both sides of
identity (\ref{generalized poisson fmla}), and remarking that only summands
satisfying $\left| n\right| \leq R+\pi \sqrt{d}$ must be considered because
of the condition on the support of $\widehat{\varphi }\widehat{u_{k}}\left(
\cdot /h_{k}\right) $. The Wigner measures of the functions
\begin{equation*}
e^{-2\pi in\cdot x/h_{k}}\overline{\widehat{\varphi }}\left(
h_{k}D_{x}\right) u_{k}\left( x\right)
\end{equation*}
are precisely (cf. Proposition \ref{Prop localization} and Remark \ref{Rmk
phihat in linf}):
\begin{equation*}
\left| \widehat{\varphi }\left( \xi +2\pi n\right) \right| ^{2}\mu \left(
x,\xi +2\pi n\right) .
\end{equation*}
By hypothesis, they are mutually singular so, by Lemma \ref{Lemma orth}, we
deduce that the measure $\mu ^{\varphi }$ obtained as the limit of $m^{h_{k}}%
\left[ T_{\delta _{0}}^{h_{k}}S_{\varphi }^{h_{k}}u_{k}\right] $ is given by
(\ref{relation meuh and MeSUh h=e}).\medskip

\textbf{Step 2: }We prove the result in the general case by taking advantage
of hypothesis (\ref{h-osc}). Let $\chi \in C_{c}^{\infty }\left( \mathbb{R}%
^{d}\right) $ be a cut-off function identically equal to one in the unit
ball $B\left( 0;1\right) .$ Denote by $S_{\varphi ,R}^{h_{k}}u_{k}$ the
truncation given by:
\begin{align*}
\widehat{S_{\varphi ,R}^{h_{k}}u_{k}}\left( \xi \right) :=\widehat{%
S_{\varphi }^{h_{k}}\chi \left( \frac{h_{k}D_{x}}{R}\right) u_{k}}\left( \xi
\right) \\
=\frac{1}{\left( h_{k}\right) ^{d}}\sum_{n\in \mathbb{Z}^{d}}\overline{%
\widehat{\varphi }\left( \xi +2\pi n\right) }\chi \left( \frac{\xi +2\pi n}{R%
}\right) \widehat{u_{k}}\left( \frac{\xi +2\pi n}{h_{k}}\right) ;
\end{align*}
Then, by the first step we have just proved, $M^{h_{k}}\left[ S_{\varphi
,R}^{h_{k}}u\right] $ converges to
\begin{equation}
\mu _{R}^{\varphi }\left( x,\xi \right) :=\sum_{n\in \mathbb{Z}^{d}}\left|
\widehat{\varphi }\left( \xi +2\pi n\right) \chi \left( \frac{\xi +2\pi n}{R}%
\right) \right| ^{2}\mu \left( x,\xi +2\pi n\right) .  \label{mesphi}
\end{equation}
We claim that (\ref{h-osc}) implies the following:
\begin{equation}
\limsup_{k\rightarrow \infty }\left\| S_{\varphi }^{h_{k}}u_{k}-S_{\varphi
,R}^{h_{k}}u_{k}\right\| _{L^{2}\left( h_{k}\mathbb{Z}^{d}\right)
}^{2}\rightarrow 0\text{\qquad as }R\rightarrow \infty \text{.}  \label{dif}
\end{equation}
It is sufficient to realize that
\begin{equation*}
S_{\varphi }^{h_{k}}u_{k}-S_{\varphi ,R}^{h_{k}}u_{k}=S_{\psi
_{R}}^{h_{k}}u_{k}
\end{equation*}
for $\widehat{\psi _{R}}:=\overline{\chi \left( \cdot /R\right) }\widehat{%
\varphi }$. The norm of $S_{\psi _{R}}^{h_{k}}$ is precisely (cf. Lemma \ref
{Lemma norm TphiUh})
\begin{equation*}
\limfunc{essup}\limits_{\xi \in Q}\sum_{n\in \mathbb{Z}^{d}}\left|
\left\langle \xi +2\pi n\right\rangle ^{s}\widehat{\varphi }\left( \xi +2\pi
n\right) \chi \left( \frac{\xi +2\pi n}{R}\right) \right| ^{2}
\end{equation*}
which tends to zero as $R\rightarrow 0$.

Lemma \ref{Lemma approx WM} then ensures that $\mu _{R}^{\varphi }$ weakly
converge to $\mu ^{\varphi }$. Identity (\ref{mesphi}) means that
\begin{equation*}
\int_{\mathbb{R}^{d}\times \mathbb{R}^{d}}a\left( x,\xi \right) d\mu
_{R}^{\varphi }\left( x,\xi \right) =\int_{\mathbb{R}^{d}\times \mathbb{R}%
^{d}}\sum_{k\in \mathbb{Z}^{d}}a\left( x,\xi +2\pi n\right) \left| \widehat{%
\varphi }\left( \xi \right) \right| ^{2}\left| \chi \left( \xi /R\right)
\right| ^{2}d\mu \left( x,\xi \right)
\end{equation*}
for every test function $a\in \mathcal{S}$. Passing to limits as $%
R\rightarrow \infty $ in the above identity we obtain the claimed result.

Notice that the same argument may be applied if, instead of condition (\ref
{h-osc}), we have that $\left( \left\langle h_{k}D_{x}\right\rangle
^{-s}u_{k}\right) $ is $h_{k}$-oscillatory. This is because (\ref{dif}) may
be estimated from above by
\begin{equation*}
\limsup_{k\rightarrow \infty }\left\| \left\langle h_{k}D_{x}\right\rangle
^{-s}\left( 1-\chi \left( \frac{h_{k}D_{x}}{R}\right) \right) u_{k}\right\|
_{L^{2}\left( \mathbb{R}^{d}\right) }^{2}\rightarrow 0\qquad \text{as }%
R\rightarrow \infty
\end{equation*}
because of Lemma \ref{Lemma norm TphiUh} and the $h_{k}$-oscillation
hypothesis.\medskip
\end{proof}

\section{\label{Sec defect}Computation of defect measures}

\subsection{Relations between defect and Wigner measures in the discrete
setting}

In this paragraph, we establish the analog of Proposition \ref{Prop e-oscill}
in the discrete setting. In particular, we present conditions that ensure
that the projection on the $x$-component of a Wigner measure may be obtained
as the limit of quadratic densities of the type:
\begin{equation*}
E^{h}\left[ U^{h}\right] \left( x\right) :=h^{d}\sum_{n\in \mathbb{Z}%
^{d}}\left| U_{n}^{h}\right| ^{2}\delta _{hn}\left( x\right) \text{.}
\end{equation*}

\begin{proposition}
\label{Prop mu and nu e=h}Let $\left( h_{k}\right) $ be a scale and $\left(
U^{h_{k}}\right) $ be an $h_{k}$-bounded sequence. Suppose that $\left(
M^{h_{k}}\left[ U^{h_{k}}\right] \right) $ converges to $\mu $ as $%
k\rightarrow \infty $. Then, for every $\phi \in C_{c}\left( \mathbb{R}%
^{d}\right) $,
\begin{equation}
\int_{\mathbb{R}^{d}\times Q}\phi \left( x\right) d\mu \left( x,\xi \right)
=\lim_{k\rightarrow \infty }\left( h_{k}\right) ^{d}\sum_{n\in \mathbb{Z}%
^{d}}\phi \left( h_{k}n\right) \left| U_{n}^{h_{k}}\right| ^{2}.
\label{marginal}
\end{equation}
If $\left( \varepsilon _{k}\right) $ is a scale such that $h_{k}\ll
\varepsilon _{k}$ and the transforms $M^{\varepsilon _{k}}\left[ U^{h_{k}}%
\right] $ converge to $\mu $, then (\ref{marginal}) holds provided $\left(
U^{h_{k}}\right) $ is $\mathbf{\varepsilon }_{k}$\textbf{-oscillatory, }%
i.e.:
\begin{equation}
\limsup_{k\rightarrow \infty }\left( h_{k}\right) ^{d}\int_{Q\setminus
B\left( 0;h_{k}/\varepsilon _{k}R\right) }\left| \widehat{U^{h_{k}}}\left(
\xi \right) \right| ^{2}d\xi \rightarrow 0\qquad \text{as }R\rightarrow
\infty .  \label{e-oscillatory D S}
\end{equation}
\end{proposition}

In view of Proposition \ref{Prop mu and nu e=h}, one could think that Wigner
measures at scales coarser than $h_{k}$ are unnecessary. However, as the
next result shows, if $\left( U^{h_{k}}\right) $ is $\varepsilon _{k}$%
-oscillatory for such a scale then the Wigner measure at scale $\left(
h_{k}\right) $ does not give any information about the oscillation effects.

\begin{proposition}
\label{Prop eoscill}Let $\left( h_{k}\right) $ and $\left( \varepsilon
_{k}\right) $ be scales such that $h_{k}\ll \varepsilon _{k}$. For every $%
\varepsilon _{k}$-oscillatory, $h_{k}$-bounded sequence $\left(
U^{h_{k}}\right) $ such that $M^{h_{k}}\left[ U^{h_{k}}\right]
\rightharpoonup \mu $ as $k\rightarrow \infty $ we have
\begin{equation*}
\mu \left( x,\xi \right) =\nu \left( x\right) \otimes \sum_{k\in \mathbb{Z}%
^{d}}\delta _{2\pi k}\left( \xi \right)
\end{equation*}
where $\nu $ is the weak limit in $\mathcal{M}_{+}\left( \mathbb{R}%
^{d}\right) $ of the measures $E^{h_{k}}\left[ U^{h_{k}}\right] $.
\end{proposition}

The Wigner measure also gathers the information on the densities $\left|
\mathcal{F}^{\varepsilon _{k}}U^{h_{k}}\left( \xi \right) \right| ^{2}$;
indeed, these converge to the projection on $\xi $ of the Wigner measure
provided that no energy is lost at infinity.

\begin{proposition}
\label{Prop mu Fe}Let $\left( h_{k}\right) $ and $\left( \varepsilon
_{k}\right) $ be scales such that $h_{k}/\varepsilon _{k}$ is bounded.
Suppose that $\left( U^{h_{k}}\right) $ is \textbf{compact at infinity}:
\begin{equation}
\limsup_{k\rightarrow \infty }\left( h_{k}\right) ^{d}\sum_{\left|
h_{k}n\right| >R}\left| U_{n}^{h_{k}}\right| ^{2}\rightarrow 0,\qquad \text{%
as }R\rightarrow \infty \text{,}  \label{compact at infinity}
\end{equation}
and that $M^{\varepsilon _{k}}\left[ U^{h_{k}}\right] \rightharpoonup \mu $
as $k\rightarrow \infty $. Then
\begin{equation*}
\int_{\mathbb{R}^{d}\times \mathbb{R}^{d}}\psi \left( \xi \right) d\mu
\left( x,\xi \right) =\lim_{k\rightarrow \infty }\int_{\mathbb{R}^{d}}\psi
\left( \xi \right) \left| \mathcal{F}^{\varepsilon _{k}}U^{h_{k}}\left( \xi
\right) \right| ^{2}d\xi
\end{equation*}
for every $\psi \in C_{c}\left( \mathbb{R}^{d}\right) $.
\end{proposition}

The proof of Propositions \ref{Prop mu and nu e=h} and \ref{Prop mu Fe}
requires the following preliminary result, which explains how the transform $%
M^{\varepsilon}\left[ U^{h}\right] $ of a discrete function $U^{h}$ can be
localized:

\begin{lemma}
\label{Lemma localization D}Let $U^{h}\in L^{2}\left( h\mathbb{Z}^{d}\right)
$ and $\varphi $, $\phi \in C_{c}^{\infty }\left( \mathbb{R}^{d}\right) $.
Then for every $a\in \mathcal{S}\left( \mathbb{R}^{d}\times \mathbb{R}%
^{d}\right) $ the following holds:
\begin{equation*}
\lim_{k\rightarrow \infty }\left| \left\langle M^{\varepsilon _{k}}\left[
U^{h_{k}}\right] ,\left| \phi \left( x\right) \right| ^{2}\varphi \left( \xi
\right) \right\rangle _{\mathcal{S}^{\prime }\times \mathcal{S}}-\left(
h_{k}\right) ^{d}\int_{\mathbb{R}^{d}}\left| \widehat{\phi U^{h_{k}}}\left(
\xi \right) \right| ^{2}\varphi \left( \frac{\varepsilon _{k}}{h_{k}}\xi
\right) \frac{d\xi }{\left( 2\pi \right) ^{d}}\right| =0.
\end{equation*}
\end{lemma}

\begin{proof}
First remark that, as a consequence of relation (\ref{relation MUh and mTUh}%
) and Lemma \ref{Lemma localization I} we have
\begin{equation}
\lim_{k\rightarrow \infty }\left| \left\langle M^{\varepsilon _{k}}\left[
U^{h_{k}}\right] ,\left| \phi \left( x\right) \right| ^{2}\varphi \left( \xi
\right) \right\rangle _{\mathcal{S}^{\prime }\times \mathcal{S}%
}-\left\langle M^{\varepsilon _{k}}\left[ \phi U^{h_{k}}\right] ,\psi \left(
x\right) \varphi \left( \xi \right) \right\rangle _{\mathcal{S}^{\prime
}\times \mathcal{S}}\right| =0  \label{loc en e-oscilante}
\end{equation}
for every test function $\psi \in C_{c}^{\infty }\left( \mathbb{R}%
^{d}\right) $ such that $\psi \left( x\right) =1$ for $x\in \limfunc{supp}%
\phi $. Now, (\ref{meu on a l2loc}.i) and (\ref{relation MUh and mTUh})
together with Plancherel's formula for the discrete Fourier transform yield:
\begin{align*}
\left\langle M^{\varepsilon _{k}}\left[ \phi U^{h_{k}}\right] ,\psi \left(
x\right) \varphi \left( \xi \right) \right\rangle _{\mathcal{S}^{\prime
}\times \mathcal{S}}& =\left\langle \phi \left( x\right) T_{\delta
_{0}}^{h_{k}}U^{h_{k}},\varphi \left( \varepsilon _{k}D_{x}\right) \phi
\left( x\right) T_{\delta _{0}}^{h_{k}}U^{h_{k}}\right\rangle _{\mathcal{S}%
^{\prime }\times \mathcal{S}} \\
& =\left( h_{k}\right) ^{2d}\int_{\mathbb{R}^{d}}\left| \widehat{\phi
U^{h_{k}}}\left( h_{k}\xi \right) \right| ^{2}\varphi \left( \varepsilon
_{k}\xi \right) \frac{d\xi }{\left( 2\pi \right) ^{d}}
\end{align*}
and the result follows.\bigskip
\end{proof}

\begin{proof}[Proof of Proposition \ref{Prop mu and nu e=h}]
Identity (\ref{marginal}) in the case $h_{k}=\varepsilon _{k}$ is a direct
consequence of the identity
\begin{equation*}
\int_{Q}M^{h}\left[ U\right] \left( x,\xi \right) d\xi =E^{h}\left[ U^{h}%
\right] \left( x\right)
\end{equation*}
and that, due to the $\Gamma $-periodicity in $\xi $ of $M^{h_{k}}\left[
U^{h_{k}}\right] $ and $\mu $, one has
\begin{equation*}
\lim_{k\rightarrow \infty }\int_{\mathbb{R}^{d}\times Q}\phi \left( x\right)
M^{h_{k}}\left[ U^{h_{k}}\right] \left( x,\xi \right) dxd\xi =\int_{\mathbb{R%
}^{d}\times Q}\phi \left( x\right) d\mu \left( x,\xi \right)
\end{equation*}
for every $\phi \in C_{c}^{\infty }\left( \mathbb{R}^{d}\right) $.

Next we analyze the case $h_{k}/\varepsilon _{k}\rightarrow 0$. Given
functions $\phi ,\chi \in C_{c}^{\infty }\left( \mathbb{R}^{d}\right) $,
using Lemma \ref{Lemma localization D} and periodization in the variable $%
\xi $ we find:
\begin{equation}
\int_{\mathbb{R}^{d}\times \mathbb{R}^{d}}\left| \phi \left( x\right)
\right| ^{2}\chi \left( \xi \right) d\mu \left( x,\xi \right)
=\lim_{k\rightarrow \infty }\left( h_{k}\right) ^{d}\int_{Q}\left| \widehat{%
\phi U^{h_{k}}}\left( \xi \right) \right| ^{2}\sum_{n\in \mathbb{Z}^{d}}\chi
\left( \frac{\varepsilon _{k}}{h_{k}}\left( \xi +2\pi n\right) \right) \frac{%
d\xi }{\left( 2\pi \right) ^{d}}.  \label{lim phi ki}
\end{equation}
Choose a function $\chi \in C_{c}^{\infty }\left( \mathbb{R}^{d}\right) $
such that
\begin{equation*}
\left.
\begin{array}{l}
\chi \left( \xi \right) =1\text{ for }\left| \xi \right| \leq 1,\medskip \\
\chi \left( \xi \right) =0\text{ for }\left| \xi \right| \geq 2,\medskip \\
0\leq \chi \left( \xi \right) \leq 1\text{ for }\xi \in \mathbb{R}^{d},
\end{array}
\right.
\end{equation*}
and set $\chi _{R}\left( \xi \right) :=\chi \left( \xi /R\right) $ for every
$R>0$. With such a test function and $h_{k}/\varepsilon _{k}<\pi /R$ we have
$\chi _{R}\left( \frac{\varepsilon _{k}}{h_{k}}\left( \xi +2\pi n\right)
\right) =\chi _{R}\left( \frac{\varepsilon _{k}}{h_{k}}\xi \right) \leq 1$
for every $\xi \in Q$. Then, taking this into account in (\ref{lim phi ki})
and using Plancherel's formula, the following is obtained:
\begin{equation*}
\lim_{k\rightarrow \infty }\left| \left( h_{k}\right) ^{d}\int_{Q}\left|
\widehat{\phi U^{h_{k}}}\left( \xi \right) \right| ^{2}\frac{d\xi }{\left(
2\pi \right) ^{d}}-\int_{\mathbb{R}^{d}\times \mathbb{R}^{d}}\left| \phi
\left( x\right) \right| ^{2}\chi _{R}\left( \xi \right) d\mu \left( x,\xi
\right) \right| \leq M\left( R\right) ,
\end{equation*}
where
\begin{equation*}
M\left( R\right) :=\limsup_{k\rightarrow \infty }\left( h_{k}\right)
^{d}\int_{\mathbb{R}^{d}}\left( 1-\chi _{R}\left( \frac{\varepsilon _{k}}{%
h_{k}}\left( \xi +2\pi n\right) \right) \right) \left| \widehat{\phi
U^{h_{k}}}\left( \xi \right) \right| ^{2}\frac{d\xi }{\left( 2\pi \right)
^{d}}.
\end{equation*}
Identity (\ref{marginal}) is obtained by letting $R$ tend to $\infty $,
noticing that (\ref{e-oscillatory D S}) implies that $M\left( R\right)
\rightarrow 0$ as $R\rightarrow \infty $.\bigskip
\end{proof}

\begin{proof}[Proof of Proposition \ref{Prop eoscill}]
Since $\left( U^{h_{k}}\right) _{k\in \mathbb{N}}$ is $\varepsilon _{k}$%
-oscillatory, we have
\begin{equation*}
\limsup_{k\rightarrow \infty }\int_{Q\setminus B\left( 0;\delta \right)
}\left| \widehat{\phi U^{h_{k}}}\left( \xi \right) \right| ^{2}d\xi =0
\end{equation*}
for every $\delta >0$ and $\phi \in C_{c}^{\infty }\left( \mathbb{R}%
^{d}\right) $. Using Lemma \ref{Lemma localization D} below, we obtain, for
every $\varphi \in C_{c}^{\infty }\left( Q\setminus \left\{ 0\right\}
\right) $,
\begin{equation*}
0=\lim_{k\rightarrow \infty }\int_{Q}\varphi \left( \xi \right) \left|
\widehat{\phi U^{h_{k}}}\left( \xi \right) \right| ^{2}\frac{d\xi }{\left(
2\pi \right) ^{d}}=\int_{\mathbb{R}^{d}\times Q}\left| \phi \left( x\right)
\right| ^{2}\varphi \left( \xi \right) d\mu \left( x,\xi \right) .
\end{equation*}
In particular, $\mu $ is concentrated on the set $\mathbb{R}^{d}\times
\left\{ 0\right\} $. Since $\mu \left( \cdot \times Q\right) =\nu \left(
x\right) $ by Proposition \ref{Prop mu and nu e=h} we find that, because of
the periodicity, $\mu \left( \mathbb{R}^{d}\times \cdot \right) =\nu \left(
\mathbb{R}^{d}\right) \sum_{k\in \mathbb{Z}^{d}}\delta _{2\pi k}\left( \xi
\right) $; and this restricts $\mu $ to be equal to $\nu \otimes \sum_{k\in
\mathbb{Z}^{d}}\delta _{2\pi k}$.\bigskip
\end{proof}

\begin{proof}[Proof of Proposition \ref{Prop mu Fe}]
Since the densities $\left| \mathcal{F}^{\varepsilon _{k}}U^{h_{k}}\right|
^{2}$ are uniformly bounded in $L^{1}\left( \mathbb{R}^{d}\right) $ (and
consequently in $\mathcal{M}\left( \mathbb{R}^{d}\right) $) it suffices to
prove the result for test functions $\psi \in C_{c}^{\infty }\left( \mathbb{R%
}^{d}\right) $. Let $\chi $ and $\chi _{R}$ be defined as in the proof of
Proposition \ref{Prop mu and nu e=h}. Because of Lemma \ref{Lemma
localization D} the following holds for every $\psi \in \mathcal{S}\left(
\mathbb{R}^{d}\right) $:
\begin{equation*}
\int_{\mathbb{R}^{d}\times \mathbb{R}^{d}}\left| \chi _{R}\left( x\right)
\right| ^{2}\psi \left( \xi \right) d\mu \left( x,\xi \right)
=\lim_{k\rightarrow \infty }\int_{\mathbb{R}^{d}}\psi \left( \xi \right)
\left| \mathcal{F}^{\varepsilon _{k}}\chi _{R}U^{h_{k}}\left( \xi \right)
\right| ^{2}d\xi .
\end{equation*}
Since $\left| \chi _{R}\left( x\right) \right| ^{2}\rightarrow 1$ as $%
R\rightarrow \infty $ for every $x\in \mathbb{R}^{d}$ we only have to show
that
\begin{equation*}
\lim_{R\rightarrow \infty }\int_{\mathbb{R}^{d}\times \mathbb{R}^{d}}\left|
\chi _{R}\left( x\right) \right| ^{2}\psi \left( \xi \right) d\mu \left(
x,\xi \right) =\lim_{k\rightarrow \infty }\int_{\mathbb{R}^{d}}\psi \left(
\xi \right) \left| \mathcal{F}^{\varepsilon _{k}}U^{h_{k}}\left( \xi \right)
\right| ^{2}d\xi .
\end{equation*}
This appears as a consequence of the identity
\begin{align*}
\int_{\mathbb{R}^{d}}\psi \left( \xi \right) \left( \left| \mathcal{F}%
^{\varepsilon _{k}}U^{h_{k}}\left( \xi \right) \right| ^{2}-\left| \mathcal{F%
}^{\varepsilon _{k}}\chi _{R}U^{h_{k}}\left( \xi \right) \right|
^{2}\right) d\xi \\
=\int_{\mathbb{R}^{d}}\psi \left( \xi \right) \left[ \mathcal{F}%
^{\varepsilon _{k}}\left( U^{h_{k}}-\chi _{R}U^{h_{k}}\right) \left( \xi
\right) \right] \overline{\mathcal{F}^{\varepsilon _{k}}U^{h_{k}}\left( \xi
\right) }d\xi \\
+\int_{\mathbb{R}^{d}}\psi \left( \xi \right) \mathcal{F}^{\varepsilon
_{k}}\chi _{R}U^{h_{k}}\left( \xi \right) \overline{\left[ \mathcal{F}%
^{\varepsilon _{k}}\left( U^{h_{k}}-\chi _{R}U^{h_{k}}\right) \left( \xi
\right) \right] }d\xi
\end{align*}
which implies
\begin{equation*}
\limsup_{k\rightarrow \infty }\left| \int_{\mathbb{R}^{d}}\psi \left( \xi
\right) \left( \left| \mathcal{F}^{\varepsilon _{k}}U^{h_{k}}\left( \xi
\right) \right| ^{2}-\left| \mathcal{F}^{\varepsilon _{k}}\chi
_{R}U^{h_{k}}\left( \xi \right) \right| ^{2}\right) d\xi \right| \leq
C_{\psi }\limsup_{k\rightarrow \infty }\left\| U^{h_{k}}-\chi
_{R}U^{h_{k}}\right\| _{L^{2}\left( h\mathbb{Z}^{d}\right) }^{2}.
\end{equation*}
Since the $U^{h_{k}}$ are compact at infinity, the second term in the above
estimate tends to zero as $R$ tends to infinity and thus:
\begin{equation*}
\left| \limsup_{k\rightarrow \infty }\int_{\mathbb{R}^{d}}\psi \left( \xi
\right) \left| \mathcal{F}^{\varepsilon _{k}}U^{h_{k}}\left( \xi \right)
\right| ^{2}-\int_{\mathbb{R}^{d}\times \mathbb{R}^{d}}\left| \chi
_{R}\left( x\right) \right| ^{2}\psi \left( \xi \right) d\mu \left( x,\xi
\right) \right| \rightarrow 0\qquad \text{as }R\rightarrow \infty .
\end{equation*}
One easily deduces from this that the measures $\left| \mathcal{F}%
^{\varepsilon _{k}}U^{h_{k}}\left( \xi \right) \right| ^{2}d\xi $ converge
in $\mathcal{M}_{+}\left( \mathbb{R}^{d}\right) $ to the measure $\int_{%
\mathbb{R}^{d}}\mu \left( dx,\cdot \right) $ as claimed.\bigskip
\end{proof}

\subsection{Defect measures of reconstructed sequences}

Let $\left( U^{h_{k}}\right) $ be $h_{k}$-bounded and $\varphi \in
H^{s}\left( \mathbb{R}^{d}\right) $ some profile satisfying (\ref{BP}). As a
consequence of Lemma \ref{Lemma norm TphiUh}, the sequence of densities
\begin{equation*}
\left| \left\langle h_{k}D_{x}\right\rangle ^{s}T_{\varphi
}^{h_{k}}U^{h_{k}}\right| ^{2}
\end{equation*}
is uniformly bounded in $L^{1}\left( \mathbb{R}^{d}\right) $. Hence, Helly's
compactness Theorem ensures that, extracting a subsequence if necessary,
there exists a measure $\nu _{\varphi }\in \mathcal{M}_{+}\left( \mathbb{R}%
^{d}\right) $ such that
\begin{equation*}
\lim_{k\rightarrow \infty }\int_{\mathbb{R}^{d}}\phi \left( x\right) \left|
\left\langle h_{k}D_{x}\right\rangle ^{s}T_{\varphi }^{h_{k}}U^{h_{k}}\left(
x\right) \right| ^{2}dx=\int_{\mathbb{R}^{d}\times Q}\phi \left( x\right)
d\nu _{\varphi }\left( x\right) .
\end{equation*}
The main issue addressed in this section is that of clarifying how $\nu
_{\varphi }$ depends of the sequence $\left( U^{h_{k}}\right) $ and the
profile $\varphi $. We shall see that a formula relating $\nu _{\varphi }$
and the limit of $E^{h_{k}}\left[ U^{h_{k}}\right] $ does not exist in
general. However such a formula may be established in terms of the Wigner
measure of $\left( U^{h_{k}}\right) $.

Suppose that $M^{h_{k}}\left[ U^{h_{k}}\right] $ converges to $\mu $. Then,
Theorem \ref{Thm WM TUh Uh} may be applied to obtain that, provided $\mu
\left( \mathbb{R}^{d}\times \overline{D_{\widehat{\varphi }}}\right) =0$,
one has
\begin{equation*}
m^{h_{k}}\left[ T_{\varphi }^{h_{k}}U^{h_{k}}\right] \rightharpoonup \left|
\widehat{\varphi }\left( \xi \right) \right| ^{2}\mu \left( x,\xi \right) .
\end{equation*}
In general, we are only able to ensure (see Proposition 1.7 in \cite{G-M-M-P}%
):
\begin{equation*}
\nu _{\varphi }\left( x\right) \geq \int_{\mathbb{R}^{d}}\left| \left\langle
\xi \right\rangle ^{s}\widehat{\varphi }\left( \xi \right) \right| ^{2}\mu
\left( x,d\xi \right) ,
\end{equation*}
but equality holds whenever $\left( \left\langle h_{k}D_{x}\right\rangle
^{s}T_{\varphi }^{h_{k}}U^{h_{k}}\right) $ is $h_{k}$-oscillatory. However,
this is not immediate since there exist profiles $\varphi $ for which $%
\left( \left\langle h_{k}D_{x}\right\rangle ^{s}T_{\varphi
}^{h_{k}}U^{h_{k}}\right) $ may fail to be $h_{k}$-oscillatory for some $%
\left( U^{h_{k}}\right) $ (an example is provided at the end of this
section). Nevertheless,

\begin{proposition}
$\left( \left\langle h_{k}D_{x}\right\rangle ^{s}T_{\varphi
}^{h_{k}}U^{h_{k}}\right) $ is $h_{k}$-oscillatory whenever (\ref{h-osc})
holds.
\end{proposition}

This immediately follows from:

\begin{lemma}
\label{Lemma hosc}$\left( \left\langle h_{k}D_{x}\right\rangle
^{s}T_{\varphi }^{h_{k}}U^{h_{k}}\right) $ is $h_{k}$-oscillatory if and
only if
\begin{equation*}
\limsup_{k\rightarrow \infty }\left( h_{k}\right) ^{d}\int_{Q}\sigma
_{\varphi }^{R}\left( \xi \right) \left| \widehat{U^{h_{k}}}\left( \xi
\right) \right| ^{2}d\xi \rightarrow 0\qquad \text{as }R\rightarrow \infty ,
\end{equation*}
where
\begin{equation*}
\sigma _{\varphi }^{R}\left( \xi \right) :=\sum_{\left| n\right| \geq
R}\left| \left\langle \xi +2\pi n\right\rangle ^{s}\widehat{\varphi }\left(
\xi +2\pi n\right) \right| ^{2}.
\end{equation*}
\end{lemma}

\begin{proof}
Start by noticing that:
\begin{equation*}
\int_{\left| \xi \right| \geq R/h_{k}}\left| \widehat{\left\langle
h_{k}D_{x}\right\rangle ^{s}T_{\varphi }^{h_{k}}U^{h_{k}}}\left( \xi \right)
\right| ^{2}d\xi =\int_{\left| \xi \right| \geq R}\left( h_{k}\right)
^{d}\left| \left\langle \xi \right\rangle ^{s}\widehat{\varphi }\left( \xi
\right) \widehat{U^{h_{k}}}\left( \xi \right) \right| ^{2}d\xi .
\end{equation*}
Periodizing in $\xi $ we get:
\begin{equation*}
\int_{Q}\sigma _{\varphi }^{R+\sqrt{d}\pi }\left( \xi \right) \left|
\widehat{U^{h_{k}}}\left( \xi \right) \right| ^{2}d\xi \leq \int_{\left| \xi
\right| \geq R}\left| \left\langle \xi \right\rangle ^{s}\widehat{\varphi }%
\left( \xi \right) \widehat{U^{h_{k}}}\left( \xi \right) \right| ^{2}d\xi
\leq \int_{Q}\sigma _{\varphi }^{R-\sqrt{d}\pi }\left( \xi \right) \left|
\widehat{U^{h_{k}}}\left( \xi \right) \right| ^{2}d\xi ,
\end{equation*}
and the claim follows.\bigskip
\end{proof}

Hence, $h_{k}$-oscillation is obtained if $\varphi $ decays at infinity at
an uniform rate. For more general $\varphi $, it is still possible to obtain
sufficient conditions; however, these depend on the particular sequence of
discrete functions to be reconstructed.

\begin{proposition}
Suppose
\begin{equation}
\begin{array}{ll}
\text{i)} & \mu \left( \mathbb{R}^{d}\times \overline{D_{\tau _{\left\langle
D_{x}\right\rangle ^{s}\varphi }}}\right) =0,\medskip \\
\text{ii)} & \left( U^{h_{k}}\right) \text{ is compact at infinity.}
\end{array}
\label{c2}
\end{equation}
Then $\left( \left\langle h_{k}D_{x}\right\rangle ^{s}T_{\varphi
}^{h_{k}}U^{h_{k}}\right) $ is $h_{k}$-oscillatory.
\end{proposition}

\begin{proof}
For the sake of simplicity, we prove the result for $s=0$; the proof in the
general case being identical. Taking into account the periodicity of the
densities involved, Proposition \ref{Prop mu Fe} ensures that
\begin{equation}
\lim_{k\rightarrow \infty }\left( h_{k}\right) ^{d}\int_{Q}\psi \left( \xi
\right) \left| \widehat{U^{h_{k}}}\left( \xi \right) \right| ^{2}d\xi =\int_{%
\mathbb{R}^{d}\times Q}\psi \left( \xi \right) d\mu \left( x,\xi \right)
\label{weak conv}
\end{equation}
for every $\psi \in C_{c}\left( \mathbb{R}^{d}\right) $ and the claim
follows. Since $\mu \left( \mathbb{R}^{d}\times \overline{D_{\tau _{\varphi
}}}\right) =0$, necessarily $\mu \left( \mathbb{R}^{d}\times \overline{%
D_{\sigma _{\varphi }^{R}}}\right) \ $is null for every $R>0$. From
classical results on weak convergence of measures, one deduces that relation
(\ref{weak conv}) also holds for $\psi =\sigma _{\varphi }^{R}$. Hence, by
the dominated convergence Theorem,
\begin{equation*}
\lim_{R\rightarrow \infty }\lim_{k\rightarrow \infty }\left( h_{k}\right)
^{d}\int_{Q}\sigma _{\varphi }^{R}\left( \xi \right) \left| \widehat{%
U^{h_{k}}}\left( \xi \right) \right| ^{2}d\xi =\lim_{R\rightarrow \infty
}\int_{\mathbb{R}^{d}\times Q}\sigma _{\varphi }^{R}\left( \xi \right) d\mu
\left( x,\xi \right) =0,
\end{equation*}
and the result follows.\medskip
\end{proof}

What we have proved so, combined with Proposition \ref{Prop e-oscill} far is
gathered in the next proposition.

\begin{proposition}
\label{Proposition defect}Suppose at least one of (\ref{h-osc}) or (\ref{c2}%
) is satisfied and that $\mu \left( \mathbb{R}^{d}\times \overline{D_{%
\widehat{\varphi }}}\right) =0$. Then
\begin{equation}
\lim_{k\rightarrow \infty }\int_{\mathbb{R}^{d}}\phi \left( x\right) \left|
\left\langle h_{k}D_{x}\right\rangle ^{s}T_{\varphi
}^{h_{k}}U^{h_{k}}\right| ^{2}dx=\int_{\mathbb{R}^{d}\times \mathbb{R}%
^{d}}\phi \left( x\right) \left| \left\langle \xi \right\rangle ^{s}\widehat{%
\varphi }\left( \xi \right) \right| ^{2}d\mu \left( x,\xi \right)
\label{identidad md}
\end{equation}
for every $\phi \in C_{c}\left( \mathbb{R}^{d}\right) $.
\end{proposition}

As anticipated above, (\ref{identidad md}) shows that the knowledge of the
weak limit of the measures $E^{h_{k}}\left[ U^{h_{k}}\right] $ and the
profile $\varphi $ are not enough, in general, to reconstruct the weak limit
of the densities $\left| \left\langle h_{k}D_{x}\right\rangle ^{s}T_{\varphi
}^{h_{k}}U^{h_{k}}\right| ^{2}dx$. However, when the sequence of discrete
functions under consideration is $\varepsilon _{k}$-oscillatory for some
scale coarser than the reconstruction step $h_{k}$, there does exist a
formula that relates both limits:

\begin{corollary}
\label{Cor defect}Let $\left( U^{h_{k}}\right) $be an $h_{k}$-bounded, $%
\varepsilon _{k}$-oscillatory sequence such that $\left( E^{h_{k}}\left[
U^{h_{k}}\right] \right) $ weakly converges to a measure $\nu $. Suppose
moreover that $\widehat{\varphi }$ is continuous at $\Gamma $ and that any
of (\ref{h-osc}) or (\ref{c2}) is satisfied. Then the densities $\left|
\left\langle h_{k}D_{x}\right\rangle ^{s}T_{\varphi
}^{h_{k}}U^{h_{k}}\right| ^{2}$ weakly converge to the measure
\begin{equation*}
\nu _{\varphi }\left( x\right) =\left( \sum_{n\in \mathbb{Z}^{d}}\left|
\left\langle 2\pi n\right\rangle ^{s}\widehat{\varphi }\left( 2\pi n\right)
\right| ^{2}\right) \nu \left( x\right) .
\end{equation*}
\end{corollary}

\begin{proof}
Using Proposition \ref{Prop eoscill}, we find that any Wigner measure at
scale $h_{k}$ of $\left( U^{h_{k}}\right) $ equals
\begin{equation*}
\mu \left( x,\xi \right) =\nu \left( x\right) \otimes \sum_{n\in \mathbb{Z}%
^{d}}\delta _{2\pi n}\left( \xi \right) .
\end{equation*}
Since $\mu \left( \mathbb{R}^{d}\times \overline{D_{\widehat{\varphi }}}%
\right) =0$, Proposition \ref{Proposition defect} is applicable and gives:
\begin{equation*}
\left| \left\langle h_{k}D_{x}\right\rangle ^{s}T_{\varphi
}^{h_{k}}U^{h_{k}}\right| ^{2}dx\rightharpoonup \tau _{\left\langle
D_{x}\right\rangle ^{s}\varphi }\left( 0\right) \nu \left( x\right) \qquad
\text{as }k\rightarrow \infty \text{,}
\end{equation*}
as claimed.\bigskip
\end{proof}

Remark that condition (\ref{c2}.i) reduces in this setting to the
requirement that $\tau _{\left\langle D_{x}\right\rangle ^{s}\varphi }$ is
continuous at $\xi =0$.

\subsection{\label{ce osc}A counterexample to $h$-oscillation}

Here we exhibit a function $\varphi \in L^{2}\left( \mathbb{R}\right) $
satisfying (\ref{BP}) with $\widehat{\varphi }$ is continuous but
\begin{equation*}
\left\| \sigma _{\varphi }^{R}\right\| _{L^{\infty }\left( Q\right)
}=1\qquad \text{for every }R>0.
\end{equation*}
With such a profile, we show that there exist a sequence of discrete
functions $\left( U^{h}\right) $ such that $\left( T_{\varphi
}^{h}U^{h}\right) $ is not $h$-oscillatory

To construct $\varphi$, define $t_{n}:=e^{-n}$ for $n=0,1,2...$ and let $%
\psi_{n}$ be the piecewise linear function given for $n\geq1$ by
\begin{equation*}
\psi_{n}\left( t\right) =\left\{
\begin{array}{ll}
\dfrac{t-t_{n+1}}{t_{n}-t_{n+1}} & \text{if }t\in\left( t_{n+1},t_{n}\right)
,\medskip \\
\dfrac{t-t_{n-1}}{t_{n}-t_{n-1}} & \text{if }t\in\left( t_{n},t_{n-1}\right)
,\medskip \\
0 & \text{otherwise.}
\end{array}
\right.
\end{equation*}
Clearly $\sum_{n=1}^{\infty}\psi_{n}\left( t\right) =1$ for $t\in\left(
0,t_{1}\right) $, and the sum vanishes for $t\leq0$. Defining
\begin{equation*}
\widehat{\varphi}\left( \xi\right) :=\sqrt{\sum_{n=1}^{\infty}\psi
_{n}\left( \xi-2\pi n\right) }
\end{equation*}
we obtain $\varphi\in L^{2}\left( \mathbb{R}^{d}\right) $, $\widehat
{\varphi}\in C\left( \mathbb{R}^{d}\right) $ and $\tau_{\varphi}\left(
\xi\right) =\sum_{n=1}^{\infty}\psi_{n}\left( \xi\right) $.

Moreover
\begin{equation*}
\sigma_{\varphi}^{n}\left( \xi\right) =\left\{
\begin{array}{ll}
1 & \text{if }\xi\in\left( 0,x_{n+1}\right) ,\medskip \\
0 & \text{if }\xi\leq0\text{.}
\end{array}
\right.
\end{equation*}
Thus $\left\| \sigma_{\varphi}^{R}\right\| _{L^{\infty}\left( Q\right) }=1$
for every $R>0$.

If we chose discrete functions $U^{h}\in L^{2}\left( h\mathbb{Z}\right) $
such that
\begin{equation*}
\widehat{U^{h}}\left( \xi \right) =h^{-1}\sum_{n\in \mathbb{Z}}\mathbf{1}%
_{\left( 0,h\right) }\left( \xi +2\pi n\right)
\end{equation*}
then for the $\varphi $ above constructed we obtain
\begin{equation*}
\lim_{h\rightarrow 0}\int_{Q}\sigma _{\varphi }^{R}\left( \xi \right)
h\left| \widehat{U^{h}}\left( \xi \right) \right| ^{2}d\xi =1\qquad \text{%
for every }R>0\text{.}
\end{equation*}
This proves that $\left( T_{\varphi }^{h}U^{h}\right) $ is not $h$%
-oscillatory.

\section{\label{Sec hlesse}High frequency analysis: $h\ll \protect\varepsilon
$}

Here we shall investigate the structure of Wigner measures at scales $\left(
\varepsilon _{k}\right) $ asymptotically coarser than the
sampling/reconstruction rate $\left( h_{k}\right) $.

In the next two Theorems, we suppose that $\varphi $ satisfies (\ref{BP})%
\textbf{\ }and $\widehat{\varphi }$ \emph{is continuous in a neighborhood of}
$\xi =0$. Moreover, $\left( h_{k}\right) $ and $\left( \varepsilon
_{k}\right) $ will be scales such that $h_{k}\ll \varepsilon _{k}$.

\begin{theorem}
\label{Thm WM TUh Uh eh}Suppose that $\left( U^{h_{k}}\right) $ is $h_{k}$%
-bounded and $M^{\varepsilon _{k}}\left[ U^{h_{k}}\right] $ converges to the
Wigner measure $\mu $. Then $m^{\varepsilon _{k}}\left[ T_{\varphi
}^{h_{k}}U^{h_{k}}\right] $ converges to a measure $\mu _{\varphi }$ given
by:
\begin{equation}
\mu _{\varphi }\left( x,\xi \right) =\left| \widehat{\varphi }\left(
0\right) \right| ^{2}\mu \left( x,\xi \right) .
\label{muphi and mu he to zero}
\end{equation}
\end{theorem}

The proof of this result is completely analogous to that of Theorem \ref{Thm
WM TUh Uh}.

Concerning the sampling operators, the situation is much similar:

\begin{theorem}
\label{Thm WM Suh uh eh}Let $\left( u_{k}\right) $ be a sequence in $%
H^{-s}\left( \mathbb{R}^{d}\right) $ such that $\left( \left\langle
h_{k}D_{x}\right\rangle ^{-s}u_{k}\right) $ is bounded in $L^{2}\left(
\mathbb{R}^{d}\right) $ and $\varepsilon _{k}$-oscillatory.\medskip

i) Then $\left( S_{\varphi }^{h_{k}}u_{k}\right) $ is $\varepsilon _{k}$%
-oscillatory.\medskip

ii) Suppose moreover that $m^{\varepsilon _{k}}\left[ u_{k}\right] $
converges to a Wigner measure $\mu $. Then $M^{\varepsilon _{k}}\left[
S_{\varphi }^{h_{k}}u_{k}\right] $ converges to the Wigner measure $\mu
^{\varphi }$ given by:
\begin{equation}
\mu ^{\varphi }=\left| \widehat{\varphi }\left( 0\right) \right| ^{2}\mu .
\label{relation meuh and MeSUh h<e}
\end{equation}
\end{theorem}

\begin{proof}
To prove the first part of the Theorem, begin by noticing that, by the
Cauchy-Schwarz inequality and Lemma \ref{Lemma FT Shu}, for almost every $%
\xi \in \mathbb{R}^{d}$:
\begin{align*}
\left| \widehat{S_{\varphi }^{h_{k}}u_{k}}\left( \xi \right) \right|
^{2}=\left| \frac{1}{\left( h_{k}\right) ^{d}}\sum_{n\in \mathbb{Z}^{d}}%
\overline{\widehat{\varphi }\left( \xi +2\pi n\right) }\widehat{u_{k}}\left(
\frac{\xi +2\pi n}{h_{k}}\right) \right| ^{2} \\
\leq \frac{\left\| \tau _{\left\langle D_{x}\right\rangle ^{s}\varphi
}\right\| _{L^{\infty }\left( Q\right) }}{\left( h_{k}\right) ^{2d}}%
\sum_{n\in \mathbb{Z}^{d}}\left| \left\langle \xi +2\pi n\right\rangle ^{-s}%
\widehat{u_{k}}\left( \frac{\xi +2\pi n}{h_{k}}\right) \right| ^{2}.
\end{align*}
Thus
\begin{align*}
\int_{Q\setminus B\left( 0;h_{k}/\varepsilon _{k}R\right) }\left(
h_{k}\right) ^{d}\left| \widehat{S_{\varphi }^{h_{k}}u_{k}}\left(
\xi \right) \right| ^{2}d\xi \leq \\\left\| \tau _{\left\langle
D_{x}\right\rangle ^{s}\varphi }\right\| _{L^{\infty }\left(
Q\right) }\int_{Q\setminus B\left( 0;R/\varepsilon _{k}\right)
}\sum_{n\in \mathbb{Z}^{d}}\left| \left\langle h_{k}\xi +2\pi
n\right\rangle ^{-s}\widehat{u_{k}}\left( \xi +2\pi n\right)
\right| ^{2}d\xi \\
\leq \left\| \tau _{\left\langle D_{x}\right\rangle ^{s}\varphi }\right\|
_{L^{\infty }\left( Q\right) }\int_{\mathbb{R}^{d}\setminus B\left(
0;R/\varepsilon _{k}\right) }\left| \left\langle h_{k}\xi \right\rangle ^{-s}%
\widehat{u_{k}}\left( \xi \right) \right| ^{2}d\xi ,
\end{align*}
and this clearly proves that $\left( S_{\varphi }^{h_{k}}u_{k}\right) $ is $%
\varepsilon _{k}$-oscillating as soon as $\left( \left\langle
h_{k}D_{x}\right\rangle ^{-s}u_{k}\right) $ is.

The proof of identity (\ref{relation meuh and MeSUh h<e}) is essentially
identical to that of Theorem \ref{Thm WM Suh uh}. A completely analogous
argument to that used in the step 2 of that proof allows us to consider only
sequences such that $\widehat{\varphi }\left( h_{k}/\varepsilon _{k}\cdot
\right) \widehat{u_{k}}\left( \cdot /\varepsilon _{k}\right) $ is supported
in a ball $B\left( 0;R\right) $. This hypothesis together with Lemma \ref
{Lemma FT Shu} implies that, for $h_{k}/\varepsilon _{k}$ small enough,
\begin{equation*}
\left( h_{k}\right) ^{d}\widehat{S_{\varphi }^{h_{k}}u_{k}}\left( \frac{h_{k}%
}{\varepsilon _{k}}\xi \right) =\overline{\widehat{\varphi }\left( \frac{%
h_{k}}{\varepsilon _{k}}\xi \right) }\widehat{u_{k}}\left( \xi /\varepsilon
_{k}\right) ,
\end{equation*}
that is, only one summand is involved. Then the result follows from
Proposition \ref{Prop localization} exactly as in the proof of Theorem \ref
{Thm WM Suh uh}.\medskip
\end{proof}

We conclude with a simple remark:

\begin{corollary}
Under the assumptions and notations of Theorems \ref{Thm WM TUh Uh eh} and
\ref{Thm WM Suh uh eh}:\medskip

i) If $\varphi $ has zero mean (i.e. $\widehat{\varphi }\left( 0\right) =0$)
then the Wigner measure at scale $\left( \varepsilon _{k}\right) \ $of any
sequence $\left( T_{\varphi }^{h_{k}}U^{h_{k}}\right) $ or $\left(
S_{\varphi }^{h_{k}}u_{k}\right) $ vanishes identically. In particular, this
is the case if $\varphi $ is a wavelet.\footnote{%
See, for instance, \cite{H-W}, Proposition 2.1.}\medskip

ii) $\mu _{\varphi }=\mu ^{\varphi }=\mu $ always holds for profiles such
that $\left| \widehat{\varphi }\left( 0\right) \right| =1$.
\end{corollary}

\section{\label{Sec Sr}Wigner measures of Sampled/Reconstructed sequences}

Now we are able to describe Wigner measures of sequences of the form $%
T_{\psi }^{h}S_{\varphi }^{h}u$. In its full generality, our result requires
several compatibility hypothesis, that we describe below. First of all,
\begin{equation}
\begin{array}{ll}
\text{i)} & \psi \text{ and }\varphi \text{ satisfy (\ref{BP}) with
exponents }s^{\prime }\text{ and }s\text{ respectively.}\bigskip \\
\text{ii)} & \varphi \text{ satisifes (\ref{h-osc}).}
\end{array}
\label{hip phipsi}
\end{equation}
The admissible sequences will be assumed to be such that:
\begin{equation}
u_{k}\in H^{-s}\left( \mathbb{R}^{d}\right) \text{ and }\left( \left\langle
h_{k}D_{x}\right\rangle ^{-s}u_{k}\right) \text{ is bounded in }L^{2}\left(
\mathbb{R}^{d}\right) ,  \label{bdd}
\end{equation}
and their Wigner measures must satisfy the following compatibility
conditions for some precise representatives of $\widehat{\psi }$ and $%
\widehat{\varphi }$:
\begin{equation}
\begin{array}{ll}
\text{i)} & \mu \text{ fulfills (\ref{ND})}.\medskip \\
\text{ii)} & \dint_{\mathbb{R}^{d}\times \mathbb{R}^{d}}\mathbf{1}_{%
\overline{D_{\widehat{\psi }}}}\left( \xi +2\pi n\right) \left| \widehat{%
\varphi }\left( \xi \right) \right| ^{2}d\mu \left( x,\xi \right) =0,\quad
n\in \mathbb{Z}^{d}.\medskip \\
\text{iii)} & \mu \text{ satisfies (\ref{singular measures}).}
\end{array}
\label{hip mu}
\end{equation}
Combining Theorems \ref{Thm WM TUh Uh} and \ref{Thm WM Suh uh} we obtain:

\begin{theorem}
\label{Thm main}Let $\psi $ and $\varphi $ be functions satisfying (\ref{hip
phipsi}); let $\left( h_{k}\right) $ be a scale and $\left( u_{k}\right) $
be a sequence satisfying (\ref{bdd}). Suppose moreover that $m^{h_{k}}\left[
u_{k}\right] $ converges to a Wigner measure $\mu $ that satisfies (\ref{hip
mu}).

Then $m^{h_{k}}\left[ T_{\psi }^{h_{k}}S_{\varphi }^{h_{k}}u_{k}\right] $
converges to the measure $\mu _{\varphi ,\psi }$ given by:
\begin{equation}
\int_{\mathbb{R}^{d}\times \mathbb{R}^{d}}a\left( x,\xi \right) d\mu
_{\varphi ,\psi }\left( x,\xi \right) =\int_{\mathbb{R}^{d}\times \mathbb{R}%
^{d}}\sum_{n\in \mathbb{Z}^{d}}a\left( x,\xi +2\pi n\right) \left| \widehat{%
\psi }\left( \xi +2\pi n\right) \right| ^{2}\left| \widehat{\varphi }\left(
\xi \right) \right| ^{2}d\mu \left( x,\xi \right)  \label{rep muphipsi}
\end{equation}
for every $a\in C_{c}\left( \mathbb{R}^{d}\times \mathbb{R}^{d}\right) $.
\end{theorem}

\begin{proof}
Hypothesis (\ref{hip mu}.ii) expresses that the closure of the set of
discontinuity points of $\widehat{\psi }$, is a null set for the Wigner
measure of $S_{\varphi }^{h_{k}}u_{k}$, $\sum_{k\in \mathbb{Z}^{d}}\left|
\widehat{\varphi }\left( \xi +2\pi n\right) \right| ^{2}\mu \left( x,\xi
+2\pi n\right) $. Hence, Theorem \ref{Thm WM TUh Uh} is applicable and we
conclude that the distributions $m^{h_{k}}\left[ T_{\psi }^{h_{k}}S_{\varphi
}^{h_{k}}u_{k}\right] $ converge to the measure
\begin{equation*}
\left| \widehat{\psi }\left( \xi \right) \right| ^{2}\sum_{n\in \mathbb{Z}%
^{d}}\left| \widehat{\varphi }\left( \xi +2\pi n\right) \right| ^{2}\mu
\left( x,\xi +2\pi n\right) .
\end{equation*}
Since $\left| \widehat{\psi }\left( \xi \right) \right| ^{2}$ is integrable
with respect to the finite measure $\left| \widehat{\varphi }\right| ^{2}\mu
$ (this is again due to (\ref{hip mu}.ii)), its periodization is integrable
as well and formula (\ref{rep muphipsi}) follows.\medskip
\end{proof}

\begin{remark}
i) When $\widehat{\psi }$ and $\widehat{\varphi }$ verify (\ref{strong shyp}%
) hypotheses (\ref{hip phipsi}), (\ref{hip mu}.i) and (\ref{hip mu}.ii) are
immediately satisfied.\smallskip

ii) (\ref{hip phipsi}.ii) may be replaced by the requirement that $\left(
\left\langle h_{k}D_{x}\right\rangle ^{-s}u_{k}\right) $ is $h_{k}$%
-oscillatory.
\end{remark}

From formula (\ref{rep muphipsi}) one sees at once that, taking $\psi
=\varphi =\delta _{0}$ one has that $\mu _{\varphi ,\psi }$ is the
periodization in $\xi $ of the Wigner measure $\mu $. Hence, $\mu _{\varphi
,\psi }$ coincides with the limit of the Wigner series corresponding to $%
\left( u_{k}\right) $.

When $\widehat{\psi }$ and $\left| \widehat{\varphi }\right| ^{2}\mu $
vanish off $Q$ it is easy to check that formula (\ref{rep muphipsi}) takes
the simple form:
\begin{equation*}
\mu _{\varphi ,\psi }\left( x,\xi \right) =\left| \widehat{\psi }\left( \xi
\right) \right| ^{2}\left| \widehat{\varphi }\left( \xi \right) \right|
^{2}\mu \left( x,\xi \right) .
\end{equation*}
It is also clear that, as soon as $\left| \widehat{\varphi }\left( \xi
\right) \right| ^{2}\mu \left( x,\xi \right) $ is not null outside $Q$, the
measures $\mu _{\varphi ,\psi }$ and $\mu $ will in general differ.

Concerning defect measures, combining Proposition \ref{Proposition defect}
and the previous theorem, we obtain:

\begin{theorem}
\label{Thm defect}Under the notations of Theorem \ref{Thm main} the
following holds: if
\begin{equation*}
\left| \left\langle h_{k}D_{x}\right\rangle ^{s^{\prime }}T_{\psi
}^{h_{k}}S_{\varphi }^{h_{k}}u_{k}\right| ^{2}dx\text{ weakly converges to a
measure }\nu _{\varphi ,\psi }
\end{equation*}
and $\psi $ verifies (\ref{h-osc}) then:
\begin{equation}
\nu _{\varphi ,\psi }\left( x\right) =\int_{\mathbb{R}_{\xi }^{d}}\sum_{n\in
\mathbb{Z}^{d}}\left| \left\langle \xi +2\pi n\right\rangle ^{s^{\prime }}%
\widehat{\psi }\left( \xi +2\pi n\right) \right| ^{2}\left| \widehat{\varphi
}\left( \xi \right) \right| ^{2}\mu \left( x,d\xi \right) .  \label{rep def}
\end{equation}
\end{theorem}

\begin{remark}
i) The conclusion of the Theorem still holds if condition ``$\psi $
satisfies (\ref{h-osc})'' is replaced by (\ref{c2}).\smallskip

ii) Theorems \ref{Thm comp Wm} and \ref{Thm Shan} follow immediately from
Theorems \ref{Thm main} and \ref{Thm defect}.
\end{remark}

With formula (\ref{rep def}) at our disposal, we are now able to answer, in
a quite general way, the questions \textbf{A}-\textbf{D} addressed in the
introduction. Of course, the answer to \textbf{A }is negative, since, in
general, $\mu $ is not trivial in its $\xi $ component; concerning the
problem of filtering, we immediately get the necessary and sufficient
condition:
\begin{equation*}
c_{\varphi ,\psi }\left( \xi \right) :=\left| \widehat{\varphi }\left( \xi
\right) \right| ^{2}\sum_{n\in \mathbb{Z}^{d}}\left| \left\langle \xi +2\pi
n\right\rangle ^{s^{\prime }}\widehat{\psi }\left( \xi +2\pi n\right)
\right| ^{2}=0\text{\qquad for }\mu \text{-a.e. }\xi \in \mathbb{R}^{d}\text{%
.}
\end{equation*}
Analogously, $c_{\varphi ,\psi }\left( \xi \right) =1$ for $\mu $-a.e. $\xi
\in \mathbb{R}^{d}$ characterizes the profiles that give $\nu _{\varphi
,\psi }=\nu $. To answer \textbf{D}, we must, of course, assume that $%
\widehat{\varphi }$ and $\tau _{\left\langle D_{x}\right\rangle ^{s^{\prime
}}\widehat{\psi }}$ are continuous (which, as we know, is the case if (\ref
{strong shyp}) holds). In that case, we have and equality $\nu _{\varphi
,\psi }=\nu $ for every admissible sequence if and only if
\begin{equation*}
\left| \widehat{\varphi }\left( \xi \right) \right| ^{2}=\frac{1}{\tau
_{\left\langle D_{x}\right\rangle ^{s^{\prime }}\psi }\left( \xi \right) }%
\qquad \text{for every }\xi \in \mathbb{R}^{d}\text{ with }\tau
_{\left\langle D_{x}\right\rangle ^{s^{\prime }}\psi }\left( \xi \right)
\neq 0\text{.}
\end{equation*}
The sampling profile $\varphi $, cannot be an $L^{2}\left( \mathbb{R}%
^{d}\right) $ function, since $\left| \widehat{\varphi }\right| ^{2}$ is
necessarily periodic. When $\varphi =\delta _{0}$ and $\psi $ generates an
orthonormal basis in the sense of Lemma \ref{Lemma Riesz ort} we always have
$\nu _{\varphi ,\psi }=\nu $. If $\psi $ merely generates a Riesz basis, $%
A\nu \leq \nu _{\varphi ,\psi }\leq B\nu $ holds instead.

The above results may be used to compute Wigner measures of the \textbf{%
orthogonal projections} $P_{\psi }^{h_{k}}u_{k}$ of a given sequence $\left(
u_{k}\right) $ on the shift-invariant space defined by the range of $T_{\psi
}^{h_{k}}$. As we have seen in Lemma \ref{Lemma ort}, $P_{\psi }^{h}$ may be
written as the composition of $T_{\psi }^{h}$ with $S_{\varphi
}^{h}\left\langle hD_{x}\right\rangle ^{s}$ for a sampling profile $\varphi
:=\widetilde{\left\langle D_{x}\right\rangle ^{s}\psi }$. Hence, Theorem \ref
{Thm main} gives:

\begin{corollary}
\label{cor proj} For $\psi $ satisfying (\ref{strong shyp}) and $\left(
u_{k}\right) $ such that (\ref{bdd}) and (\ref{singular measures}) holds,
the defect measures of the sequence $\left( P_{\psi }^{h_{k}}u_{k}\right) $
is given by:
\begin{equation*}
\nu _{P_{\psi }}\left( x\right) =\int_{\mathbb{R}^{d}}\frac{\mathbf{1}_{\psi
}\left( \xi \right) }{\tau _{\left\langle D_{x}\right\rangle ^{s}\psi
}\left( \xi \right) }\left| \left\langle \xi \right\rangle ^{s}\widehat{\psi
}\left( \xi \right) \right| ^{2}\mu \left( x,d\xi \right) ,
\end{equation*}
where $\mathbf{1}_{\psi }\left( \xi \right) $ denotes the characteristic
function of the set of $\xi \in \mathbb{R}^{d}$ such that $\tau
_{\left\langle D_{x}\right\rangle ^{s}\psi }\left( \xi \right) \neq 0$.
\end{corollary}

In particular, when $\psi $ gives rise to an orthonormal family, we obtain
the simple formula (cf. Lemma \ref{Lemma Riesz ort}):
\begin{equation*}
\nu _{P_{\psi }}\left( x\right) =\int_{\mathbb{R}^{d}}\left| \left\langle
\xi \right\rangle ^{s}\widehat{\psi }\left( \xi \right) \right| ^{2}\mu
\left( x,d\xi \right) .
\end{equation*}

To conclude, we shall see how the above results may be refined when the
sequence $\left( \left\langle h_{k}D_{x}\right\rangle ^{-s}u_{k}\right) $ is
assumed to be $\varepsilon _{k}$-oscillatory at some scale $h_{k}\ll
\varepsilon _{k}$. The assumptions of $\varphi $ and $\psi $ are weaker:
\begin{equation}
\begin{array}{ll}
\text{i)} & \psi \text{ and }\varphi \text{ satisfy (\ref{BP}) with
exponents }s^{\prime }\text{ and }s\text{ respectively}.\medskip \\
\text{ii)} & \widehat{\psi }\text{, }\widehat{\varphi }\text{ are continuous
in a neighborhood of }\xi =0.\text{ }\medskip \\
\text{iii)} & \tau _{\left\langle D_{x}\right\rangle ^{s^{\prime }}\psi }%
\text{ is continuous at }\xi =0.
\end{array}
\label{hip phipsi w}
\end{equation}
Theorems \ref{Thm WM TUh Uh eh}, \ref{Thm WM Suh uh eh} and Corollary \ref
{Cor defect} give then:

\begin{theorem}
\label{Thm main2}Let $\psi $ and $\varphi $ be functions satisfying (\ref
{hip phipsi w}); let $\left( h_{k}\right) $, $\left( \varepsilon _{k}\right)
$ be scales with $h_{k}\ll \varepsilon _{k}$ and let $\left( u_{k}\right) $
be a sequence such that (\ref{bdd}) holds and $\left( \left\langle
h_{k}D_{x}\right\rangle ^{-s}u_{k}\right) $ is $\varepsilon _{k}$%
-oscillatory. Suppose moreover that $m^{\varepsilon _{k}}\left[ u_{k}\right]
$ converges to a Wigner measure $\mu $.

Then $m^{\varepsilon _{k}}\left[ T_{\psi }^{h_{k}}S_{\varphi }^{h_{k}}u_{k}%
\right] $ converges to the measure $\mu _{\varphi ,\psi }$ given by:
\begin{equation*}
\mu _{\varphi ,\psi }\left( x,\xi \right) =\left| \widehat{\psi }\left(
0\right) \right| ^{2}\left| \widehat{\varphi }\left( 0\right) \right|
^{2}\mu \left( x,\xi \right) .
\end{equation*}
Moreover, if $\left| \left\langle h_{k}D_{x}\right\rangle ^{s^{\prime
}}T_{\psi }^{h_{k}}S_{\varphi }^{h_{k}}u_{k}\right| ^{2}dx$ weakly converges
to a measure $\nu _{\varphi ,\psi }$ then:
\begin{equation*}
\nu _{\varphi ,\psi }\left( x\right) =\sum_{n\in \mathbb{Z}^{d}}\left|
\left\langle 2\pi n\right\rangle ^{s^{\prime }}\widehat{\psi }\left( 2\pi
n\right) \right| ^{2}\left| \widehat{\varphi }\left( 0\right) \right|
^{2}\nu \left( x\right) ,
\end{equation*}
where $\nu $ is the weak limit of the densities $\left| \left\langle
h_{k}D_{x}\right\rangle ^{-s}u_{k}\right| ^{2}dx$.
\end{theorem}

Hence, when a sequence possesses a characteristic oscillation scale $\left(
\varepsilon _{k}\right) $ (that is the meaning of the $\varepsilon _{k}$%
-oscillation condition), choosing a sampling/ reconstruction rate
$\left( h_{k}\right) $ asymptotically finer than $\left(
\varepsilon _{k}\right) $ allows to completely capture its
oscillation/concentration behavior (modulo a constant that only
depends on $\psi $ and $\varphi $).

Filtering in that case can only be achieved by means of a sampling profile $%
\varphi $ with zero mean ($\widehat{\varphi }\left( 0\right) =0$) or a
reconstruction profile such that $\widehat{\psi }$ vanishes at $\Gamma $.

\section{\label{Sec Ap}Tools from the theory of Wigner measures}

The main tools from the theory of Wigner measures used in this article are
Propositions \ref{Prop definition wm l2loc} and \ref{Prop localization}
below. The first of these is an extension of Theorem \ref{Thm definition wm}
to bounded sequences in Sobolev spaces:

\begin{proposition}
\label{Prop definition wm l2loc}Let $\left( \varepsilon _{k}\right) $ be a
scale and $\left( u_{k}\right) $ be a sequence of functions in $H^{-s}\left(
\mathbb{R}^{d}\right) $ for some $s\geq 0$ satisfying:
\begin{equation}
\left\| \left\langle \varepsilon _{k}D_{x}\right\rangle ^{-s}u_{k}\right\|
_{L^{2}\left( \mathbb{R}^{d}\right) }\text{ are uniformly bounded in }k\text{%
.}  \label{uk is Hsh bounded}
\end{equation}

Then the sequence of distributions $\left( m^{\varepsilon _{k}}\left[ u_{k}%
\right] \right) $ is uniformly bounded in $\mathcal{S}^{\prime }$. Moreover,
any of its weakly converging subsequences tends to a positive measure.
\end{proposition}

As we have done so far, a measure $\mu\in\mathcal{M}_{+}\left( \mathbb{R}%
^{d}\times\mathbb{R}^{d}\right) $ will be called the \textbf{Wigner measure
at scale }$\left( \varepsilon_{k}\right) $ of a sequence $\left(
u_{k}\right) $ (satisfying the hypotheses of Proposition \ref{Prop
definition wm l2loc}) provided $m^{\varepsilon_{k}}\left[ u_{k}\right]
\rightharpoonup \mu$ in $\mathcal{S}^{\prime}$ as $k\rightarrow\infty$.

\begin{remark}
i) When $s>0$, condition (\ref{uk is Hsh bounded}) is stronger than just
requiring that $\left( u_{k}\right) $ is bounded in $H^{-s}\left( \mathbb{R}%
^{d}\right) $.\medskip

ii) Let $\left( h_{k}\right) $ be a scale such that $h_{k}\ll \varepsilon
_{k}$. If $\left\| \left\langle h_{k}D_{x}\right\rangle ^{-s}u_{k}\right\|
_{L^{2}\left( \mathbb{R}^{d}\right) }\leq C$ for every $k\in \mathbb{N}$
then $\left\| \left\langle \varepsilon _{k}D_{x}\right\rangle
^{-s}u_{k}\right\| _{L^{2}\left( \mathbb{R}^{d}\right) }$ is uniformly
bounded as well.

iii) The same result holds if $m^{\varepsilon }\left[ \cdot \right] $ is
replaced by the Wigner transform (\ref{definition WT}).
\end{remark}

The second main result of this section a localization formula for Wigner
measures which is used several times in this article:

\begin{proposition}
\label{Prop localization}Let $\left( \varepsilon _{k}\right) $, $\left(
h_{k}\right) $ be scales and let $\left( u_{k}\right) $ be a sequence in $%
H^{-s}\left( \mathbb{R}^{d}\right) $, $s\geq 0$, satisfying (\ref{uk is Hsh
bounded}). Suppose that $\phi $ is a Borel function such that $\phi \in
L^{\infty }\left( \mathbb{R}^{d};\left\langle \xi \right\rangle ^{r}\right) $
for some $r\in \mathbb{R}$. If $m^{\varepsilon _{k}}\left[ u_{k}\right] $
converges to $\mu $ then $m^{\varepsilon _{k}}\left[ \phi \left(
h_{k}D_{x}\right) u_{k}\right] $ converges to a Wigner measure $\mu _{\phi }$
which has the following properties:\medskip

i) If $h_{k}=\varepsilon _{k}$ and $\mu \left( \mathbb{R}^{d}\times
\overline{D_{\phi }}\right) =0$, $D_{\phi }$ being the set of points where $%
\phi $ is not continuous, then
\begin{equation*}
\mu _{\phi }\left( x,\xi \right) =\left| \phi \left( \xi \right) \right|
^{2}\mu \left( x,\xi \right) .
\end{equation*}

ii) If $h_{k}\ll \varepsilon _{k}$ and $\phi $ is continuous in a
neighborhood of $\xi =0$ then
\begin{equation*}
\mu _{\phi }=\left| \phi \left( 0\right) \right| ^{2}\mu .
\end{equation*}
\end{proposition}

When applied to $\phi\left( \xi\right) :=\left\langle \xi\right\rangle ^{s}$%
, this result gives:

\begin{remark}
\label{Rmk finite}Let $\left( \varepsilon _{k}\right) $, $\left(
h_{k}\right) $ and $\left( u_{k}\right) $ be as in Proposition \ref{Prop
localization}. Suppose $m^{\varepsilon _{k}}\left[ u_{k}\right] $ converges
to $\mu $. Then $m^{\varepsilon _{k}}\left[ \left\langle
h_{k}D_{x}\right\rangle ^{-s}u_{k}\right] $ converges to the measure $\mu
_{s}$ given by:
\begin{equation*}
\begin{array}{ll}
\mu _{s}\left( x,\xi \right) =\left\langle \xi \right\rangle ^{-2s}\mu
\left( x,\xi \right) , & \text{if }h_{k}=\varepsilon _{k},\medskip \\
\mu _{s}=\mu , & \text{if }h_{k}\ll \varepsilon _{k}\text{.}
\end{array}
\end{equation*}
In particular (cf. Theorem \ref{Thm definition wm}), $\left\langle \xi
\right\rangle ^{-2s}\mu $ (resp. $\mu $) is a finite measure when $%
h_{k}=\varepsilon _{k}$ (resp. $h_{k}\ll \varepsilon _{k}$).
\end{remark}

For the convenience of the reader, we give detailed proofs of both results;
they follow the ideas present in the existing literature on the subject (
\cite{Ge91c, Li-Pau, Ge-Lei, G-M-M-P}). Proposition \ref{Prop definition wm
l2loc} will be proved in paragraph \ref{sbs bdd}. We shall essentially show
that truncation of the high frequencies of a sequence satisfying (\ref{uk is
Hsh bounded}) implies $\xi$-variable localization of the corresponding $%
m^{\varepsilon}\left[ \cdot\right] $. Then we conclude by applying Theorem
\ref{Thm definition wm} to the localized sequence.

Proposition \ref{Prop localization} is proved in paragraph \ref{sbs ploc};
to conclude this section, we describe two results useful for the computation
of Wigner measures (Lemmas \ref{Lemma approx WM} and \ref{Lemma orth}).

\subsection{First properties of $m^{\protect\varepsilon}\left[ u\right] $}

We begin by discussing three alternative ways of computing $m^{\varepsilon }%
\left[ u\right] $ that may be used when $u$ is merely a tempered
distribution. First remark that, given a $u\in\mathcal{S}^{\prime}\left(
\mathbb{R}^{d}\right) $, it makes sense to consider the distribution $%
m^{\varepsilon}\left[ u\right] $ given by (\ref{meu}), since the Fourier
transform of $u$ is well-defined. Actually $m^{\varepsilon}\left[ u\right]
\in\mathcal{S}^{\prime}$.

\textbf{1. }The action of $m^{\varepsilon}\left[ u\right] $ on a test
function $a\in\mathcal{S}$ is given by any of the formulas (see \cite{Ge96}%
):
\begin{equation}
\left\langle m^{\varepsilon}\left[ u\right] ,a\right\rangle _{\mathcal{S}%
^{\prime}\times\mathcal{S}}=\left\{
\begin{array}{ll}
\left\langle \overline{u},a\left( x,\varepsilon D_{x}\right) u\right\rangle
_{\mathcal{S}^{\prime}\left( \mathbb{R}^{d}\right) \times\mathcal{S}\left(
\mathbb{R}^{d}\right) },\bigskip & \qquad\text{(i)} \\
\dint _{\mathbb{R}^{d}}\dint _{\mathbb{R}^{d}}\dfrac{1}{\varepsilon^{d}}%
k_{a}\left( x,\dfrac{x-p}{\varepsilon}\right) u\left( p\right) \overline{%
u\left( x\right) }dpdx. & \qquad\text{(ii)}
\end{array}
\right.  \label{meu on a l2loc}
\end{equation}
where $a\left( x,\varepsilon D_{x}\right) $ is the \textbf{semiclassical
pseudodifferential operator }of symbol $a$:
\begin{equation}
a\left( x,\varepsilon D_{x}\right) u\left( x\right) =\int_{\mathbb{R}%
^{d}}a\left( x,\varepsilon\xi\right) \widehat{u}\left( \xi\right)
e^{ix\cdot\xi}\frac{d\xi}{\left( 2\pi\right) ^{d}},  \label{definition opea}
\end{equation}
and the kernel $k_{a}\left( x,p\right) $ is the inverse Fourier transform of
$a$ with respect to $\xi$:
\begin{equation*}
k_{a}\left( x,p\right) :=\int_{\mathbb{R}^{d}}a\left( x,\xi\right)
e^{ip\cdot\xi}\frac{d\xi}{\left( 2\pi\right) ^{d}}.
\end{equation*}
Formula (\ref{meu on a l2loc}.i) makes sense because the operator $a\left(
x,\varepsilon D_{x}\right) $ maps continuously $\mathcal{S}^{\prime}\left(
\mathbb{R}^{d}\right) $ into $\mathcal{S}\left( \mathbb{R}^{d}\right) $
whenever $a\in\mathcal{S}$ (see, for instance, \cite{Mart}). The integral in
(\ref{meu on a l2loc}.ii) must, of course, be understood in distributional
sense.

\textbf{2.} The distribution $m^{\varepsilon}\left[ u\right] $ may be
computed through the rescaled Fourier transform
\begin{equation}
\mathcal{F}^{\varepsilon}u\left( \xi\right) :=\frac{1}{\left(
2\pi\varepsilon\right) ^{d/2}}\widehat{u}\left( \frac{\xi}{\varepsilon }%
\right) ,  \label{definition rescaled FT}
\end{equation}
using the identity:
\begin{equation}
m^{\varepsilon}\left[ u\right] \left( x,\xi\right) =\overline
{m^{\varepsilon}\left[ \mathcal{F}^{\varepsilon}u\right] \left(
\xi,-x\right) }.  \label{relation meuh and meFuh}
\end{equation}
This follows from a direct computation from the definition (\ref{meu}).

\textbf{3. }Now we present two localization formulas:

\begin{lemma}
\label{Lemma localization I}Let $u\in \mathcal{S}^{\prime }\left( \mathbb{R}%
^{d}\right) $, $\phi \in C^{\infty }\left( \mathbb{R}^{d};\left\langle
x\right\rangle ^{r}\right) $ for some $r\in \mathbb{R}$ and $a\in \mathcal{S}
$. Then there exists $r_{1}^{\sigma },r_{2}^{\sigma }\in \mathcal{S}$ such
that:
\begin{equation*}
\left.
\begin{array}{l}
\left\langle m^{\varepsilon }\left[ \phi u\right] ,a\right\rangle _{\mathcal{%
S}^{\prime }\times \mathcal{S}}=\left\langle \left| \phi \left( x\right)
\right| ^{2}m^{\varepsilon }\left[ u\right] ,a\right\rangle _{\mathcal{S}%
^{\prime }\times \mathcal{S}}+\varepsilon \left\langle m^{\varepsilon }\left[
u\right] ,r_{1}^{\varepsilon }\right\rangle _{\mathcal{S}^{\prime }\mathcal{%
\times S}},\bigskip \\
\left\langle m^{\varepsilon }\left[ \phi \left( hD_{x}\right) u\right]
,a\right\rangle _{\mathcal{S}^{\prime }\times \mathcal{S}}=\left\langle
\left| \phi \left( \dfrac{h}{\varepsilon }\xi \right) \right|
^{2}m^{\varepsilon }\left[ u\right] ,a\right\rangle _{\mathcal{S}^{\prime
}\times \mathcal{S}}+\dfrac{h}{\varepsilon }\left\langle m^{\varepsilon }%
\left[ u\right] ,r_{2}^{h/\varepsilon }\right\rangle _{\mathcal{S}^{\prime }%
\mathcal{\times S}}.
\end{array}
\right.
\end{equation*}
Moreover, the test functions $r_{1}^{\sigma }$, $r_{2}^{\sigma }$ are
uniformly bounded in $\mathcal{S}$ for $0<\sigma \leq 1$.
\end{lemma}

This holds as a consequence of standard results on symbolic calculus for
semiclassical pseudodifferential operators; see for instance \cite{Mart}.
Remark that Proposition \ref{Prop localization} is not a consequence of this
result, since the multiplier $\phi\left( hD_{x}\right) $ there may have a
non-smooth symbol.

\subsection{\label{sbs bdd}Boundedness of the transforms $m^{\protect%
\varepsilon }\left[ u\right] $}

The next lemmas are used to establish the boundedness in $\mathcal{S}%
^{\prime }$ of the sequence $\left( m^{\varepsilon_{k}}\left[ u_{k}\right]
\right) $ provided $\left( u_{k}\right) $ satisfies the hypotheses of
Proposition \ref{Prop definition wm l2loc}.

\begin{lemma}
For every $u\in L^{2}\left( \mathbb{R}^{d};\left\langle x\right\rangle
^{r}\right) $ and $a\in \mathcal{S}$ the following estimate holds:
\begin{equation*}
\left| \left\langle m^{\varepsilon }\left[ u\right] ,a\right\rangle _{%
\mathcal{S}^{\prime }\times \mathcal{S}}\right| \leq \left\| u\right\|
_{L^{2}\left( \mathbb{R}^{d};\left\langle x\right\rangle ^{r}\right)
}^{2}\int_{\mathbb{R}^{d}}\sup_{x\in \mathbb{R}^{d}}\left| k_{a}\left(
x,p\right) \left\langle x-\varepsilon p\right\rangle ^{-r/2}\left\langle
x\right\rangle ^{-r/2}\right| dp,
\end{equation*}
\end{lemma}

\begin{proof}
Use formula (\ref{meu on a l2loc}.ii) to write
\begin{equation*}
\left\langle m^{\varepsilon }\left[ u\right] ,a\right\rangle _{\mathcal{S}%
^{\prime }\times \mathcal{S}}=\dint_{\mathbb{R}^{d}}\dint_{\mathbb{R}%
^{d}}k_{a}\left( x,p\right) u\left( x-\varepsilon p\right) \overline{u\left(
x\right) }dpdx,
\end{equation*}
noticing that this integral makes sense as $k_{a}\in \mathcal{S}$. Multiply
and divide the integrand above by $\left\langle x-\varepsilon p\right\rangle
^{r/2}\left\langle x\right\rangle ^{r/2}$ to obtain, by H\"{o}lder's
inequality,
\begin{equation*}
\left| \left\langle m^{\varepsilon }\left[ u\right] ,a\right\rangle _{%
\mathcal{S}^{\prime }\times \mathcal{S}}\right| \leq \int_{\mathbb{R}%
^{d}}\sup_{x\in \mathbb{R}^{d}}\left| k_{a}\left( x,p\right) \left\langle
x-\varepsilon p\right\rangle ^{-r/2}\left\langle x\right\rangle
^{-r/2}\right| \int_{\mathbb{R}^{d}}\left| u_{r}\left( x-\varepsilon
p\right) \overline{u_{r}\left( x\right) }\right| dxdp,
\end{equation*}
where we have set $u_{r}\left( x\right) :=\left\langle x\right\rangle
^{r/2}u\left( x\right) $. The conclusion follows from another application of
H\"{o}lder's inequality.\medskip
\end{proof}

If $u\in H^{-s}\left( \mathbb{R}^{d}\right) $ then $\mathcal{F}%
^{\varepsilon}u\in L^{2}\left( \mathbb{R}^{d};\left\langle \xi\right\rangle
^{-2s}\right) $. Clearly,
\begin{equation}
\left\| \left\langle \varepsilon D_{x}\right\rangle ^{-s}u\right\|
_{L^{2}\left( \mathbb{R}^{d}\right) }^{2}=\left\| \mathcal{F}^{\varepsilon
}u\right\| _{L^{2}\left( \mathbb{R}^{d};\left\langle \xi\right\rangle
^{-2s}\right) }^{2}.  \label{norm identity}
\end{equation}
Thus, taking identity (\ref{relation meuh and meFuh}) into account, we
obtain using the preceding lemma:
\begin{equation}
\left| \left\langle m^{\varepsilon}\left[ u\right] ,a\right\rangle _{%
\mathcal{S}^{\prime}\times\mathcal{S}}\right| \leq\left\| \left\langle
\varepsilon D_{x}\right\rangle ^{-s}u\right\| _{L^{2}\left( \mathbb{R}%
^{d}\right) }^{2}\int_{\mathbb{R}^{d}}\sup_{\xi\in\mathbb{R}^{d}}\left|
\widehat{a}\left( q,\xi\right) \left\langle \xi+\varepsilon q\right\rangle
^{s}\left\langle \xi\right\rangle ^{s}\right| \frac{dq}{\left( 2\pi\right)
^{d}},  \label{est exact}
\end{equation}
where $\widehat{a}\left( q,\xi\right) $ denotes the Fourier transform in $x$
of the function $a\left( x,\xi\right) $.

\begin{lemma}
\label{Lemma meuh bdd Hs}For every $s\geq 0$ there exists a constant $%
C_{s,d}>0$ such that
\begin{equation}
\left| \left\langle m^{\varepsilon }\left[ u\right] ,a\right\rangle _{%
\mathcal{S}^{\prime }\times \mathcal{S}}\right| \leq C_{s,d}\left\|
\left\langle \varepsilon D_{x}\right\rangle ^{-s}u\right\| _{L^{2}\left(
\mathbb{R}^{d}\right) }^{2}\int_{\mathbb{R}^{d}}\sup_{\xi \in \mathbb{R}%
^{d}}\left| \widehat{a}\left( q,\xi \right) \left\langle \xi \right\rangle
^{2s}\right| \left\langle \varepsilon q\right\rangle ^{s}dq,
\label{estimate meu on a Hs}
\end{equation}
holds for every $u\in H^{-s}\left( \mathbb{R}^{d}\right) $ and every $a\in
\mathcal{S}$.
\end{lemma}

\begin{proof}
This is obtained through the simple inequality $\left\langle \xi
+q\right\rangle ^{s}\leq C_{s,d}\left\langle \xi \right\rangle
^{s}\left\langle q\right\rangle ^{s}$, which holds when $s\geq 0$.
\end{proof}

Notice that whenever $a\in\mathcal{S}$, the integrals $\int_{\mathbb{R}%
^{d}}\sup_{\xi\in\mathbb{R}^{d}}\left| \widehat{a}\left( q,\xi\right)
\left\langle \xi\right\rangle ^{2s}\right| \left\langle \varepsilon
q\right\rangle ^{s}dq$ are uniformly bounded for $0<\varepsilon\leq1$.
Consequently,

\begin{corollary}
\label{Corollary bdd}Let $\left( \varepsilon _{k}\right) $ and $\left(
u_{k}\right) $ satisfy the hypotheses of Proposition \ref{Prop definition wm
l2loc}. Then the sequence $\left( m^{\varepsilon _{k}}\left[ u_{k}\right]
\right) $ is bounded in $\mathcal{S}^{\prime }$.
\end{corollary}

Estimate (\ref{estimate meu on a Hs}) gives immediately the following:

\begin{remark}
\label{Rmk larger test functions}Lemma \ref{Lemma meuh bdd Hs} shows that $%
m^{\varepsilon }\left[ u\right] $ acts continuously on test functions $a$ in
the closure of $\mathcal{S}$ for the norm:
\begin{equation}
\left[ a\right] _{s}:=\int_{\mathbb{R}^{d}}\sup_{\xi \in \mathbb{R}%
^{d}}\left| \widehat{a}\left( q,\xi \right) \left\langle \xi \right\rangle
^{2s}\right| \left\langle q\right\rangle ^{s}dq<\infty .  \label{norma a}
\end{equation}
This closure contains the space
\begin{equation}
\Sigma ^{s}:=\left\{ \left\langle D_{x}\right\rangle ^{s}\left\langle \xi
\right\rangle ^{2s}a\in C_{0}\left( \mathbb{R}^{d}\times \mathbb{R}%
^{d}\right) :\left[ a\right] _{s}<\infty \right\} .
\label{definition sigmas}
\end{equation}
\end{remark}

\begin{remark}
\label{Rmk convergence larger test functions}Consequently, if $\left(
u_{k}\right) $ is as in Proposition \ref{Prop definition wm l2loc} and $%
\left( m^{\varepsilon _{k}}\left[ u_{k}\right] \right) $ converges weakly in
$\mathcal{S}^{\prime }$ then $\left\langle m^{\varepsilon _{k}}\left[ u%
\right] ,a\right\rangle $ converges as well for every $a\in \Sigma ^{s}$.
\end{remark}

\begin{proof}[Proof of Proposition \ref{Prop definition wm l2loc}]
The boundedness of the sequence $\left( m^{\varepsilon _{k}}\left[ u_{k}%
\right] \right) $ was proved in Corollary \ref{Corollary bdd}. Suppose now
that the distributions $m^{\varepsilon _{k}}\left[ u_{k}\right] $ weakly
converge to some $\mu \in \mathcal{S}^{\prime }$. We next show by means of a
localization argument that $\mu $ is a positive distribution and thus, due
to Schwartz's Theorem, a positive Radon measure.

Take $\phi \in \mathcal{S}\left( \mathbb{R}_{\xi }^{d}\right) $; Lemma \ref
{Lemma localization I} gives
\begin{equation*}
\lim_{k\rightarrow \infty }\left\langle m^{\varepsilon _{k}}\left[ \phi
\left( \varepsilon _{k}D_{x}\right) u_{k}\right] ,a\right\rangle _{\mathcal{S%
}^{\prime }\times \mathcal{S}}=\int_{\mathbb{R}^{d}\times \mathbb{R}%
^{d}}a\left( x,\xi \right) \left| \phi \left( \xi \right) \right| ^{2}d\mu
\left( x,\xi \right)
\end{equation*}
for every $a\in \mathcal{S}$. Since $\left( \phi \left( \varepsilon
_{k}D_{x}\right) u_{k}\right) $ is a bounded sequence in $L^{2}\left(
\mathbb{R}^{d}\right) $, Theorem \ref{Thm definition wm} ensures that $%
\left| \phi \left( \xi \right) \right| ^{2}\mu $ is a positive Radon measure
(and hence a positive distribution). But $\phi \in \mathcal{S}\left( \mathbb{%
R}_{\xi }^{d}\right) $ is arbitrary, so $\mu $ itself is positive and we
obtain the desired result.\medskip
\end{proof}

Notice that a very similar proof would give a version of Proposition \ref
{Prop definition wm l2loc} in the context of weighted spaces $L^{2}\left(
\mathbb{R}^{d};\left\langle x\right\rangle ^{r}\right) $.

\subsection{\label{sbs ploc}Proof of Proposition \ref{Prop localization}}

The key ingredient in the proof of the Proposition is the following
auxiliary result:

\begin{lemma}
\label{Lemma localization II}Under the assumptions of Proposition \ref{Prop
localization} and for every $a\in \mathcal{S}$, if any of the following
conditions hold:\medskip

i) $h_{k}=\varepsilon _{k}$ and $a$ vanishes on the set of discontinuity
points of $\phi $.\medskip

ii) $h_{k}\ll \varepsilon _{k}$ and $\phi $ is continuous at $\xi =0$%
.\medskip

Then
\begin{equation}
\lim_{k\rightarrow \infty }\left| \left\langle m^{\varepsilon _{k}}\left[
\phi \left( h_{k}D_{x}\right) u_{k}\right] -\left| \phi \left( \frac{h_{k}}{%
\varepsilon _{k}}\xi \right) \right| ^{2}m^{\varepsilon _{k}}\left[ u_{k}%
\right] ,a\right\rangle \right| =0.  \label{lim loc}
\end{equation}
\end{lemma}

\begin{proof}
Take $a\in \mathcal{S}$ and set $\Phi _{k}\left( \xi \right) :=\phi \left(
h_{k}/\varepsilon _{k}\xi \right) $. From relations (\ref{relation meuh and
meFuh}), (\ref{meu on a l2loc}.i) and (\ref{est exact}) we obtain:
\begin{equation*}
\left| \left\langle m^{\varepsilon _{k}}\left[ \phi \left( h_{k}D_{x}\right)
u_{k}\right] -\left| \phi \left( \frac{h_{k}}{\varepsilon _{k}}\xi \right)
\right| ^{2}m^{\varepsilon _{k}}\left[ u_{k}\right] ,a\right\rangle \right|
\leq M_{k}\left( a\right) \left\| \left\langle \varepsilon
_{k}D_{x}\right\rangle ^{-s}u_{k}\right\| _{L^{2}\left( \mathbb{R}%
^{d}\right) }^{2}
\end{equation*}
where
\begin{equation}
M_{k}\left( a\right) :=\int_{\mathbb{R}^{d}}\sup_{\xi \in \mathbb{R}%
^{d}}\left| \widehat{a}\left( q,\xi \right) \Phi _{k}\left( \xi \right) %
\left[ \Phi _{k}\left( \xi +\varepsilon _{k}q\right) -\Phi _{k}\left( \xi
\right) \right] \left\langle \xi +\varepsilon _{k}q\right\rangle
^{s}\left\langle \xi \right\rangle ^{s}\right| \frac{dq}{\left( 2\pi \right)
^{d}};  \label{remainder}
\end{equation}
recall that $\widehat{a}\left( q,\xi \right) $ stands for the Fourier
transform of $a\left( x,\xi \right) $ in $x$.

We now must prove that $M_{k}\left( a\right) \rightarrow 0$ as $k\rightarrow
\infty $. This will be done by first checking that for test functions $a$
belonging to the smaller class:
\begin{equation*}
\widehat{\mathcal{D}}:=\left\{ a\in \mathcal{S}:\widehat{a}\in C_{c}^{\infty
}\left( \mathbb{R}^{d}\times \mathbb{R}^{d}\right) \right\} .
\end{equation*}
Take $R>0$ such that $\limfunc{supp}a$ is contained in $B\left( 0;R\right)
\times B\left( 0;R\right) $.

When $k\in \mathbb{N}$ is sufficiently large, $\varepsilon _{k}\leq 1$ and
\begin{equation}
\frac{h_{k}}{\varepsilon _{k}}\left( \xi +\varepsilon _{k}q\right) \in
B\left( 0;2R\sup h_{k}/\varepsilon _{k}\right) \qquad \text{for every }q,\xi
\in B\left( 0;R\right) .  \label{xq in support}
\end{equation}
Suppose now that i) holds. If $C_{\phi }$ denotes the set of points where $%
\phi $ is continuous, then $\Phi _{k}=\phi $ is uniformly continuous over $%
C_{\phi }\cap B\left( 0;R\right) $ and, consequently,
\begin{equation*}
\sup_{q,\xi \in B\left( 0;R\right) }\mathbf{1}_{C_{\phi }}\left( \xi \right)
\left| \phi \left( \xi +\varepsilon _{k}q\right) -\phi \left( \xi \right)
\right| \rightarrow 0\qquad \text{as }k\rightarrow \infty
\end{equation*}
because of (\ref{xq in support}).

On the other hand, when $h_{k}/\varepsilon _{k}\rightarrow 0$ and $\phi $ is
continuous at $\xi =0$, again as a consequence of (\ref{xq in support}),
\begin{equation*}
\sup_{\xi ,q\in B\left( 0;R\right) }\left| \phi \left( \frac{h_{k}}{%
\varepsilon _{k}}\left( \xi +\varepsilon _{k}q\right) \right) -\phi \left(
\frac{h_{k}}{\varepsilon _{k}}\xi x\right) \right| \leq 2\sup_{\xi \in
B\left( 0;2h_{k}/\varepsilon _{k}R\right) }\left| \phi \left( \xi \right)
-\phi \left( 0\right) \right| \rightarrow 0\qquad \text{as }k\rightarrow
\infty .
\end{equation*}
Thus, in either case,
\begin{equation*}
\sup_{\xi \in \mathbb{R}^{d}}\left| \widehat{a}\left( q,\xi \right) \Phi
_{k}\left( \xi \right) \left[ \Phi _{k}\left( \xi +\varepsilon _{k}q\right)
-\Phi _{k}\left( \xi \right) \right] \left\langle \xi +\varepsilon
_{k}q\right\rangle ^{s}\left\langle \xi \right\rangle ^{s}\right|
\rightarrow 0\qquad \text{as }k\rightarrow \infty ,
\end{equation*}
for every $q\in \mathbb{R}^{d}$. Lebesgue's dominated convergence Theorem
gives the convergence to zero of the integrals (\ref{remainder}). The
density of $\widehat{\mathcal{D}}$ in $\mathcal{S}$ concludes the proof of
the Lemma.\bigskip
\end{proof}

\begin{proof}[Proof of Proposition \ref{Prop localization}]
To prove i) and ii) it only needs to be checked that, for any $a\in
C_{c}^{\infty }\left( \mathbb{R}^{d}\times \mathbb{R}^{d}\right) $ (if $%
\varepsilon _{k}=h_{k}$, we further require that $a|_{\mathbb{R}^{d}\times
\overline{D_{\phi }}}\equiv 0$), the functions $\left| \phi \left(
h_{k}/\varepsilon _{k}\xi \right) \right| ^{2}a\left( x,\xi \right) $ belong
to the class $\Sigma ^{s}$. If so, then
\begin{equation*}
\lim_{k\rightarrow \infty }\left\langle \left| \phi \left( \frac{h_{k}}{%
\varepsilon _{k}}\xi \right) \right| ^{2}m^{\varepsilon _{k}}\left[ u_{k}%
\right] ,a\right\rangle _{\mathcal{S}^{\prime }\times \mathcal{S}}=\int_{%
\mathbb{R}^{d}\times \mathbb{R}^{d}}\left| \phi \left( c\xi \right) \right|
^{2}a\left( x,\xi \right) d\mu ,
\end{equation*}
holds with $c:=\lim h_{k}/\varepsilon _{k}$, because of Remark \ref{Rmk
convergence larger test functions}. The conclusion would then follow from
identity (\ref{lim loc}).

First, notice that $\left| \phi \left( h_{k}/\varepsilon _{k}\cdot \right)
\right| ^{2}a$ are compactly supported and infinitely differentiable in $x$.
When $\varepsilon _{k}=h_{k}$ we must verify that $\left| \phi \right|
^{2}a\in \Sigma ^{s}$ which is clearly the case if $a|_{\mathbb{R}^{d}\times
\overline{D_{\phi }}}\equiv 0$, for then $\left| \phi \right| ^{2}a$ is
continuous in $\xi $.

On the other hand, if $h_{k}\ll \varepsilon _{k}$ and $\phi $ is merely
continuous in a ball $B\left( 0;\delta \right) $ then, for $k$ large enough,
$\limfunc{supp}a\subset B\left( 0;h_{k}/\varepsilon _{k}\delta \right) $ and
consequently, $\phi \left( h_{k}/\varepsilon _{k}\cdot \right) $ is
continuous on $\limfunc{supp}a$.\medskip
\end{proof}

\subsection{\label{sbs propm}Additional properties.}

The next approximation result is sometimes useful in the computation of
Wigner measures:

\begin{lemma}
\label{Lemma approx WM}Let $\left( u_{k}\right) $ and $\left(
u_{k}^{N}\right) $ be sequences in $H^{-s}\left( \mathbb{R}^{d}\right) $, $%
s\geq 0$, satisfying (\ref{uk is Hsh bounded}) with the same bound and
\begin{equation*}
\limsup_{k\rightarrow \infty }\left\| \left\langle \varepsilon
_{k}D_{x}\right\rangle ^{-s}\left( u_{k}-u_{k}^{N}\right) \right\|
_{L^{2}\left( \mathbb{R}^{d}\right) }\rightarrow 0\text{\qquad as }%
N\rightarrow \infty \text{.}
\end{equation*}
Suppose that $m^{\varepsilon _{k}}\left[ u_{k}\right] $ and $m^{\varepsilon
_{k}}\left[ u_{k}^{N}\right] $ converge respectively to $\mu $ and $\mu _{N}$%
. Then
\begin{equation*}
\mu _{N}\rightharpoonup \mu \text{\qquad in }\mathcal{M}_{+}\left( \mathbb{R}%
^{d}\times \mathbb{R}^{d}\right) \text{ as }N\rightarrow \infty .
\end{equation*}
\end{lemma}

\begin{proof}
This is a simple consequence of the identity:
\begin{align*}
\left\langle m^{\varepsilon _{k}}\left[ u_{k}\right] -m^{\varepsilon _{k}}%
\left[ u_{k}^{N}\right] ,a\right\rangle _{\mathcal{S}^{\prime }\times
\mathcal{S}}& =\left\langle \overline{u_{k}^{N}},a\left( x,\varepsilon
_{k}D_{x}\right) \left( u_{k}-u_{k}^{N}\right) \right\rangle _{\mathcal{S}%
^{\prime }\times \mathcal{S}}+ \\
& +\left\langle \overline{\left( u_{k}-u_{k}^{N}\right) },a\left(
x,\varepsilon _{k}D_{x}\right) u_{k}\right\rangle _{\mathcal{S}^{\prime
}\times \mathcal{S}}.
\end{align*}
This gives an estimate:
\begin{equation*}
\left| \left\langle m^{\varepsilon _{k}}\left[ u_{k}\right] -m^{\varepsilon
_{k}}\left[ u_{k}^{N}\right] ,a\right\rangle _{\mathcal{S}^{\prime }\times
\mathcal{S}}\right| \leq C\left\| \left\langle \varepsilon
_{k}D_{x}\right\rangle ^{-s}\left( u_{k}-u_{k}^{N}\right) \right\|
_{L^{2}\left( \mathbb{R}^{d}\right) };
\end{equation*}
taking limits as $k\rightarrow \infty $ we obtain:
\begin{equation*}
\left| \int_{\mathbb{R}^{d}\times \mathbb{R}^{d}}a\left( x,\xi \right)
\left( d\mu -d\mu _{N}\right) \right| \leq C\limsup_{k\rightarrow \infty
}\left\| \left\langle \varepsilon _{k}D_{x}\right\rangle ^{s}\left(
u_{k}-u_{k}^{N}\right) \right\| _{L^{2}\left( \mathbb{R}^{d}\right) }
\end{equation*}
and the result follows, since the measures $\mu _{N}$ and $\mu $ are
equibounded.\medskip
\end{proof}

We conclude this section with an almost orthogonality result:

\begin{lemma}
\label{Lemma orth}Let $\left( u_{k}\right) $ and $\left( v_{k}\right) $ be
sequences in $H^{-s}\left( \mathbb{R}^{d}\right) $, $s\geq 0$, satisfying (%
\ref{uk is Hsh bounded}) for some scale $\left( \varepsilon _{k}\right) $.
Suppose their Wigner measures at scale $\left( \varepsilon _{k}\right) $, $%
\mu $ and $\nu $ are mutually singular. Then $m^{\varepsilon _{k}}\left[
u_{k}+v_{k}\right] $ converges to $\mu +\nu $.
\end{lemma}

\begin{proof}
A proof of this result for $s=0$ may be found in \cite{Ge91c} or \cite
{Li-Pau}. For the general case, it suffices to take into account Remark \ref
{Rmk finite} to conclude that the Wigner measures of $\left\langle
\varepsilon _{k}D_{x}\right\rangle ^{-s}u_{k}$ and $\left\langle \varepsilon
_{k}D_{x}\right\rangle ^{-s}v_{k}$ are $\left\langle \xi \right\rangle
^{-2s}\mu $ and $\left\langle \xi \right\rangle ^{-2s}\nu $. These are
clearly mutually singular and thus the aforementioned $L^{2}$-version of the
present result gives
\begin{equation*}
m^{\varepsilon _{k}}\left[ \left\langle \varepsilon _{k}D_{x}\right\rangle
^{-s}(u_{k}+v_{k})\right] \rightharpoonup \left\langle \xi \right\rangle
^{-2s}\mu +\left\langle \xi \right\rangle ^{-2s}\nu
\end{equation*}
and finally
\begin{equation*}
m^{\varepsilon _{k}}\left[ u_{k}+v_{k}\right] \rightharpoonup \left\langle
\xi \right\rangle ^{2s}\left( \left\langle \xi \right\rangle ^{-2s}\mu
+\left\langle \xi \right\rangle ^{-2s}\nu \right) =\mu +\nu
\end{equation*}
as claimed.
\end{proof}

\bigskip

\textbf{Acknowledgments: }This article extends and improves some of the
results of the author's Ph.D. Thesis. He wishes to thank the guidance of his
Ph.D. advisor Enrique Zuazua. He also wishes to thank Patrick G\'{e}rard,
for many helpful discussions and suggestions.

This work has been supported by projects BFM02-03345 of MCyT (Spain) and
HYKE (ref. HPRN-CT-2002-00282), HMS2000 (ref. HPRN-CT-2000-00109) of the
European Union.

\end{document}